\documentclass[times]{elsarticle}
\usepackage[margin=1in]{geometry}
\usepackage[utf8]{inputenc}
\usepackage{amsmath}
\usepackage{xcolor}
\usepackage{graphicx}
\usepackage{url}
\usepackage{hyperref}
\usepackage[ruled,vlined]{algorithm2e}
\usepackage{comment}
\usepackage{placeins}
\usepackage{float}

\graphicspath{ {./Plots/} }

\makeatletter
\def\ps@pprintTitle{%
  \let\@oddhead\@empty
  \let\@evenhead\@empty
  \let\@oddfoot\@empty
  \let\@evenfoot\@oddfoot
}
\makeatother

\begin{document}

\begin{frontmatter}

\title{A deep learning approach for adaptive zoning}

\author[rvt]{Massimiliano Lupo Pasini\corref{cor1}} \ead{lupopasinim@ornl.gov}
\author[rvt2]{ Luka Maleni\c{c}a}

\author[rvt]{Kwitae Chong}

\author[rvt]{Stuart Slattery}

\address[rvt]{Oak Ridge National Laboratory, Computational Sciences and Engineering Division, 1 Bethel Valley Road, Oak Ridge, TN, USA, 37831}

\address[rvt2]{ETH Zurich, Institute for Building Materials, Zurich, 8093, Switzerland}

\cortext[cor1]{Corresponding author}

\date{}

\begin{abstract}
    We propose a supervised deep learning (DL) approach to perform adaptive zoning on time dependent partial differential equations that model the propagation of 1D shock waves in a compressible medium. We train a neural network on a dataset composed of different static shock profiles associated with the corresponding adapted meshes computed with standard adaptive zoning techniques. We show that the trained DL model learns how to capture the presence of shocks in the domain and generates at each time step an adapted non-uniform mesh that relocates the grid nodes to improve the accuracy of Lax-Wendroff and fifth order weighted essentially non-oscillatory (WENO5) space discretization schemes. We also show that the surrogate DL model reduces the computational time to perform adaptive zoning by at least a 2x factor with respect to standard techniques without compromising the accuracy of the reconstruction of the physical quantities of interest. 
\end{abstract}

\begin{keyword}
Artificial Intelligence,
Deep Learning,
Numerical Partial Differential Equations,
Compressible Flow, 
Mesh Adaptivity,
Adaptive Zoning
\end{keyword}

\end{frontmatter}

{\footnotesize \noindent This manuscript has been authored in part by UT-Battelle, LLC, under contract DE-AC05-00OR22725 with the US Department of Energy (DOE). The US government retains and the publisher, by accepting the article for publication, acknowledges that the US government retains a nonexclusive, paid-up, irrevocable, worldwide license to publish or reproduce the published form of this manuscript, or allow others to do so, for US government purposes. DOE will provide public access to these results of federally sponsored research in accordance with the DOE Public Access Plan (\url{http://energy.gov/downloads/doe-public-access-plan}).}

\section*{Introduction}
Mesh generation lies at the core of many areas in advanced computational sciences and engineering. To both improve computational performance as well as automate some aspects of the mesh generation process, adaptive mesh optimization and refinement play an important role in complex high-fidelity simulation fields. The accuracy and convergence of solutions using mesh-based numerical methods are strongly dependent on the quality of the mesh being used and generation of an appropriate quality mesh to support these calculations is a challenging problem.

Over the years, adaptive zoning (also called adaptive moving mesh method or $r$-adaptivity) has been used to generate meshes that can help the numerical schemes to attain higher accuracy by producing a higher concentration of mesh points in regions of the domain where a higher resolution is desired \cite{brackbill82, brackbill93,  ceniceros_efficient_2001, lapenta2003, budd_adaptivity_2009, huang2010, huang_variational_2011, delzanno2011, browne2014, remski_balanced_2014, hu2015, huang_comparative_2015, budd2015, pathak2016, heaton_ian_2018, fritzen_--fly_2019}.
An efficient adaptive zoning is relevant for applications that need to model physical quantities characterized by steep gradients that propagate in space over time such as shock waves in a compressible medium. An accurate and numerically stable reconstruction of shocks requires a higher concentration of grid nodes in regions with high steep gradients, whereas fewer nodes can be used in regions where the profile has a smoother shape.

The approaches proposed in the literature to perform adaptive zoning rely on the equidistribution principle to ensure that a quantity of interest is equally distributed across all the cells. Relocating the grid nodes according to the equidistribution principle on one-dimensional problems can be performed using Bohr's algorithm that is inexpensive to apply \cite{yang2012}, but cannot be generalized to multiple dimensions. Performing adaptive zoning strategies according to the equidistribution principle in multiple dimensions requires solving an auxiliary non-linear partial differential equation (PDE) with homogeneous Dirichlet boundary conditions, also denoted in the literature as moving mesh PDE (MMPDE) \cite{brackbill82, brackbill93}. Since the size of the grid used to solve the auxiliary MMPDE is as large as the size of the grid used to solve the original physics problem, performing adaptive zoning with MMPDE in multiple dimensions becomes computationally expensive for large scale simulations. Existing PDE solvers combined with adaptive zoning either (i) perform computations on the physical domain by using flexible discretization schemes that can handle non-uniform meshes \cite{CrnjaricZic2007EfficientIO, CAPDEVILLE20082977, CORALIC201495, huang2018} or (ii) remap the equations on a computational domain where the grid is always uniform and keep track of the relocation of the nodes on the physical mesh by adjusting the coefficients of the PDEs with Jacobian terms that describe the local deformation of the physical domain \cite{ cao_study_1999, cao_approaches_2003, chacon2006, lapenta2006, chacon2011}.   
Performing computations directly on the physical domain has the advantage of solving simpler PDEs, but can only handle non-uniform meshes without curvatures. On the contrary, performing computations on the computational domain requires solving more complex PDEs with additional terms that keep track of the remapping, but the advantage is an enhanced flexibility of the procedure that can handle adapted physical meshes with curvatures. 

For situations where adaptive zoning is performed on the physical domain for one-dimensional PDE solvers, we present a novel deep learning (DL) approach to perform fast and accurate adaptive zoning as a computationally convenient alternative to solving MMPDEs. 
Our approach starts by training a neural network (NN) on a dataset consisting of different staircase shock profiles, where the number of jumps and the location of each jump changes arbitrarily from case to case. Each staircase shock profile is associated with an adapted mesh computed by solving the MMPDE. The DL model learns how to relocate the grid nodes based on the presence of an arbitrary number of regions where the quantity of interest is characterized by steep gradients in the domain. Once the predictive performance of the NN is validated, the NN is deployed as a surrogate DL model to perform adaptive zoning. Since the use of a trained NN avoids solving expensive auxiliary PDEs, our approach provides a path to efficiently perform adaptive zoning in situations with limited computational budget, which is an important aspect in multi-physics modeling for two reasons: the available computational power and the period of time when the computational power is available. The former imposes obvious intrinsic limitations, the latter becomes important when critical decisions have to be made in a timely and accurate manner. 
Since the training dataset contains only staircase profiles where the gradients are either zero or infinite, the DL models learns to identify regions with infinite gradients and concentrate more grid nodes around those regions and it disregards regions where there are non-zero gradients associated with a smoother profile, where a low concentration of nodes is sufficient for the numerical scheme to attain high accuracy. 
Numerical results performed on 1D Euler equations for the Sod shock tube test \cite{sod1978}, the Taylor–von Neumann–Sedov blast wave test \cite{sedov1946}, and the Woodward-Colella blast waves test case \cite{woodward} show that the DL model outperforms the standard adaptive zoning in efficiently capturing the presence of shocks and better reproduces the discontinuities of the physical quantities by minimizing the numerical oscillation of the discretization scheme close to the shock. Moreover, the DL approach is faster than solving the MMPDE and reduces the computational time to perform adaptive zoning by at least a 2x factor. 

This work is structured as follows. Section \ref{background} describes the 1D Euler equations used to model the propagation of waves in an ideal polytropic gas and the standard approach to perform adaptive zoning via solving a global elliptic problem to calculate the node remapping. Section \ref{deep_learning_section} describes the architecture of the DL model and the procedure to generate the training and validation data. In Section \ref{numerical_section}, a first set of numerical results in 1D are presented that assess the accuracy of the DL model with respect to different volumes of training data, showing an improvement in accuracy when the training data increases in size. The second set of numerical results describes the performance of the DL approach for adaptive zoning on test cases described by the 1D Euler equations.
Section \ref{conclusion_section} contains concluding remarks and future developments. 

\section{Background}
\label{background}
\subsection{1D Euler equations}
The physics of compressible flow systems addressed in this work is described by the 1D time-dependent Euler equations \cite{toro99} that result from neglecting the effect of viscosity, heat conduction, and body forces on a compressible medium. These equations model the behavior of an ideal gas and describe the evolution of rarefactions, contacts, and shock waves.  
The conservative formulation of the 1D Euler equations reads:
\begin{equation}
\begin{cases}
\frac{\partial \rho (x,t)}{\partial t} + \frac{\partial}{\partial x}(\rho(x,t) u(x,t)) = 0 \;\;\; \qquad \qquad \qquad \qquad \qquad\text{in} \quad (a,b)\times(0,T) \\
\frac{\partial }{\partial t} [\rho(x,t) u(x,t)] + \frac{\partial}{\partial x}[\rho(x,t) u^2(x,t) + p(x,t)] = 0 \qquad \; \text{in} \quad (a,b)\times(0,T) \\
\frac{\partial E(x,t)}{\partial t} + \frac{\partial}{\partial x}[(E(x,t) + p(x,t))u(x,t)] = 0  \qquad \qquad \qquad \text{in} \quad (a,b)\times(0,T) \\
\rho(x,0) = \rho_0, \quad u(x,0) = u_0, \quad p(x,0) = p_0  \qquad \qquad \qquad \text{for} \quad x\in [a,b]\\
\rho(a,t) = \rho_a, \quad u(a,t) = u_a, \;\quad p(a,t) = p_a \qquad \qquad \quad \quad \text{for} \quad t\in [0,T)\\
\rho(b,t) = \rho_b, \quad u(b,t) = u_b, \qquad p(b,t) = p_b \qquad \qquad \qquad \text{for} \quad t\in [0,T)\\
\end{cases}
\label{euler}
\end{equation}
where $\rho$, $u$, $E$, and $p$ are the gas density, velocity, total energy, and pressure, respectively. The total energy has two components: the internal energy density, $\rho e$, and the kinetic energy density, $\rho u^2$: 
\begin{equation}
E = \rho e + \rho u^2.
\end{equation}
Since there are three conservation laws in Equations \eqref{euler} and four unknowns (density, momentum, energy, and pressure), we need to introduce an equation of state to close the system. For a polytropic ideal gas, the state equation relates the specific internal energy $e$ to the pressure and density as follows:
\begin{equation}
    p = \rho e(\gamma - 1),
\end{equation}
where $\gamma$ is the polytropic constant, i.e. the ratio of the specific heat at constant pressure, $c_p$, and the specific heat at constant volume, $c_v$, $\displaystyle \gamma = \frac{c_p}{c_v}$.

\subsection{Numerical discretization of 1D scalar hyperbolic PDEs}
All three scalar equations in the PDE system \eqref{euler} are 1D scalar hyperbolic PDEs that can be generalized in the following form:
\begin{equation}
    \frac{\partial u (x,t)}{\partial t} + \frac{\partial f (u(x,t))}{\partial x} = 0,
    \label{general_pde}
\end{equation}
where $u$ can either represent the density $\rho$, the momentum $\rho u$, or the total energy $E$. 
We discretize the 1D space domain $[a,b]$ into a sequence of cells $[x_{i-\frac{1}{2}}, x_{i+\frac{1}{2}}]$ with $i = 1, \ldots, n$ such that
\begin{equation}
    a = x_\frac{1}{2} < x_1 \ldots < x_{n-1} < x_{n+\frac{1}{2}} = b.
\end{equation}
and we discretize the time interval $[0,T]$ with subintervals of length $\Delta t$ such that $t^i = i \Delta t$.
The average value of $f(u(x,t))$ in cell $i$, at a fixed time $t = t^m = m\Delta t$ is defined as 
\begin{equation}
    f_i^m = \frac{1}{[x_{i+\frac{1}{2}} - x_{i-\frac{1}{2}}]} \int_{x_{i-\frac{1}{2}}}^{x_{i+\frac{1}{2}}}f(u(x,t^m))dx
\end{equation}
 We shall assign the value $f_i^m$ at the centre of the cell, which gives rise to the \textit{conservative cell-centered methods} also known as \textit{finite volume (FV) methods}. 
The specifics of how the discretization values in space and time are combined impacts the overall accuracy of the numerical scheme. 

The Lax-Wendroff scheme can be derived in different ways. The one we use for the numerical discretization of PDEs in this work is a two-step procedure. The first step in the Richtmyer Lax–Wendroff method calculates values for $u(x,t)$ at half time steps, $t^{m+\frac{1}{2}}$ and at the cell center $x_{i+\frac{1}{2}}$ as follows:
\begin{equation}
    u_{i+\frac{1}{2}}^{m+\frac{1}{2}} = \frac{1}{2}(u_{i+1}^m + u_{i}^m) - \frac{\Delta t}{2[x_{i+1}-x_{i}]} (f(u_{i+1}^m)-f(u_{i}^m))
\end{equation}
\begin{equation}
    u_{i-\frac{1}{2}}^{m+\frac{1}{2}} = \frac{1}{2}(u_{i}^m + u_{i-1}^m) - \frac{\Delta t}{2[x_{i}-x_{i-1}]} (f(u_{i}^m)-f(u_{i-1}^m)).
\end{equation}
The second step calculates the value for $u(x,t)$ at $t^m$ using the data for $t^m$ and $t^{m+\frac{1}{2}}$
\begin{equation}
    u_i^{m+1} = u_i^m - \frac{\Delta t}{\frac{1}{2}[x_{i+1} - x_{i}] + \frac{1}{2}[x_{i} - x_{i-1}]}\bigg( f(u_{i+\frac{1}{2}}^{m+\frac{1}{2}}) - f(u_{i-\frac{1}{2}}^{m+\frac{1}{2}}) \bigg).
\end{equation}
The Lax-Wendroff scheme combines first-order time derivatives and first-order space derivatives in such a way that the resulting overall scheme is second-order accurate both in space and time. For a comprehensive treatment of the family of Lax-Wendroff schemes see \cite{hirsch}. 

\begin{figure}
    \centering
    \includegraphics[width=0.5\textwidth]{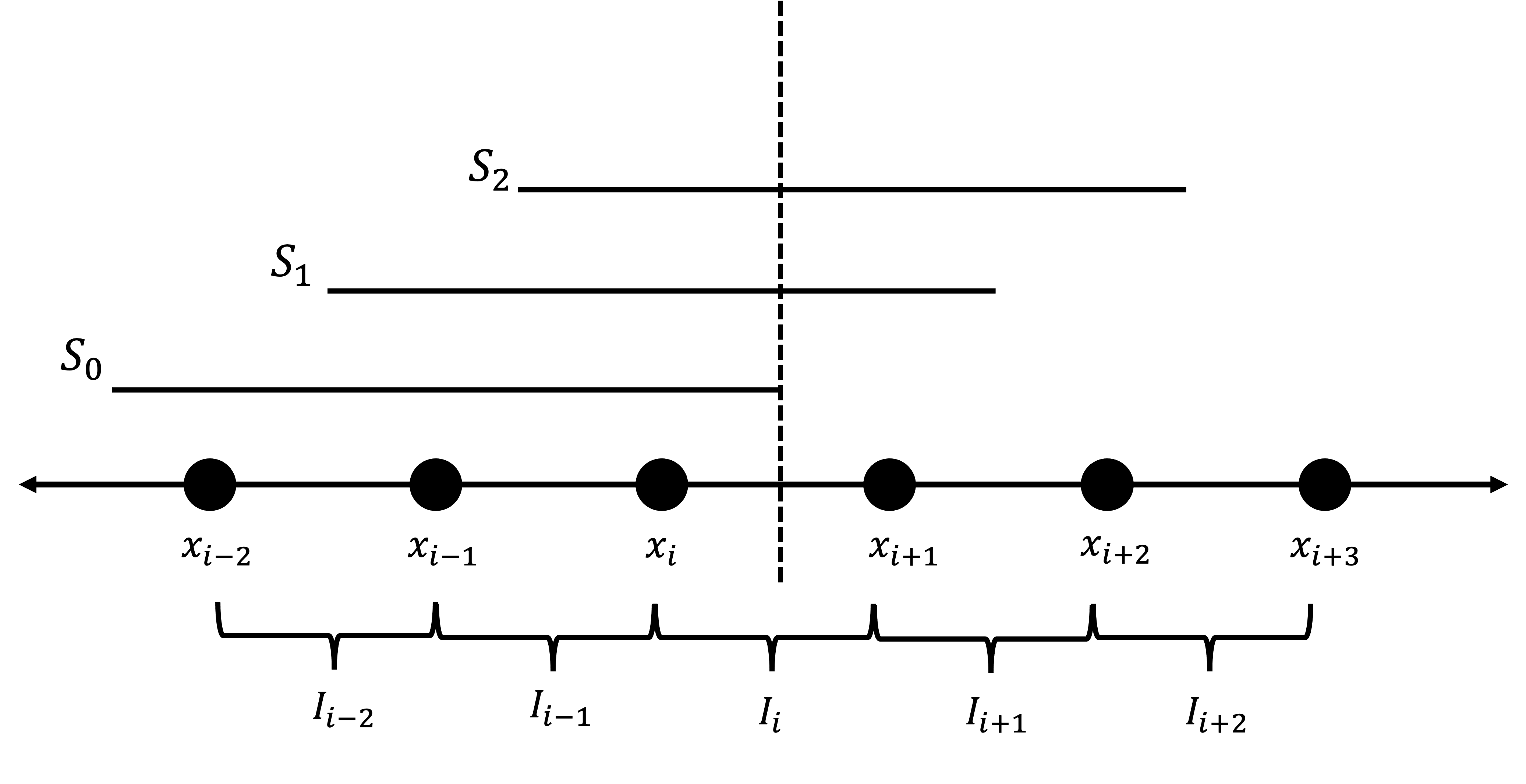}
    \caption{Grid stencil used by WENO5. }
    \label{stencil}    
\end{figure}

To improve the accuracy order of the numerical scheme, one can adopt different discretization techniques for space and time. 
For space discretization, one option is the use of weighted essentially non-oscillatory (WENO) methods that represent a class of high-order schemes. We used a variant described in \cite{huang2018} that handles non-uniform meshes and stems from the original implementation proposed in \cite{shu}.
In the fifth-order WENO method (WENO5), there are three candidate stencils, as shown in Figure \ref{stencil}. The $i$-th cell is denoted as $I_i$. The three candidate stencils are, respectively, $S_0=\{I_{i-2},I_{i-1},I_i\}$, $S_1=\{I_{i-1},I_i,I_{i+1}\}$, and $S_2=\{I_i,I_{i+1},I_{i+2}\}$.
We define the average value of $u(x,t)$ on the cell $I_i$ at time $t^m$ as
\begin{equation}
    \overline{u}_i^m = \frac{1}{x_{i+\frac{1}{2}} - x_{i-\frac{1}{2}}} \int_{x_{i-\frac{1}{2}}}^{x_{i+\frac{1}{2}}} u(x, t^m)dx
\end{equation}
and
\begin{equation}
    \tilde{u}_{i+\frac{1}{2}}^m = (1-\alpha)\overline{u}_i^m + \alpha \overline{u}_{i+1}^m, \quad \alpha = \frac{x_{i+\frac{1}{2}} - x_{i-\frac{1}{2}}}{ ( x_{i+\frac{1}{2}} - x_{i-\frac{1}{2}} ) + ( x_{i+\frac{3}{2}} - x_{i+\frac{1}{2}} ) }.
\end{equation}
On the stencils $S_0$, $S_1$, $S_2$ we calculate the interfacial states $q^{0,m}_{i+\frac{1}{2}}$ , $q^{1,m}_{i+\frac{1}{2}}$, $q^{2,m}_{i+\frac{1}{2}}$ and the smoothness indicators $IS^{0,m}_{i+\frac{1}{2}}$, $IS^{1,m}_{i+\frac{1}{2}}$, $IS^{2,m}_{i+\frac{1}{2}}$ at time $t=t^m$. For specifics about the definition of these quantities, we refer to \cite{huang2018}.
The final $5$-th order polynomial reconstructed interfacial state of the quantity $u(x_{i+\frac{1}{2}},t^m)$ is
\begin{equation}
    \hat{u}_{i+\frac{1}{2}}^m = \sum_{k=0}^2\Omega_k^m q_{i+\frac{1}{2}}^k,
\end{equation}
where 
\begin{equation}
\Omega_k^m = \frac{\alpha_k^m}{\alpha_0^m + \alpha_1^m + \alpha_2^m}, \quad \alpha_k^m =\frac{C_k}{{(\epsilon + IS_k^m)}^p}, \quad
    C_k = \begin{cases}
    \frac{1}{10}\quad k=0\\
    \frac{6}{10}\quad k=1\\
    \frac{3}{10}\quad k=2
    \end{cases} 
    \quad k=0,1,2
\end{equation}
and $\epsilon = 1\mathrm{e}-6$, $p=2$.
For high order time stepping, we used third-order RK also abbreviated as RK-3 as described in \cite{runge, kutta}.

To increase the spatial resolution in regions that are identified with large solution gradients, we allow the nodes of the computational mesh to move in time while keeping the domain boundaries fixed. The relocation of the nodes of the grid is performed using adaptive zoning which we describe next. 

\subsection{1D r-adaptive zoning by solving an auxiliary MMPDE}
\label{adaptive_zoning_section}
We define a \textit{monitor function} $\omega(x)$ \cite{huang_variational_2011, pathak2016, cao_study_1999, cao_approaches_2003} as
\begin{equation}
\begin{split}
\omega(x) = \sqrt{1+ \alpha_1 \bigg \lvert \frac{\frac{d\rho(x, t)}{dx}}{\lVert \frac{d\rho(x, t)}{dx} \rVert_{L^{\infty}} }\bigg \rvert^2+\alpha_2 \bigg \lvert \frac{\frac{d(\rho(x, t) u(x, t))}{dx}}{\lVert \frac{d(\rho(x, t) u(x, t))}{dx} \rVert_{L^{\infty}}} \bigg \rvert^2 + \alpha_3 \bigg \lvert \frac{\frac{d e(x, t)}{dx}}{\lVert \frac{d e(x, t)}{dx}\rVert_{L^{\infty}} } \bigg \rvert^2 + \alpha_4 \bigg \lvert \frac{\frac{dp(x, t)}{dx}}{\lVert \frac{dp(x, t)}{dx}\rVert_{L^{\infty}}} \bigg \rvert^2 }, \\ \quad t\in [0,T]  
\end{split}
\label{monitor_function}
\end{equation}
where $\alpha_1$, $\alpha_2$, $\alpha_3$, and $\alpha_4$ are adjustable positive constants. The goal of the adaptive zoning is to compute a bijective mapping $\hat{x}(x)$ to transform the original coordinates $x$ of the mesh into a new set of coordinates $\hat{x}$ that equidistributes $\omega(\hat{x})$ over each cell of the mesh. 
This is achieved by solving the following auxiliary non-linear elliptic PDE problem called moving mesh PDE (MMPDE) \cite{brackbill82, cao_study_1999, cao_approaches_2003, brackbill82, chacon2006, lapenta2006, chacon2011}
\begin{equation}
\begin{cases}
    \frac{d}{d x}\bigg [\omega(\hat{x})\frac{d \hat{x}}{d x}\bigg ] = 0\\
    \hat{x}(a) = a\\
    \hat{x}(b) = b
\end{cases}.
\label{elliptic_bvp}
\end{equation}
Since the monitor function depends on the solution of the physical PDE itself, the PDE problem in \eqref{elliptic_bvp} is non-linear and it is solved with a non-linear fixed-point iteration scheme where monitor function and adapted mesh are alternatively updated until self-consistency between the two quantities is reached. The fixed-point scheme stops when the discrepancy between two consecutive updates drops below a tolerance defined by the user.
If the quantities modeled by the physical PDE abruptly change in a short time, the adaptive zoning performed by solving Equation \eqref{elliptic_bvp} may cause the grid nodes to change their location too aggressively and this may affect the accuracy of the PDE solver. 
One way to mitigate the change of the grid across consecutive time steps is to replace the elliptic MMPDE Equation \eqref{elliptic_bvp} with a parabolic MMPDE
\begin{equation}
\begin{cases}
    \frac{d\hat{x}}{d\tau}-\frac{d}{d x}\bigg [\omega(\hat{x})\frac{d \hat{x}}{d x}\bigg ] = 0\\
    \hat{x}(a) = a\\
    \hat{x}(b) = b 
\end{cases},
\label{parabolic_bvp}
\end{equation}
where the pseudo time $\tau$ determines the timescale over which the solution of Equation \eqref{parabolic_bvp} relaxes to a steady state.
After performing adaptive zoning with the parabolic MMPDE in Equation \eqref{parabolic_bvp}, the grid nodes occupy an intermediate location between the initial location and the one occupied if adaptive zoning were performed by solving the elliptic MMPDE in Equation \eqref{elliptic_bvp}. When Equation \eqref{parabolic_bvp} is used to perform adaptive zoning, the time stepping for the auxiliary PDE is always chosen the be smaller than the time stepping of the physics PDE solver, and the time horizon reached in the auxiliary MMPDE is defined by the user. 

\begin{figure}
    \centering
    \includegraphics[width=0.8\textwidth]{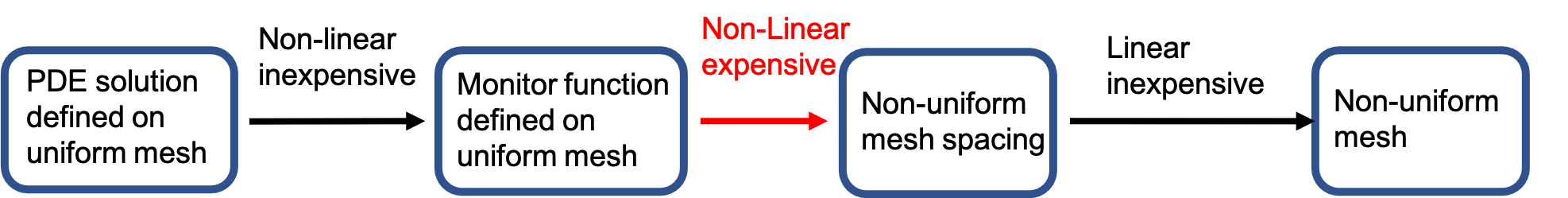}
    \caption{Sequence of computational tasks performed to map the solution to the PDE on a uniform mesh with the final adapted mesh. Solving the auxiliary nonlinear MMPDE is the dominant term in the total computational time.}
    \label{pipeline}    
\end{figure}

Steep gradients in the profile $u$ lead to high values in the monitor function, which in turn cause a high concentration of grid nodes in the local region of the domain where the the monitor function $\omega(x)$ exhibits steep gradients. 
To avoid nodes from collapsing, one can either (i) tune the scalars $\alpha_1$, $\alpha_2$, $\alpha_3$, and $\alpha_4$ and/or (ii) apply a kernel smoothing (e.g Gaussian smoothing, moving average, weighted averaging) over the monitor function via a convolutional operation as follows
\begin{equation}
    \hat{\omega}(x) = \int_{a}^b K(x-y)\omega(y)dy
\end{equation}
where $K$ is a smoothing (regularizing) convolutional kernel operator such as Gaussian kernel or moving average \cite{Savitzky64, kecs82, Enders2004}.
The support of the smoothing kernel is tunable: the wider the support, the smoother (less aggressive) is the re-zoning of the grid nodes. 
 The auxiliary MMPDEs in Equations \eqref{elliptic_bvp} and \eqref{parabolic_bvp} are discretized using a first order FV method with the same mesh that discretizes the PDE system in Equation \eqref{euler}. The computational cost to perform adaptive zoning using an MMPDE thus increases with a higher resolution of the mesh to discretize the PDE system in Equation \eqref{euler}, which causes non-negligible computational overheads for large scale simulations. 
The schematic in Figure \ref{pipeline} summarizes the sequential steps performed by standard adaptive zoning techniques to map an initially uniform mesh into a non-uniform adapted mesh by solving the non-linear auxiliary MMPDE \eqref{elliptic_bvp} or \eqref{parabolic_bvp}.

\section{Deep learning for adaptive zoning}
\label{deep_learning_section}
Since solving Equation \eqref{elliptic_bvp} and \eqref{parabolic_bvp} for large scale simulations is be computationally expensive, we use a DL model as surrogate to generate an adapted non-uniform mesh and use the DL model to replace some of the steps in the procedure described in Figure \ref{pipeline}.

\subsection{Deep feedforward networks}
\label{neural_network}
A deep feedforward network, also called feedforward NN or multilayer perceptron (MLP) \cite{mlp, goodfellow}, is a predictive statistical model to approximate some unknown function $f$ of the form
\[
y=f(x),
\]
where $x\in\mathbb{R}^a$, $y\in \mathbb{R}^b$ and $f:\mathbb{R}^a\rightarrow \mathbb{R}^b$. 
Given a collection of $m$ vectors for $x$ stored in $\mathbf{x}\in \mathbb{R}^{am}$ and the corresponding $m$ vectors of $f(x)$ stored in $\mathbf{y}=f(\mathbf{x})\in \mathbb{R}^{bm}$, feedforward NNs search a mapping that best approximates the unknown function $f$. In this work, $\mathbf{y}$ can be either the coordinates of the adapted mesh or the non-uniform mesh spacing of the adapted mesh. The input quantities $\mathbf{x}$ can be either the shock profile or the monitor function.  

A feedforward NN provides an approximation of $f(x)$ as a composite function: 
\begin{equation}
\hat{f}(x) = f_\ell(f_{\ell-1}(\ldots f_0(x))),
\label{composition}
\end{equation}
where $\hat{f}:\mathbb{R}^a\rightarrow \mathbb{R}^b$, $f_0:\mathbb{R}^a\rightarrow \mathbb{R}^{k_1}$, $f_\ell:\mathbb{R}^{k_\ell}\rightarrow \mathbb{R}^{b}$ and $f_p:\mathbb{R}^{k_{p}}\rightarrow \mathbb{R}^{k_{p+1}}$ for $p=1,\ldots,\ell-1$.
The proper number $\ell$ will be identified so that the composition in Equation \eqref{composition} resembles the unknown function $f$. The number $\ell$ quantifies the complexity of the composition and is equal to the number of hidden layers in the NN. Therefore, $f_1$ corresponds to the first hidden layer of the NN, $f_2$ is the second hidden layer, and so on. {Each one of the functions $f_i$'s combines the regression coefficients between consecutive hidden layers through nonlinear activation functions that allow deep feedforward networks to create nonlinear relations between input $x$ and output $y$.} 

\begin{figure}[H]
    \centering
    \includegraphics[width=0.5\textwidth]{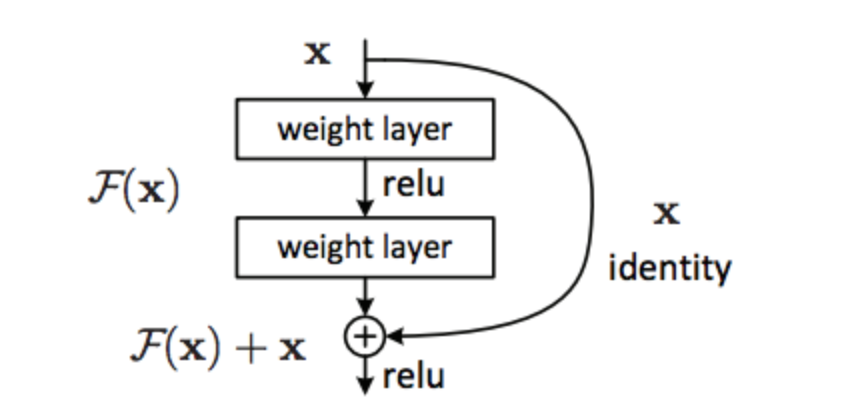}
    \caption{Structure of a residual block used in ResMLP. The residual block has two main components: (i) a corrective sequence of layers and (ii) a skip connection. The corrective sequence is made of two hidden layers that aim to improve the predictive performance of the DL model obtained by the previous set of hidden layers through an additive correction. The skip connection operates as an identity block by preserving the quantity $\mathbf{x}$ as is.}
    \label{residual_block}    
\end{figure}

The training of deep NNs is affected by numerical artifacts like vanishing gradients which limit the predictive capacity of the model. This problem has been addressed and solved with the introduction of residual neural networks (ResNet) \cite{resnet}, where successive hidden layers are used to learn features with incremental complexity in a numerically stable fashion, and this leads to an improved final accuracy for the model. Since the neurons of adjacent hidden layers are fully connected via a multilayer perceptron (MLP) structure, from now on we refer to the ResNets used in this work as ResMLP to distinghuish them from the most common version of ResNets that are characterized by convolutional layers. 
The main idea behind ResMLP is the use of a residual block illustrated in Figure \ref{residual_block}. The skip connection models the identity operator between layers and adds the outputs from previous layers to the outputs of stacked layers. This results in the ability to train much deeper networks by reducing the risk of vanishing gradients. 
The use of ResMLP allows us to increase the number of hidden layers without introducing the effect of vanishing gradients which allows ResNet to model highly non-linear mappings such as the ones arising from adaptive zoning.
The architecture of the ResMLP we use in this work is made of 5 residual blocks, each containing a correction and an identity hidden layer. Both correction and identity hidden layers have 100 neurons and the rectified linear unit (ReLU) is used at the end of each residual block to enforce non-linearity of the model. 

\begin{figure}[H]
    \centering
    \includegraphics[width=0.8\textwidth]{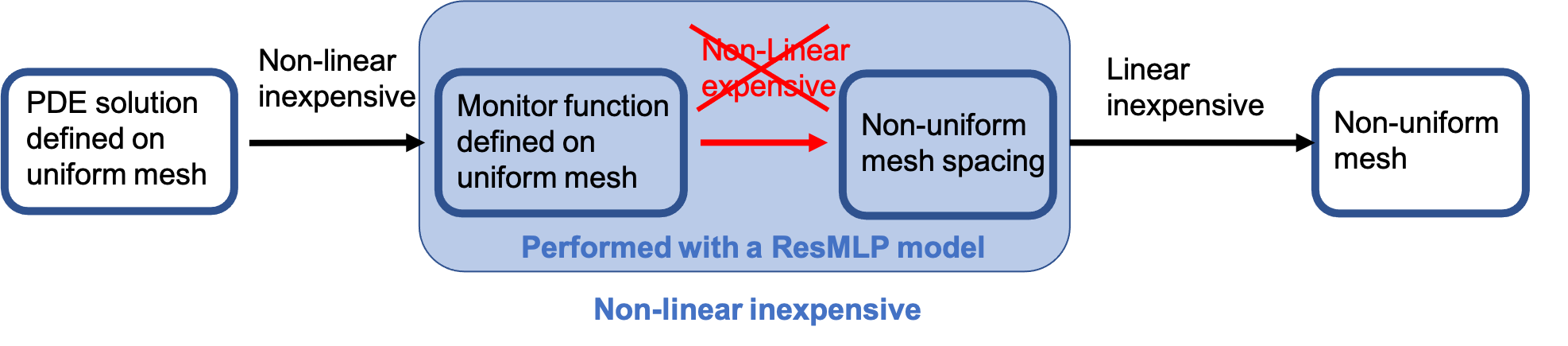}
    \caption{Sequence of computational tasks performed to map the solution to the PDE on a uniform mesh with the final adapted mesh. A ResMLP architecture is used as a surrogate DL model for the expensive task that maps the monitor function to the mesh spacing of the adapted mesh.}
    \label{pipeline_DL}    
\end{figure}

\subsection{Setup of deep learning training }
The dataset used to train and validate the ResMLP model is generated by running the standard adaptive zoning for several artificial shock profiles, each characterized by an arbitrary number of jumps at arbitrary locations in the domain. The size of the 1D mesh is set to 201 nodes (extreme points of the domain included) for each test cases included in the training dataset. 
The choice between Equation \eqref{elliptic_bvp} and \eqref{parabolic_bvp}, as well as the values for $\alpha_1$, $\alpha_2$, $\alpha_3$, and $\alpha_4$ depend on the physics of the problem. 
The training of a ResMLP model as surrogate for adaptive zoning can be performed using different options as input and output data. Different choices lead to conglomerating different steps in the workflow pipeline described in Figure \ref{pipeline_DL}. 
The input data in the training set can either be the shock profile itself or the monitor function evaluated on the initial uniform mesh. 
The output data in the training set can either be the coordinates of the nodes on the adapted mesh or the non-uniform mesh spacing. This leads to four possible combinations of input and output data that are completely arbitrary and independent of the physics problem.
We tested the performance of the trained DL model using all the four possible combinations of input data and output data in the training set, and the combination that resulted into the best performing ResMLP model uses the monitor function as input and the non-unifor mesh spacing as output. 
The set of parameters varied during the generation of the dataset with artificial shock profiles are described in Table \ref{tab:parameters_data}.
The dataset is split in training and validation portions with the percentages of $80\%$ and $20\%$ of the total dataset, respectively. 
Denoting the vector that collects the output data with $\mathbf{y}$ and the vector that collects the predictions produced by the NN with $\hat{\mathbf{y}}$, the training of the NN is performed using the mean absolute error (MAE) 
\begin{equation}
MAE(\mathbf{y}, \mathbf{\hat{y}}) = \frac{1}{N} \sum_{i=1}^N \lvert y_i -\hat{y}_i\rvert 
\end{equation}
as loss function and Adam \cite{adam} is used as stochastic optimizer. The total number of training epochs is $5,000$ and the data is iteratively processed in batches of size $100$. 

\begin{table}[H]
    \centering
    \begin{tabular}{|c|c|}
    \hline
    \textbf{Parameter} & \textbf{Range}\\
    \hline
         number of randomly located shocks & $\{1,2,3,4\}$  \\
         range of values for each shock profile & $[0,+1]$\\
    \hline
    \end{tabular}
    \caption{Description of the parameter space to generate the dataset.}
    \label{tab:parameters_data}
\end{table}

\section{Numerical results}
\label{numerical_section}

\subsection{Hardware and software description}
The numerical results are obtained using a \textregistered{MacBookPro15,3}.
The 1D physics PDE solvers are implemented in \textregistered{Matlab} \cite{MATLAB:2019}, and all the DL computational tasks are performed using \texttt{Python 3.6}. The 1D physics PDE solvers in \textregistered{Matlab} interacts with the validated DL model implemented in \texttt{Python} using the \textregistered{Matlab} function \texttt{py.importlib.import\_module}. 

\subsection{Communication overhead between \textregistered{Matlab} and \texttt{Python}}
The physics PDE solver implemented in \textregistered{Matlab} loads from the hard disk the DL model implemented in \texttt{Python} at each time step and this introduces communication overhead that impacts the total computational time. We measure the communication overhead by subtracting the average time spent for an evaluation of the DL model in \texttt{Python} from the average time spent to perform the same DL evaluation in \textregistered{Matlab}. Performing 1,000 DL evaluations in \texttt{Python} took 21.10 wall-clock seconds, whereas performing it in \textregistered{Matlab} took 168.70 wall-clock seconds. Therefore, we estimate the communication overhead per time step to be (168.70 - 21.10)/1,000 = 0.1476 wall-clock seconds.

\subsection{Sensitivity of DL performance on volume of training data}
\begin{figure}[H]
    \centering
    \includegraphics[width=0.45\textwidth]{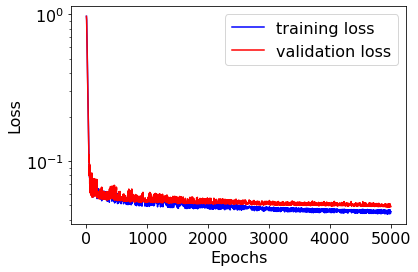}
    \includegraphics[width=0.45\textwidth]{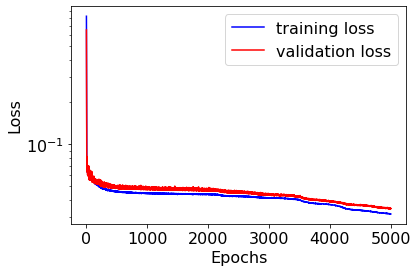}
    \includegraphics[width=0.45\textwidth]{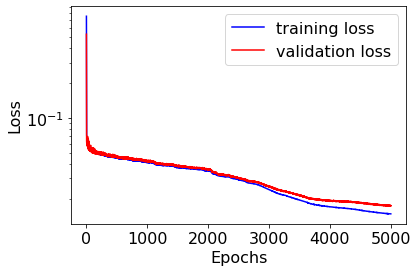}    
    \includegraphics[width=0.45\textwidth]{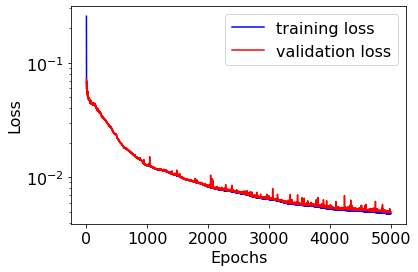}
    \caption{Training and validation MAE for training of DL models on dataset generated with 500 artificial shock profiles (top-left), 5,000 artificial shock profiles (top-right), 10,000 artificial shock profiles (bottom-left), and 50,000 artificial shock profiles (bottom-right). Training input data: monitor function. Training output data: mesh spacing. }
    \label{loss_functions}
\end{figure}

We first study how the predictive accuracy of ResMLP improves when the DL model is trained on datasets of increasing size. We generate the dataset using the elliptic MMPDE in Equation \eqref{elliptic_bvp} and the monitor function in Equation \eqref{monitor_function} to perform the standard adaptive zoning. The parameters used for the standard adaptivity and generate the training dataset for DL are described in Table \ref{parameters_mmpde1}.
The numerical results in Figure \ref{loss_functions} compare the trend of the training and validation loss functions for datasets with $500$, $5,000$, $10,000$, and $50,000$ data samples. The final accuracy improves for increasing sizes of the dataset, and the validation loss function follows the trend of the training loss function, which confirms good generalizability of the DL model. 

\subsection{Square wave test case}
\begin{figure}
    \includegraphics[width=0.32\textwidth]{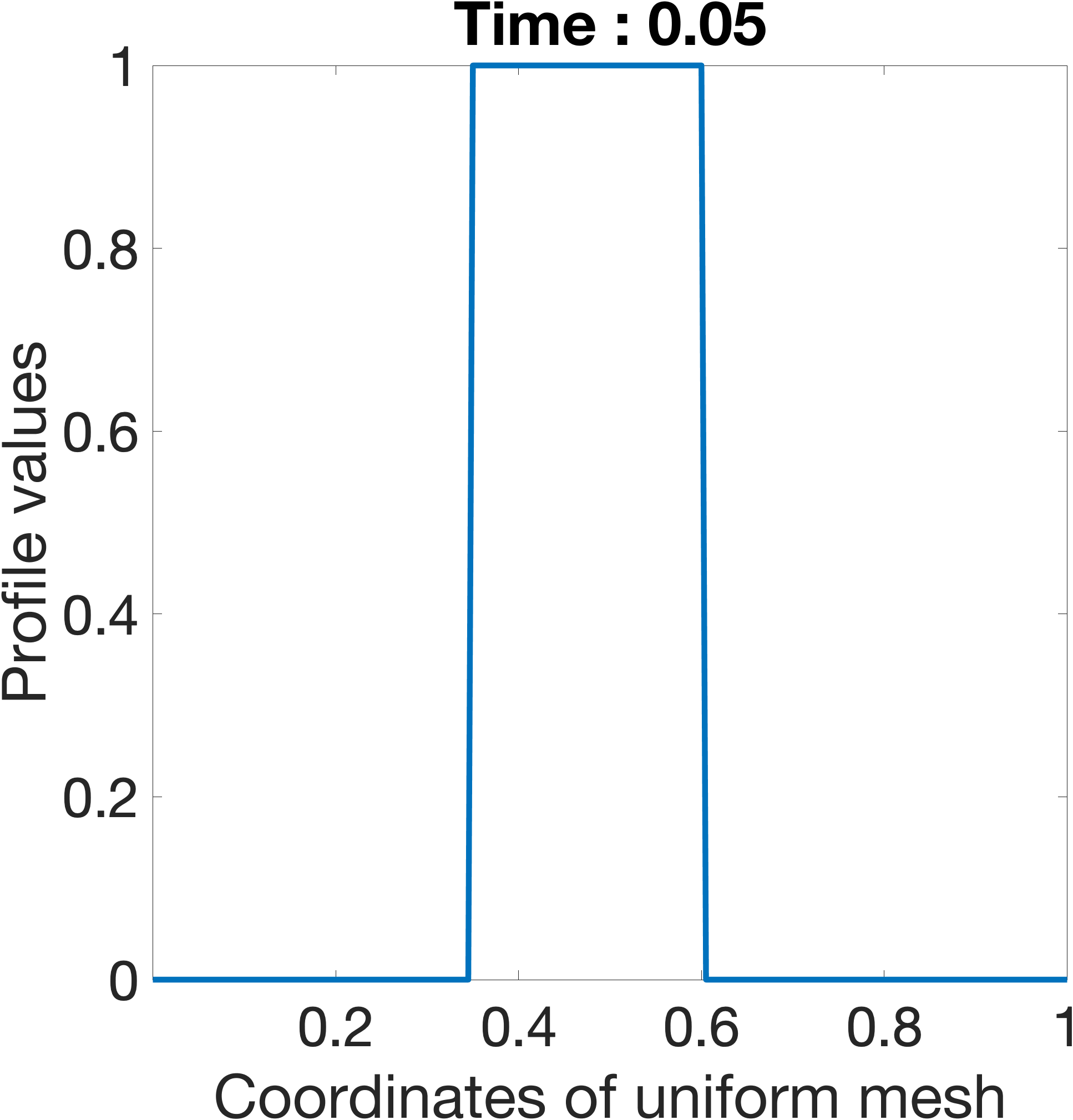}
    \includegraphics[width=0.32\textwidth]{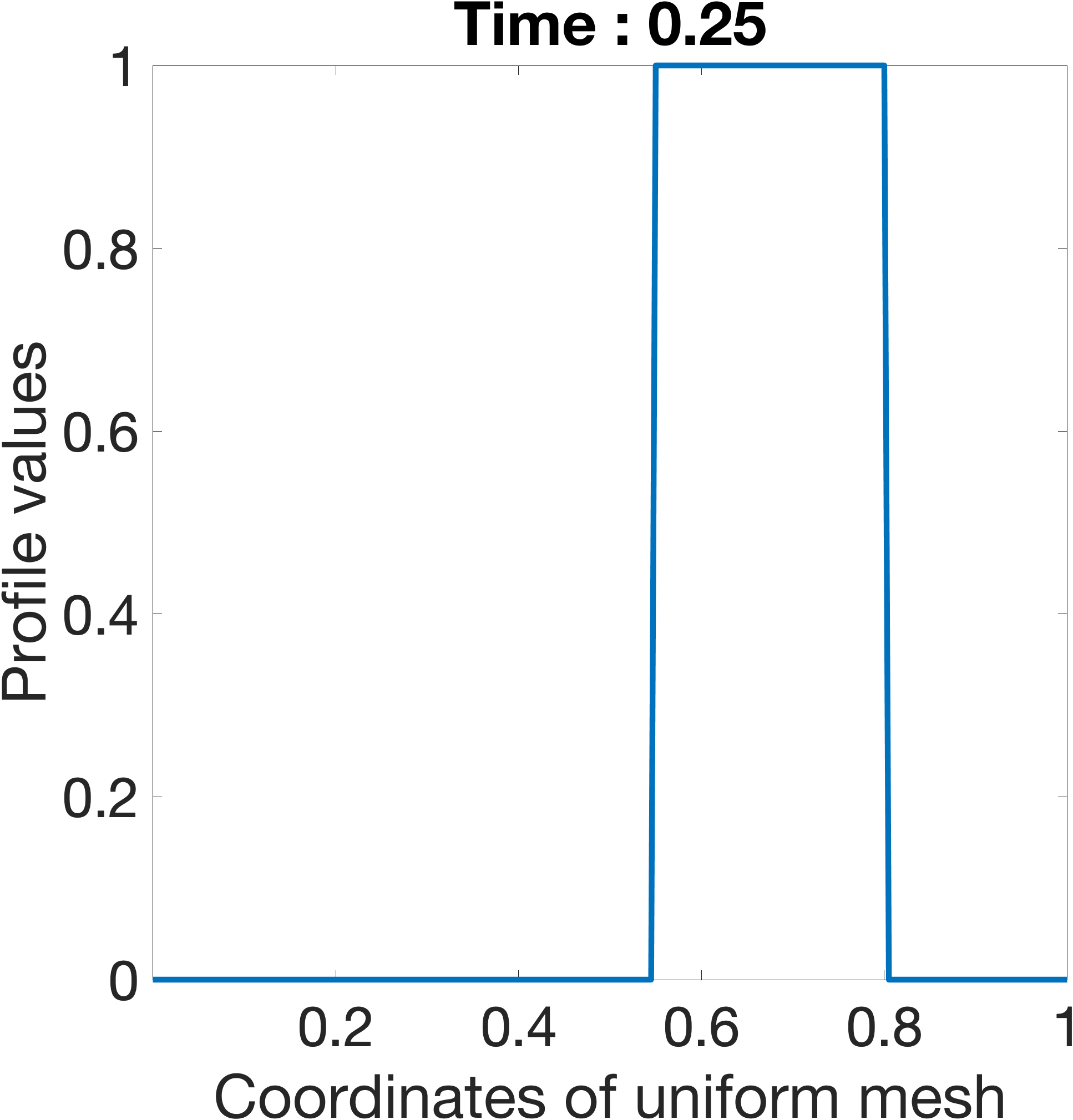}
    \includegraphics[width=0.32\textwidth]{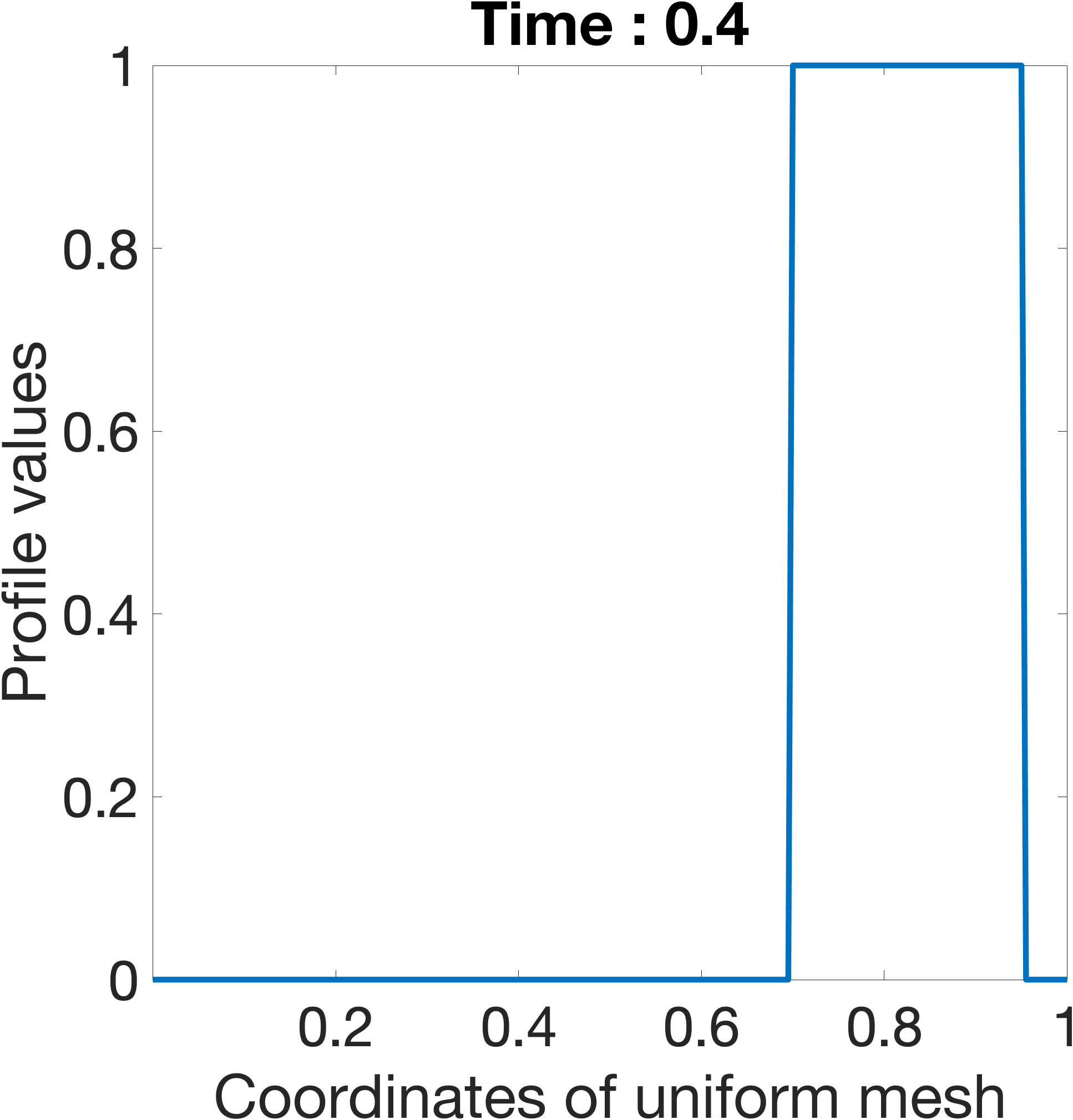}
    \includegraphics[width=0.32\textwidth]{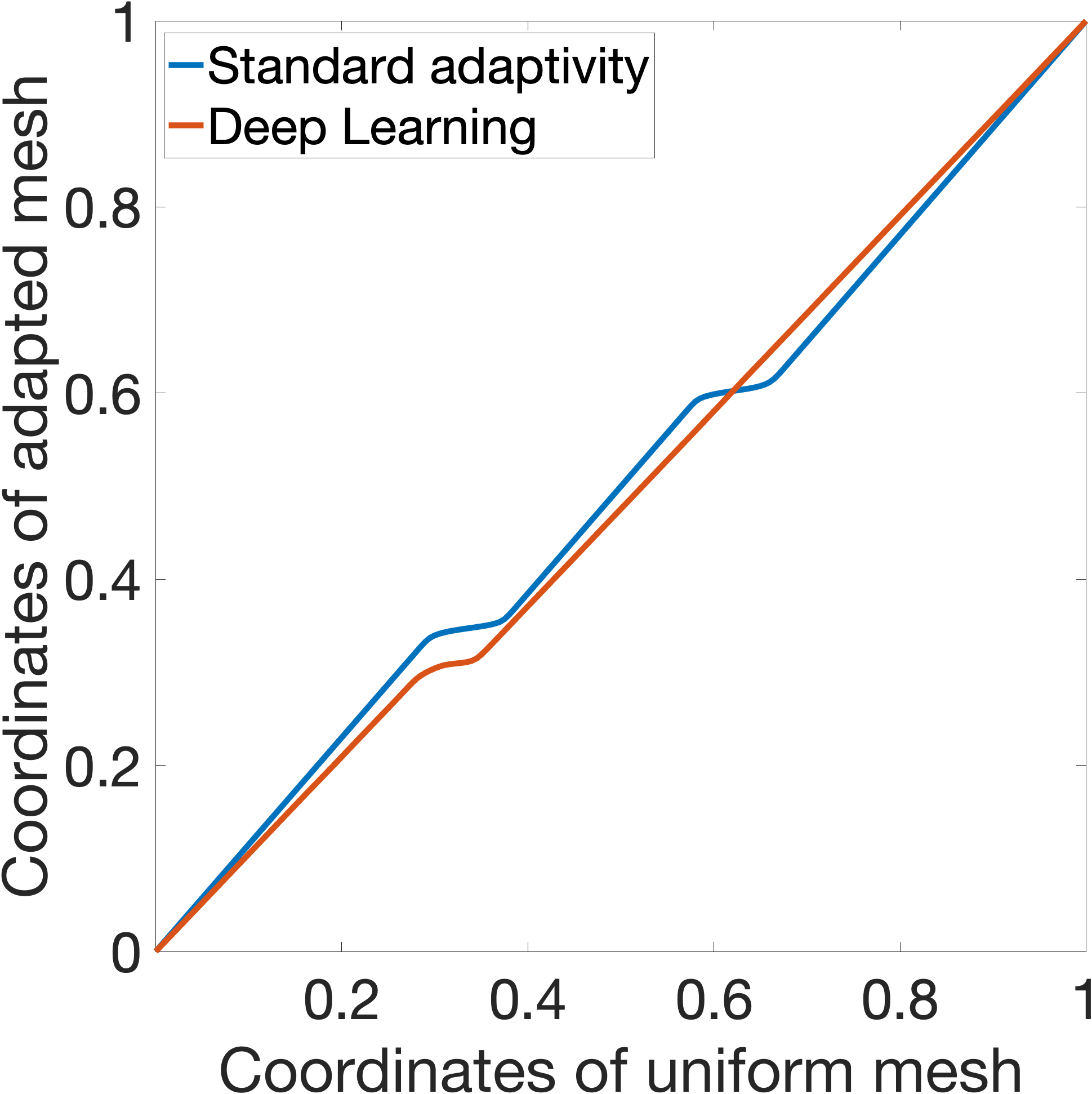}
    \includegraphics[width=0.32\textwidth]{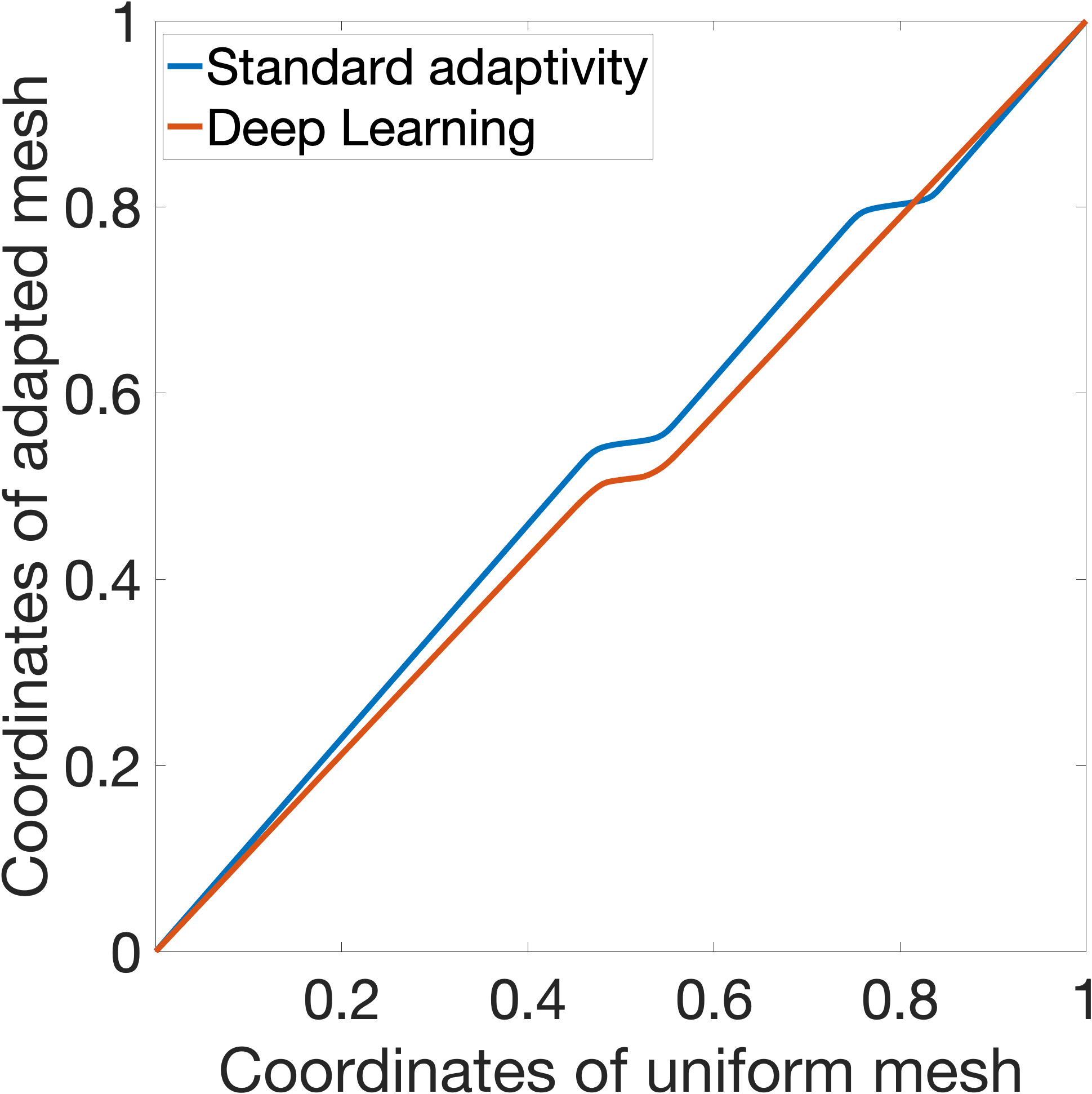}
    \includegraphics[width=0.32\textwidth]{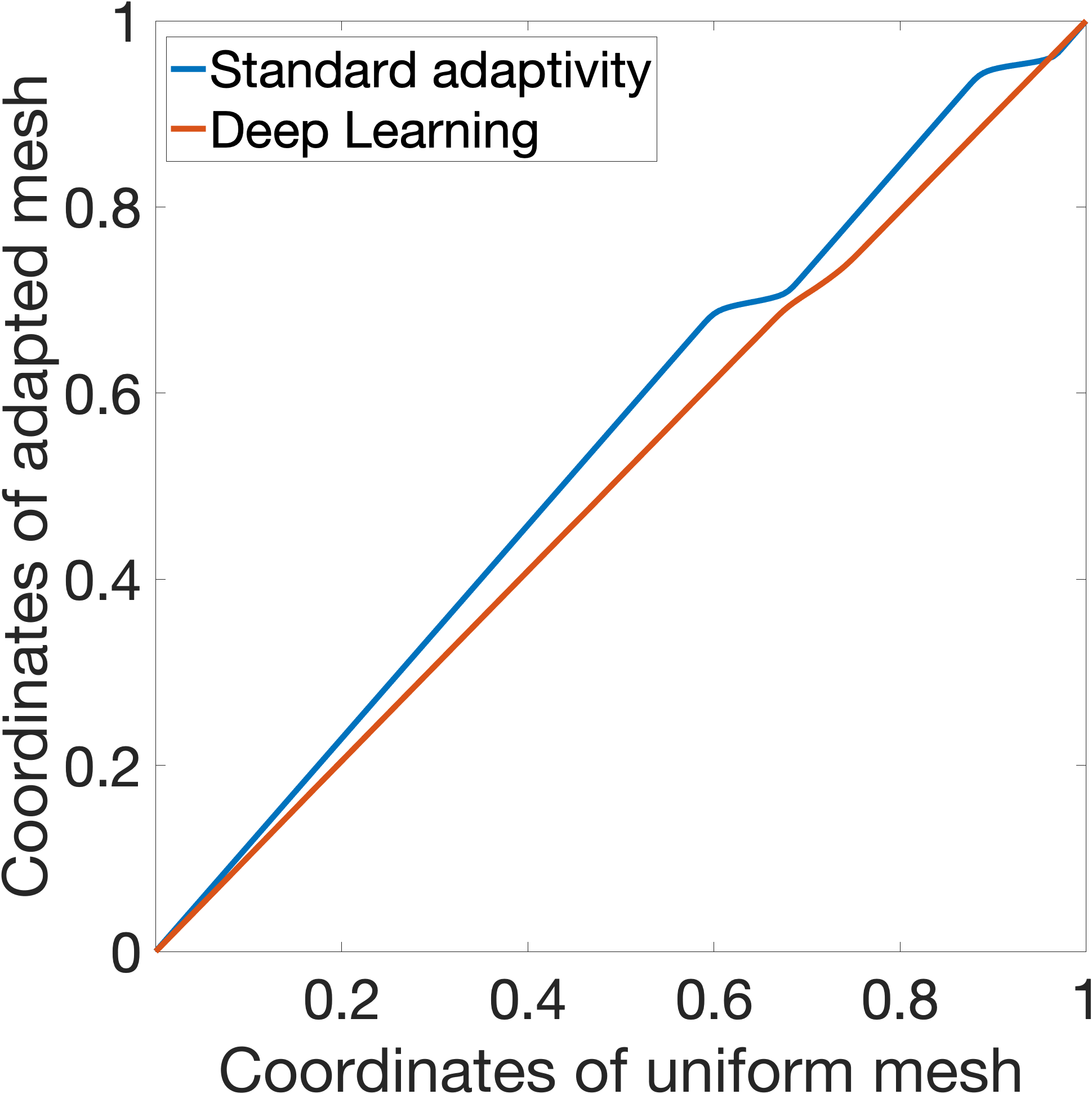}
    \includegraphics[width=0.32\textwidth]{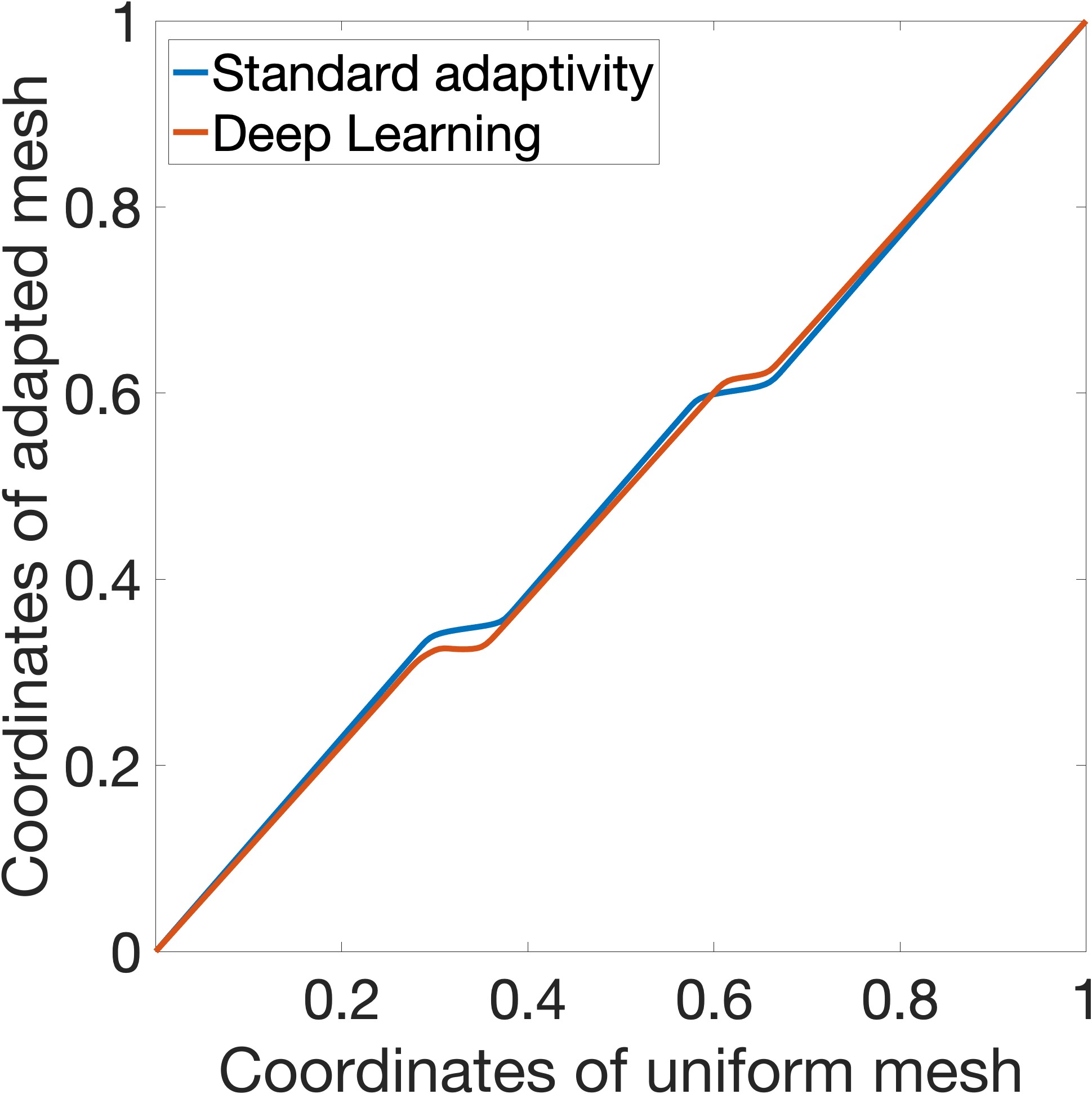}
    \includegraphics[width=0.32\textwidth]{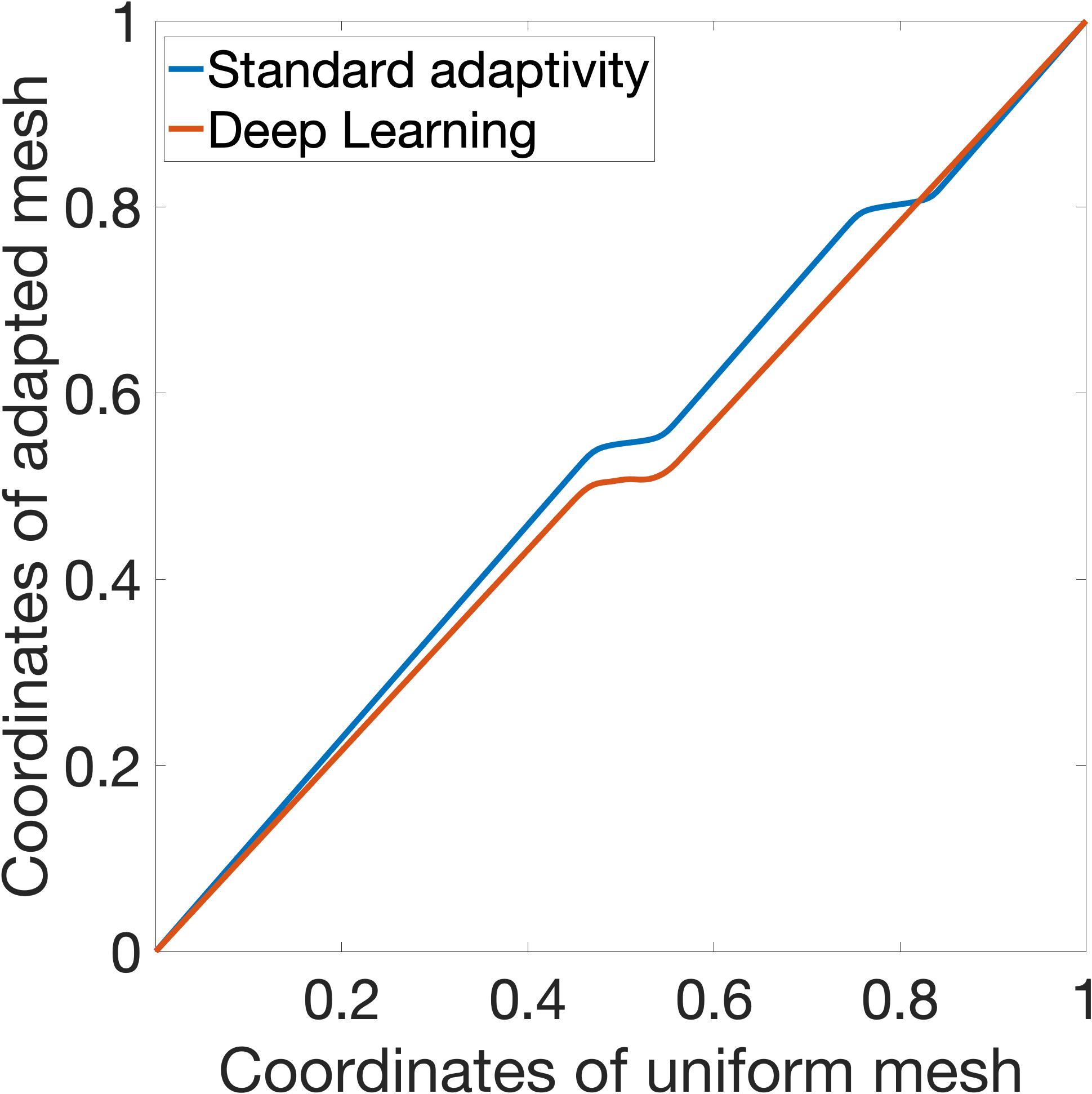}
    \includegraphics[width=0.32\textwidth]{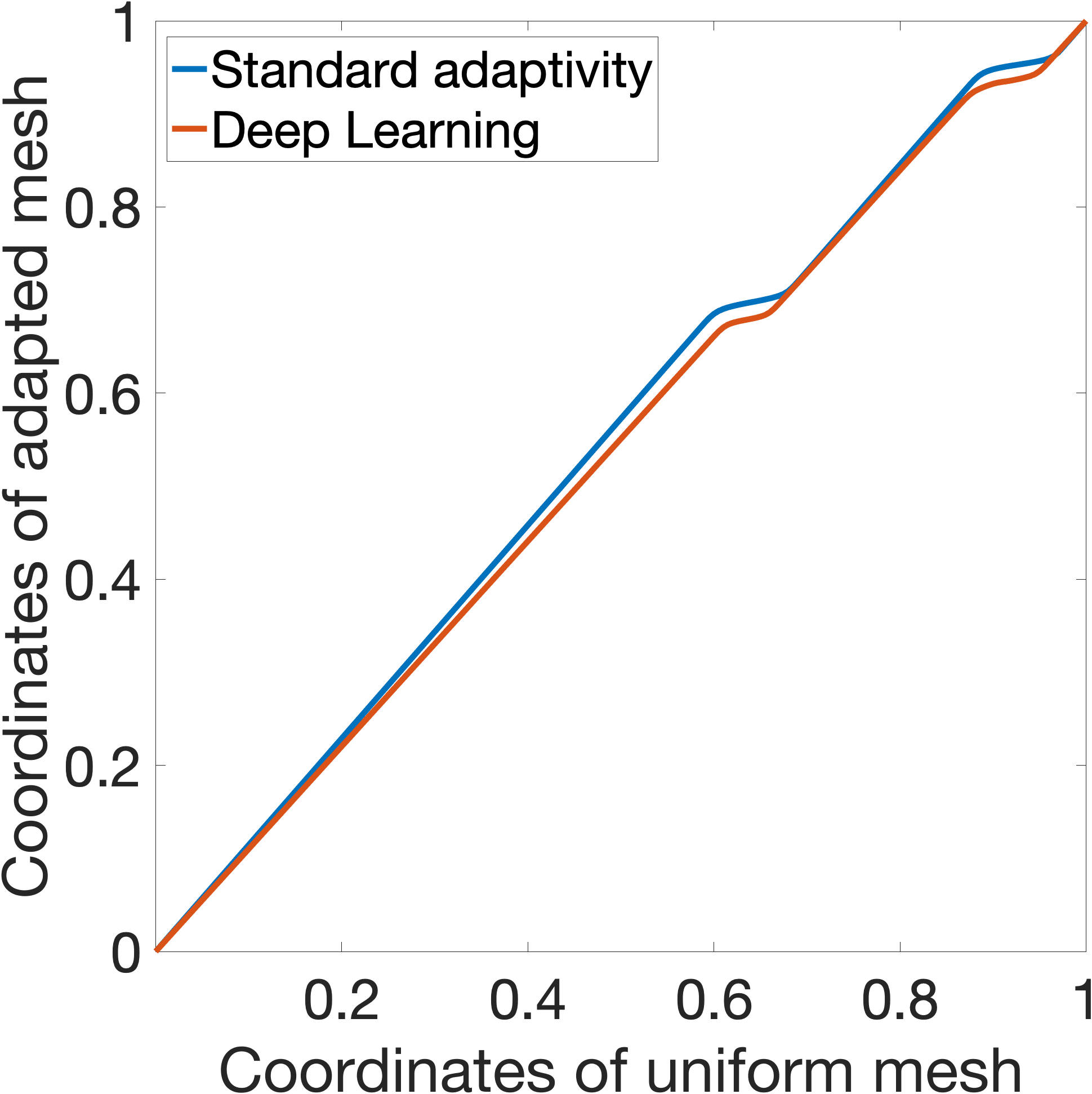}
    \includegraphics[width=0.32\textwidth]{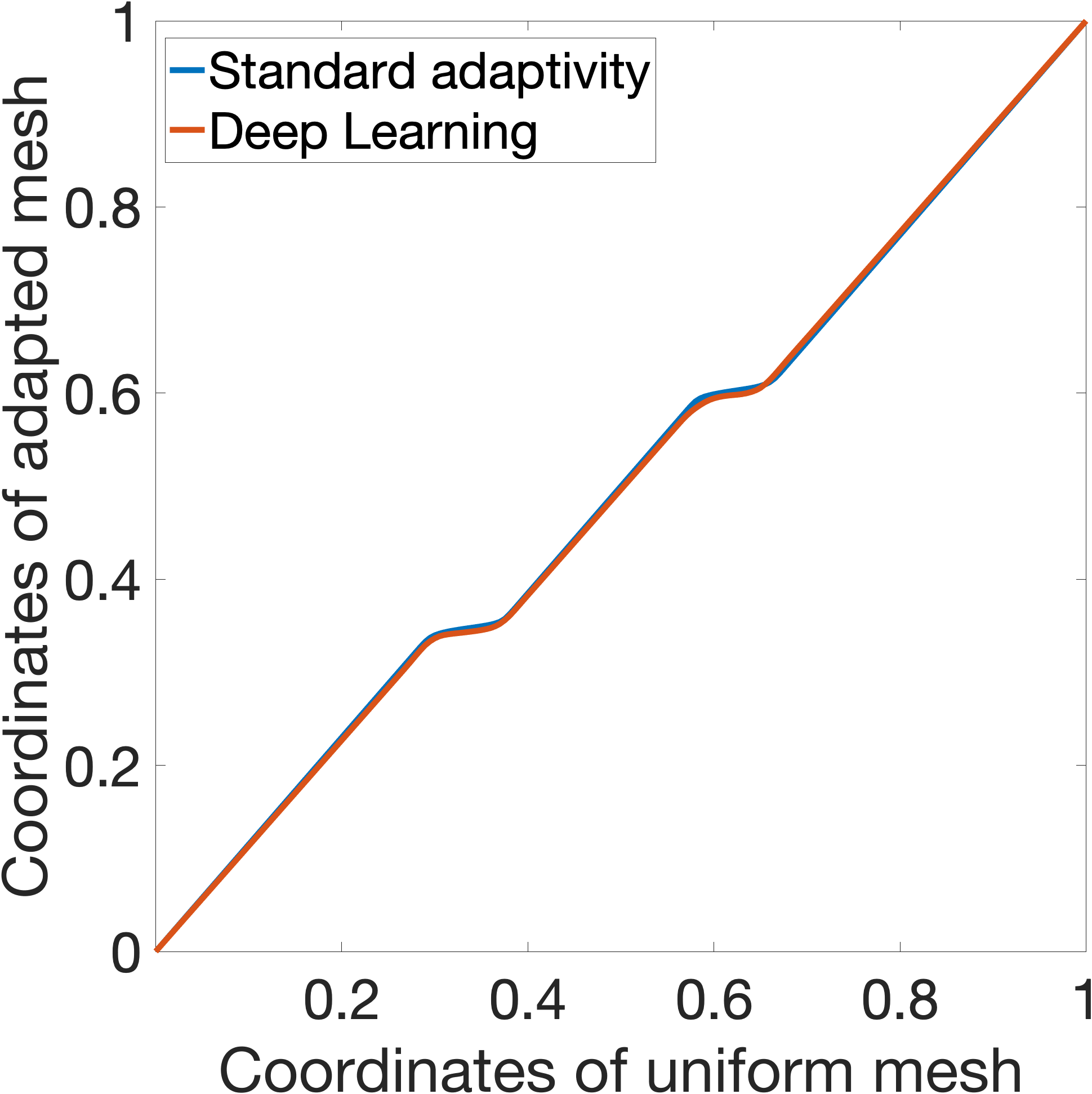}
    \includegraphics[width=0.32\textwidth]{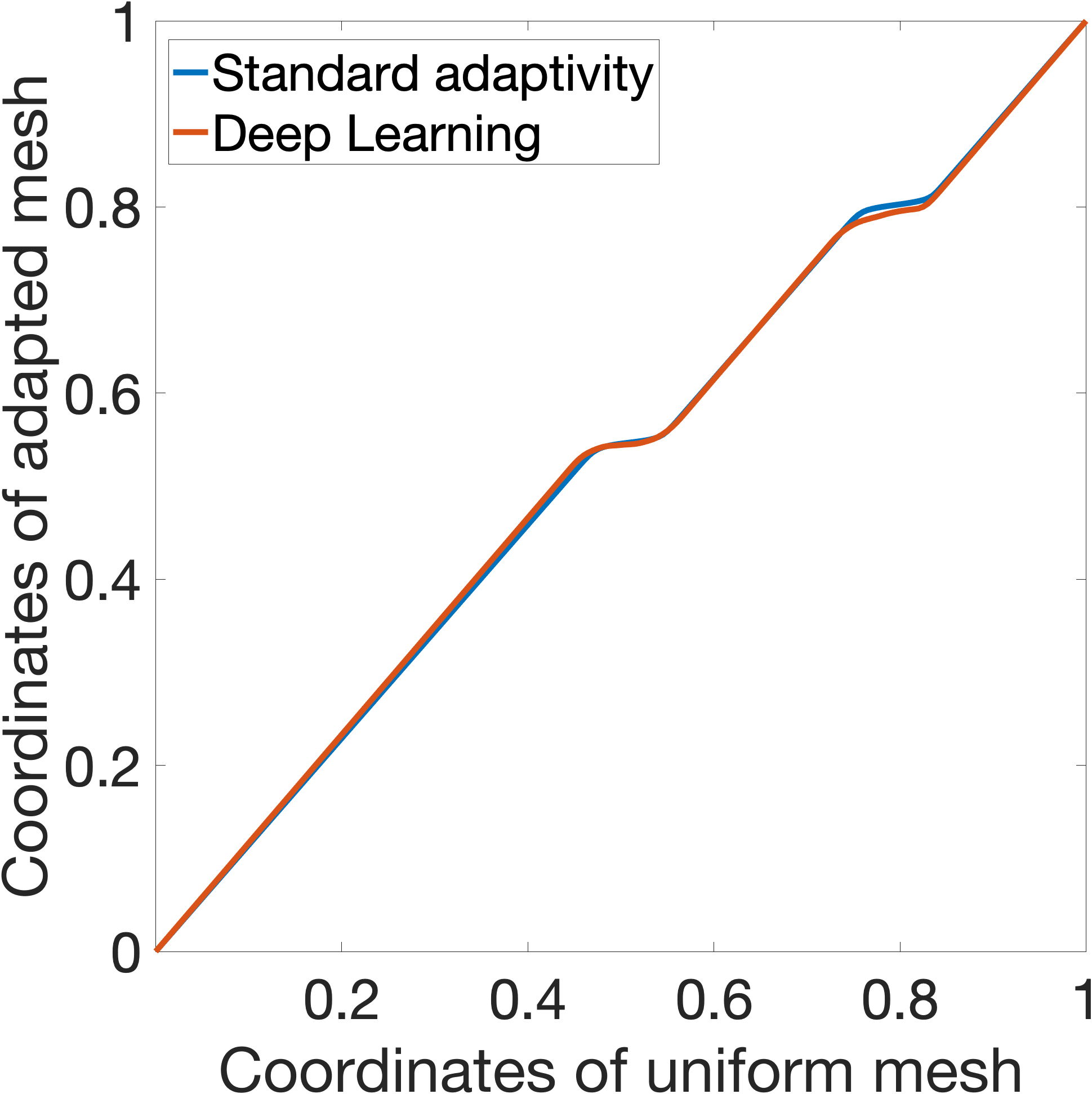}
    \includegraphics[width=0.32\textwidth]{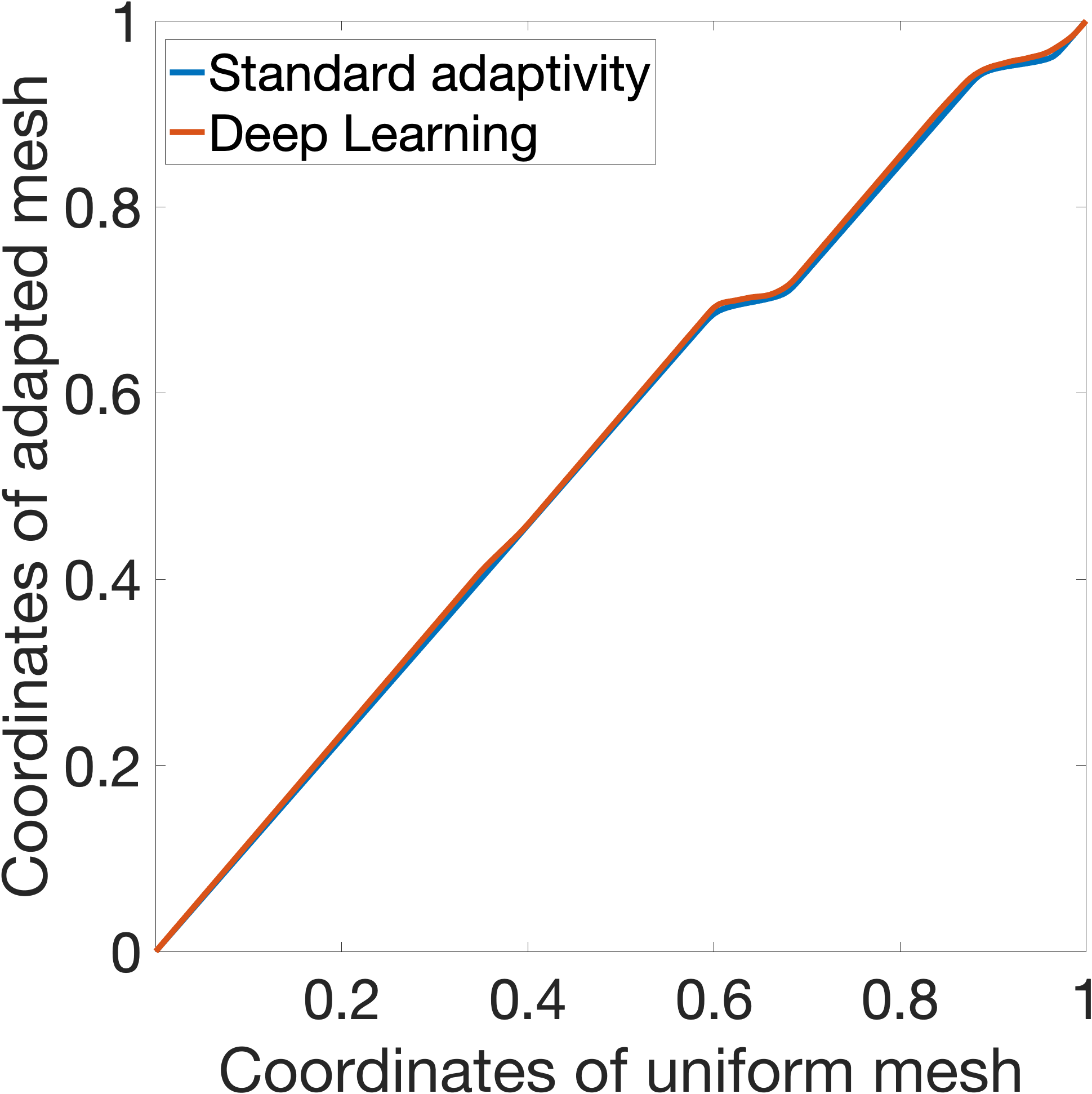}    
    \caption{Training input data: monitor function. Training output data: mesh spacing. DL adaptive zoning (right) compared with standard adaptive zoning (center) for propagating square wave at time $t=0.05$ (left), $t = 0.25$ (center) and $t=0.4$ (right). The comparison is performed by training the DL model on a data set with 500, $5,000$, and $50,000$ artificial shock profiles (from top down).}
    \label{meshes_changesize}
\end{figure}

Figure \ref{meshes_changesize} shows the performance of the DL model trained on different volumes of data to adaptively relocate the mesh nodes based on the location of the discontinuitites for a propagating square wave, reconfirming that the DL model improves its performance in capturing the discontinuities at different locations for larger volumes of training data. Although training the ResMLP model on larger training datasets makes the training more computationally expensive, this offline cost has to be paid only once upfront, and it is then amortized by the online deployment of the trained model.

\begin{table}
\begin{center}
\begin{tabular}{ |c| c| }
\hline
\multicolumn{2}{|c|}{\textbf{Elliptic MMPDE }}\\
\hline
Definition of monitor function & Equation \eqref{monitor_function}\\
\hline
 $\alpha_1$ & 0.0 \\ 
 \hline
 $\alpha_2$ & $600.0$ \\ 
 \hline
 $\alpha_3$ & 0.0 \\  
 \hline
 $\alpha_4$ & 0.0 \\  
 \hline
 smoothing & None\\
 \hline
 fixed point tolerance & $1\mathrm{e}-8$ \\
 \hline
 maximum number of fixed point iterations & 1,000 \\
 \hline 
\end{tabular}
\caption{Parameters of the monitor function and fixed point solver to perform standard adaptivity and generate the training dataset for DL using the elliptic MMPDE in the square wave propagation test case, Sod shock tube test case, and Taylor–von Neumann–Sedov blast wave test case.}
\label{parameters_mmpde1}
\end{center}
\end{table}

\subsection{Test cases with 1D Euler equations}
We present numerical results for three test cases described by 1D Euler equations in Equation \eqref{euler}: the Sod shock tube test case \cite{sod1978}, the Taylor–von Neumann–Sedov blast wave test case \cite{sedov1946}, and the Woodward-Colella blast waves test case \cite{woodward}. The DL adaptive zoning for these test cases is performed using the same ResMLP model used for the square wave test case. The ResMLP model to perform DL adaptivity on the Sod shock tube test case and the Taylor–von Neumann–Sedov blast wave test case is trained on a dataset obtained by solving the elliptic MMPDE Equation \eqref{elliptic_bvp} with the monitor function defined in Equation \eqref{monitor_function} for a set of 50,000 staircase profiles, whereas the dataset used to train the ResMLP for the Woodward-Colella blast waves test case used the parabolic MMPDE \eqref{parabolic_bvp} with the monitor function \eqref{monitor_function}. The time steps are adaptively adjusted to respect the Courant–Friedrichs–Lewy (CFL) condition \cite{cfl} using a CFL number equal to 0.6.

\subsubsection{Sod shock tube test case}
\label{sod_test_case}
We now compare the performance of DL against the elliptic MMPDE for the 1D Sod shock tube problem, which is described by the 1D Euler equations \eqref{euler} in the spatial domain [0,1] with the following initial conditions
\begin{equation*}
\begin{aligned}
    &\rho(x, t=0) = 
    \begin{cases}
     1.0\quad \text{if} \quad x<=0.5\\
     \frac{1}{8}\quad \text{if} \quad x>0.5\\
    \end{cases},\\&
    u(x, t=0) = 0.0, \\&
   p(x, t=0) = 
    \begin{cases}
     1.0\quad \text{if} \quad x<=0.5\\
     \frac{1}{10}\quad \text{if} \quad x>0.5\\
    \end{cases}\\    
\end{aligned}
\end{equation*}
and non-homogeneous Dirichlet boundary conditions 
\begin{equation*}
\rho(0,t) = 1.0, \quad u(0,t) = 0.0, \quad p(0,t) = 1.0, \quad 
\rho(1,t) = \frac{1}{8}, \quad u(1,t) = 0.0, \quad p(1,t) = \frac{1}{10} \quad t\in [0,T).
\end{equation*}
The initial and boundary values of the energy are obtained from the state equations described in Section \ref{background}. The domain [0,1] is discretized with 200 finite volume cells and we present results using both the Lax-Wendroff and the WENO5 schemes for the space discretization combine with RK-3 for time stepping. The time integration is performed using RK-3 and is stopped at $t=0.2$.
The parameters used for the standard adaptivity are described in Table \ref{parameters_mmpde1}. 
The DL model is the same as for the square wave propagation test case. 

\begin{figure}[H]
\centering
\includegraphics[width=0.45\textwidth]{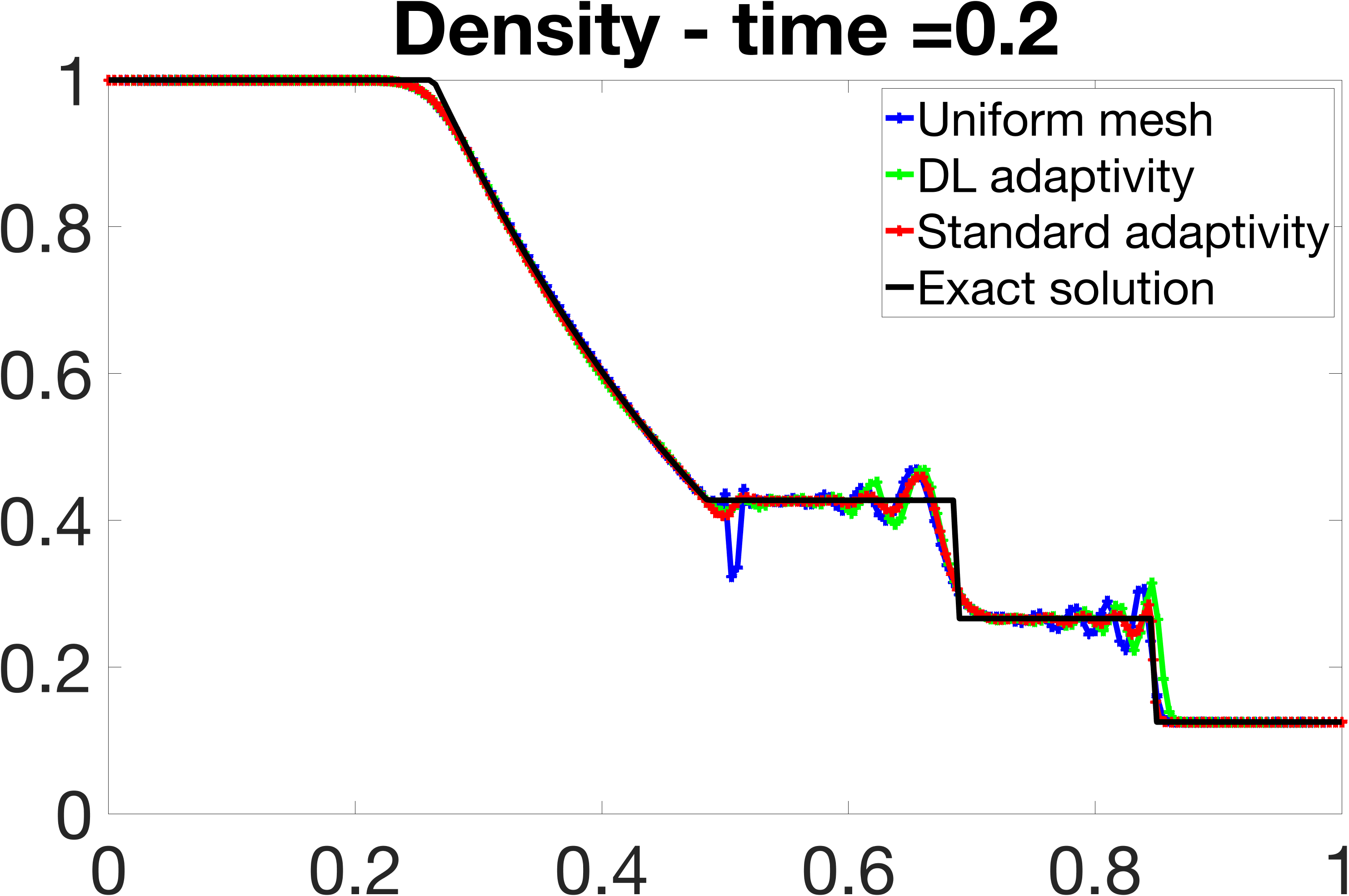}
\hspace{1cm}
\includegraphics[width=0.45\textwidth]{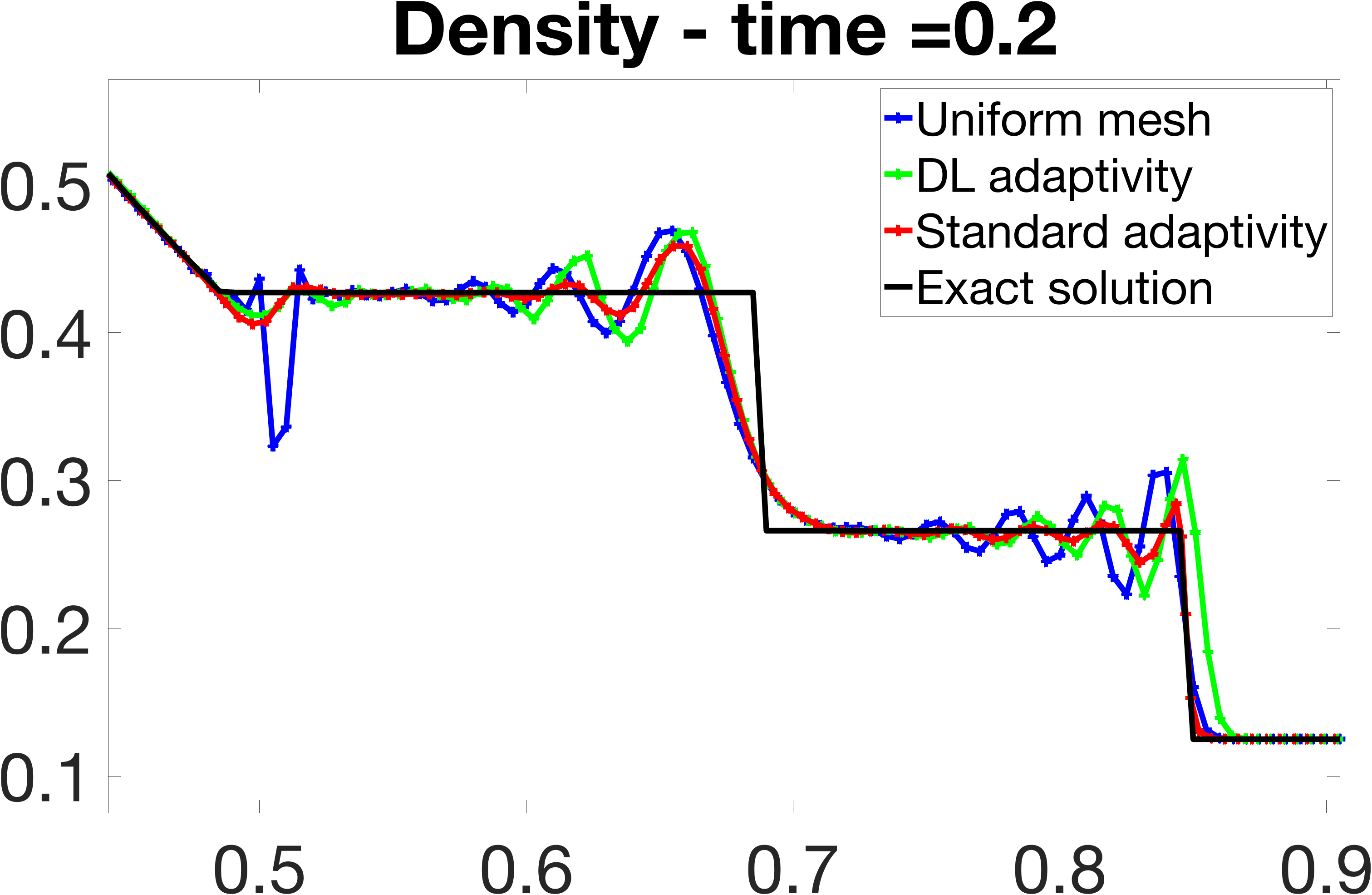}   
\caption{Density profile of 1D Sod test case at $t=0.2$. Lax Wendroff scheme used for space discretization and time stepping. Entire profile shown at left and zoom-in in the $x$ range [0.5, 0.9] shown at the right. }
\label{euler_lw_density}
\end{figure}

\begin{figure}[H]
\centering
\includegraphics[width=0.45\textwidth]{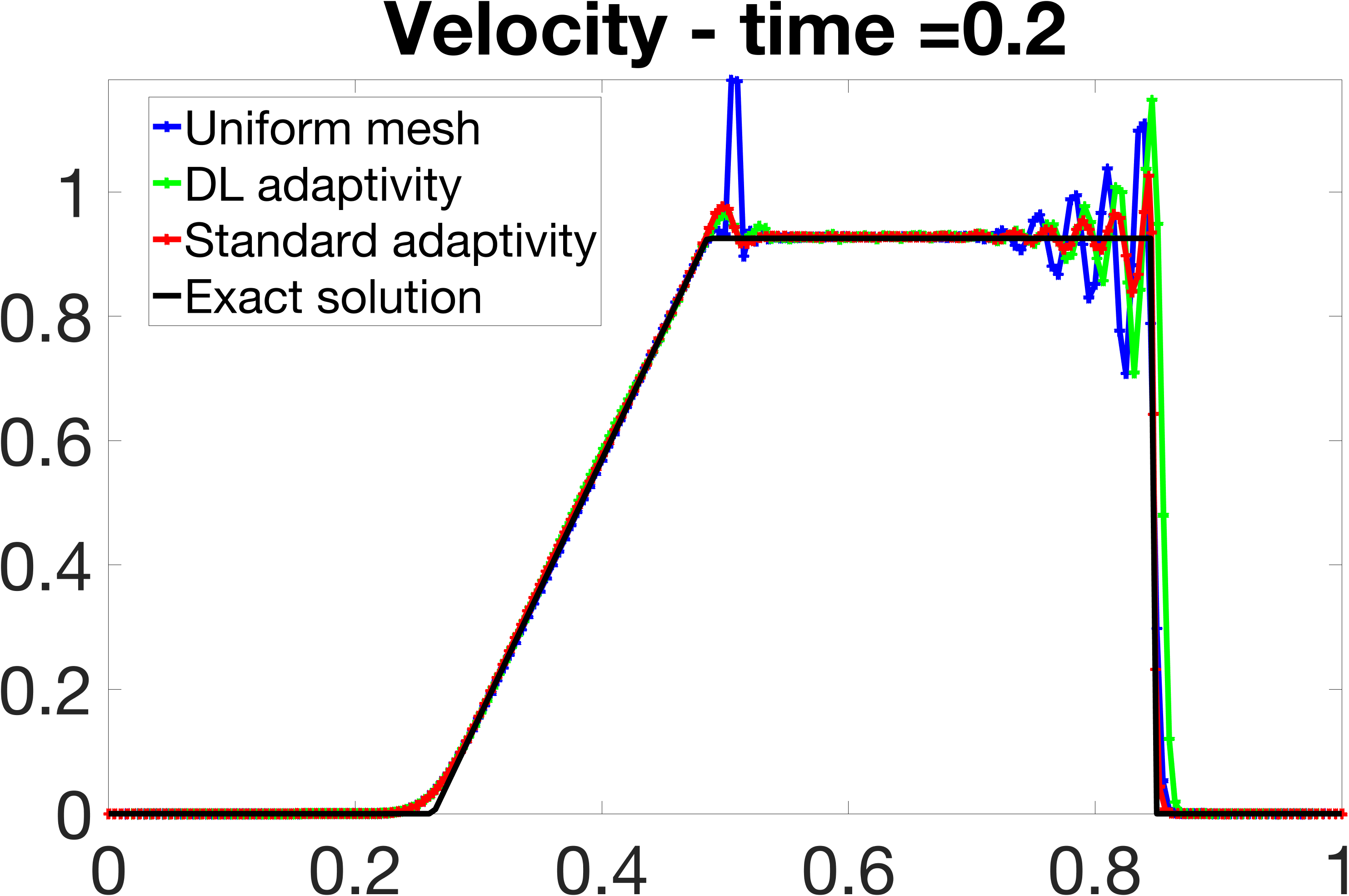}   
\hspace{1cm}
\includegraphics[width=0.45\textwidth]{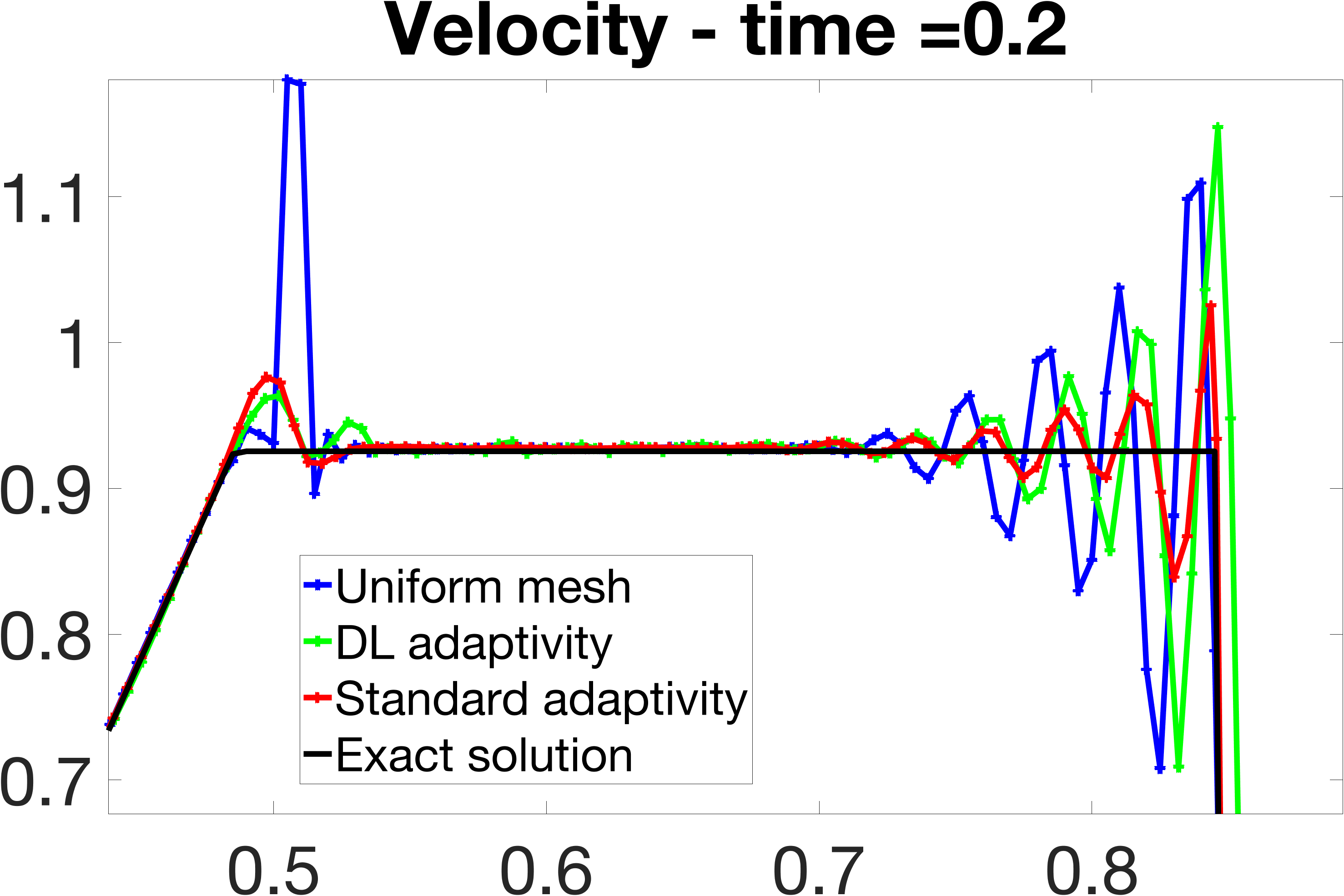}   
\caption{Velocity profile of 1D Sod test case at $t=0.2$. Lax Wendroff scheme used for space discretization and time stepping. Entire profile shown at left and zoom-in in the $x$ range [0.5, 0.9] shown at right. }
\label{euler_lw_velocity}
\end{figure}

\begin{figure}[H]
\centering
\includegraphics[width=0.45\textwidth]{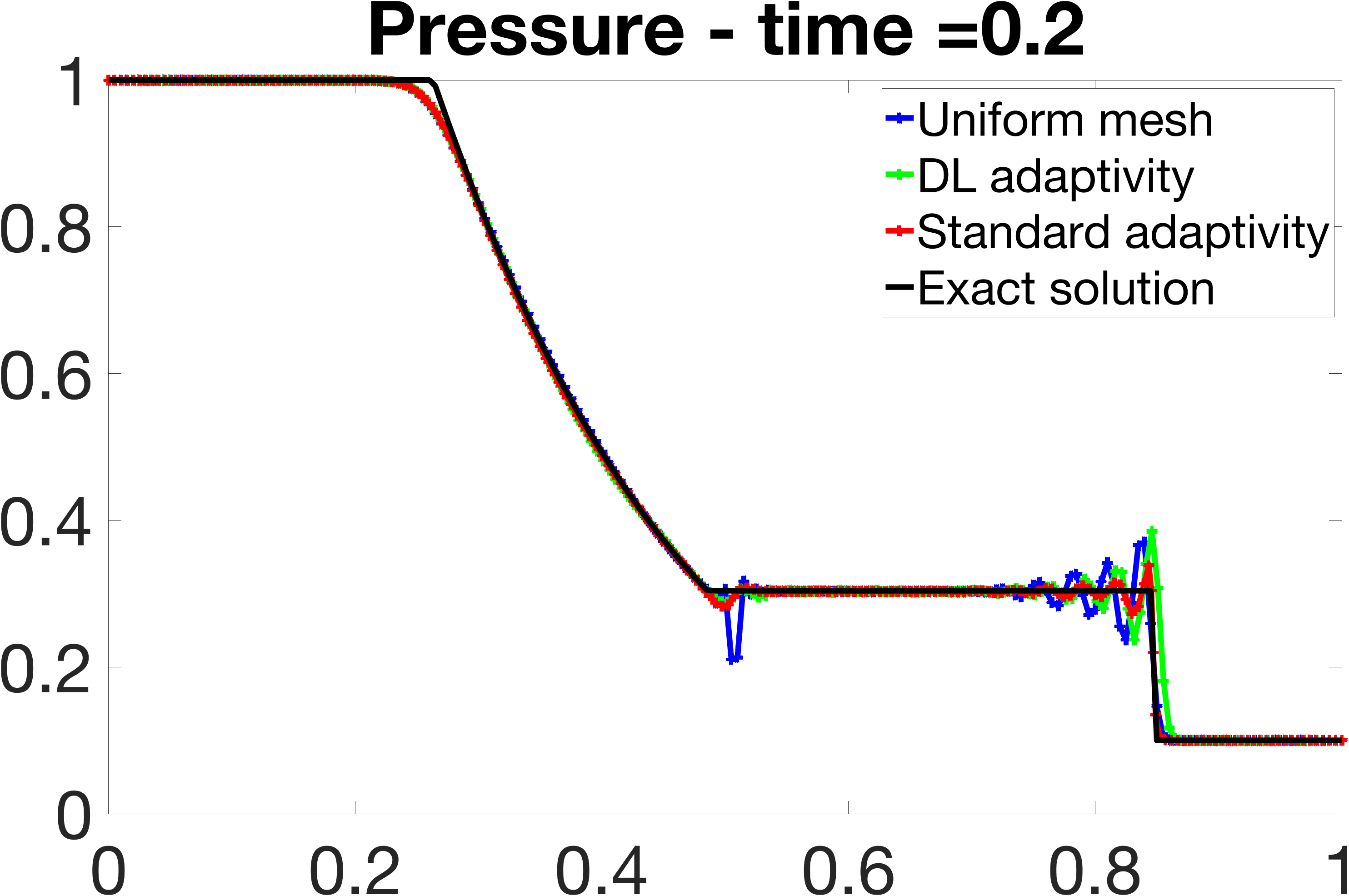}   
\hspace{1cm}
\includegraphics[width=0.45\textwidth]{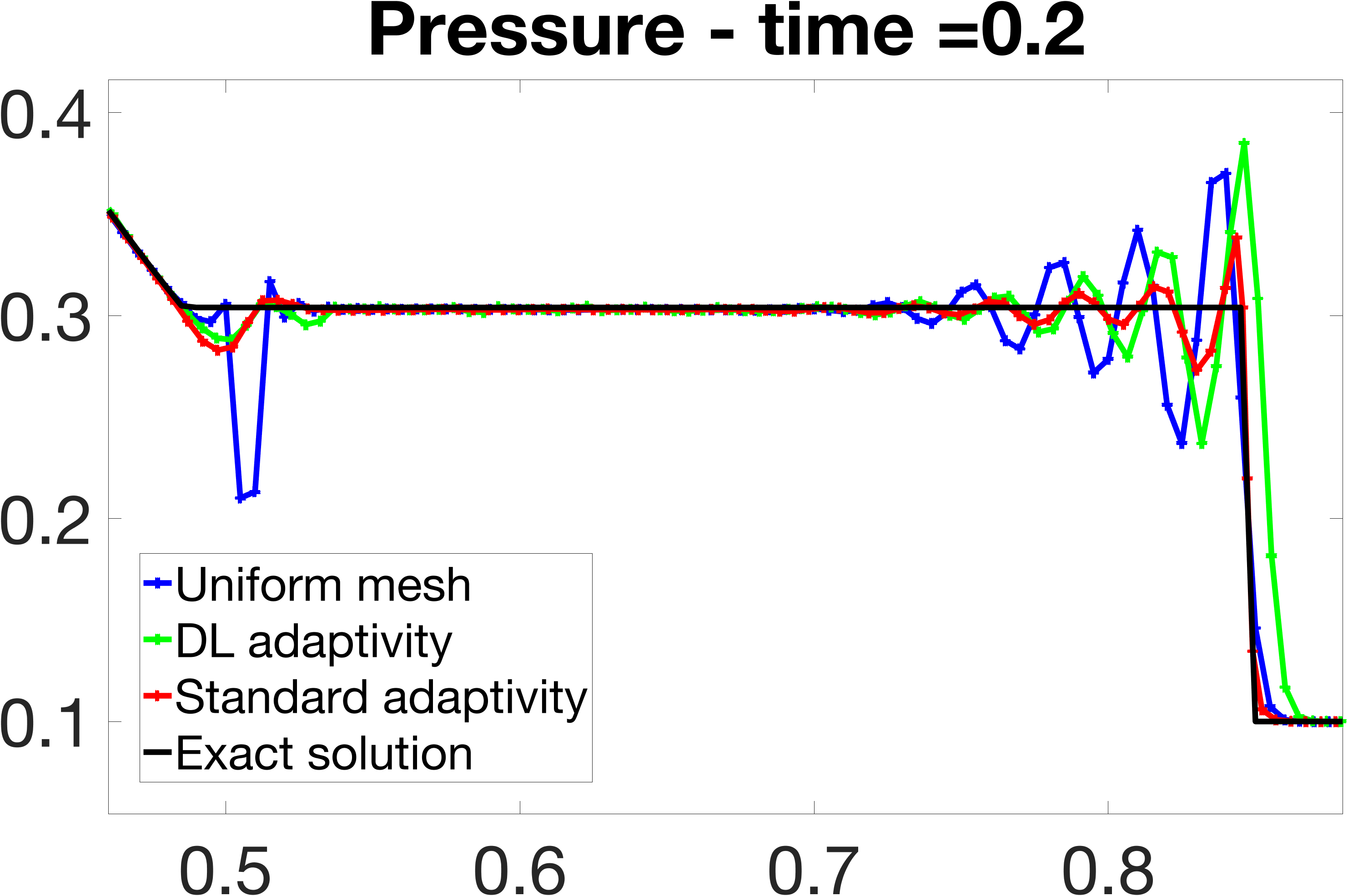}   
\caption{Pressure profile of 1D Sod test case at $t=0.2$. Lax Wendroff scheme used for space discretization and time stepping. Entire profile shown at left and zoom-in in the $x$ range [0.5, 0.9] shown at right. }
\label{euler_lw_pressure}
\end{figure}

\begin{figure}[H]
\centering
\includegraphics[width=0.45\textwidth]{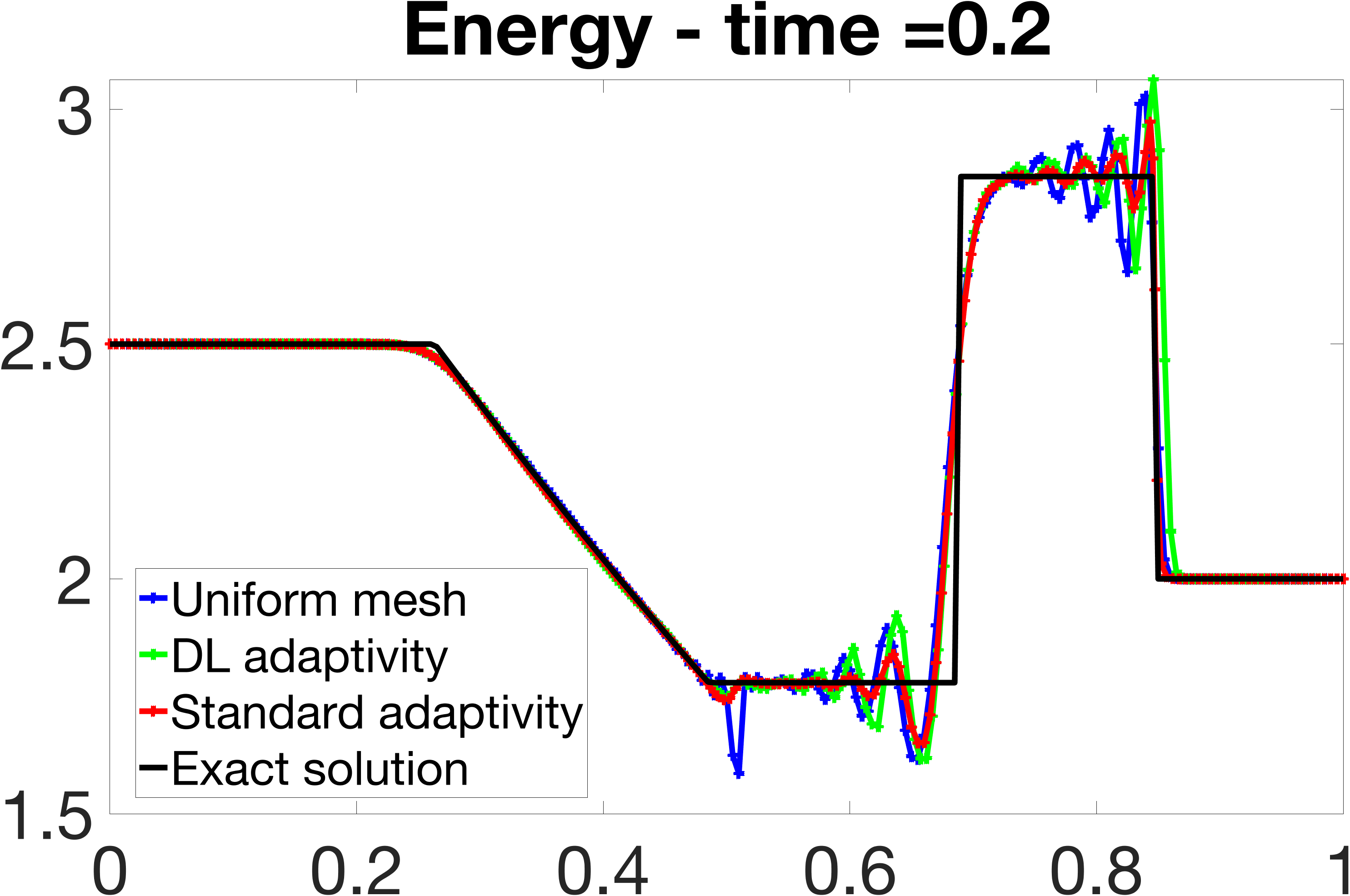}  
\hspace{1cm}
\includegraphics[width=0.45\textwidth]{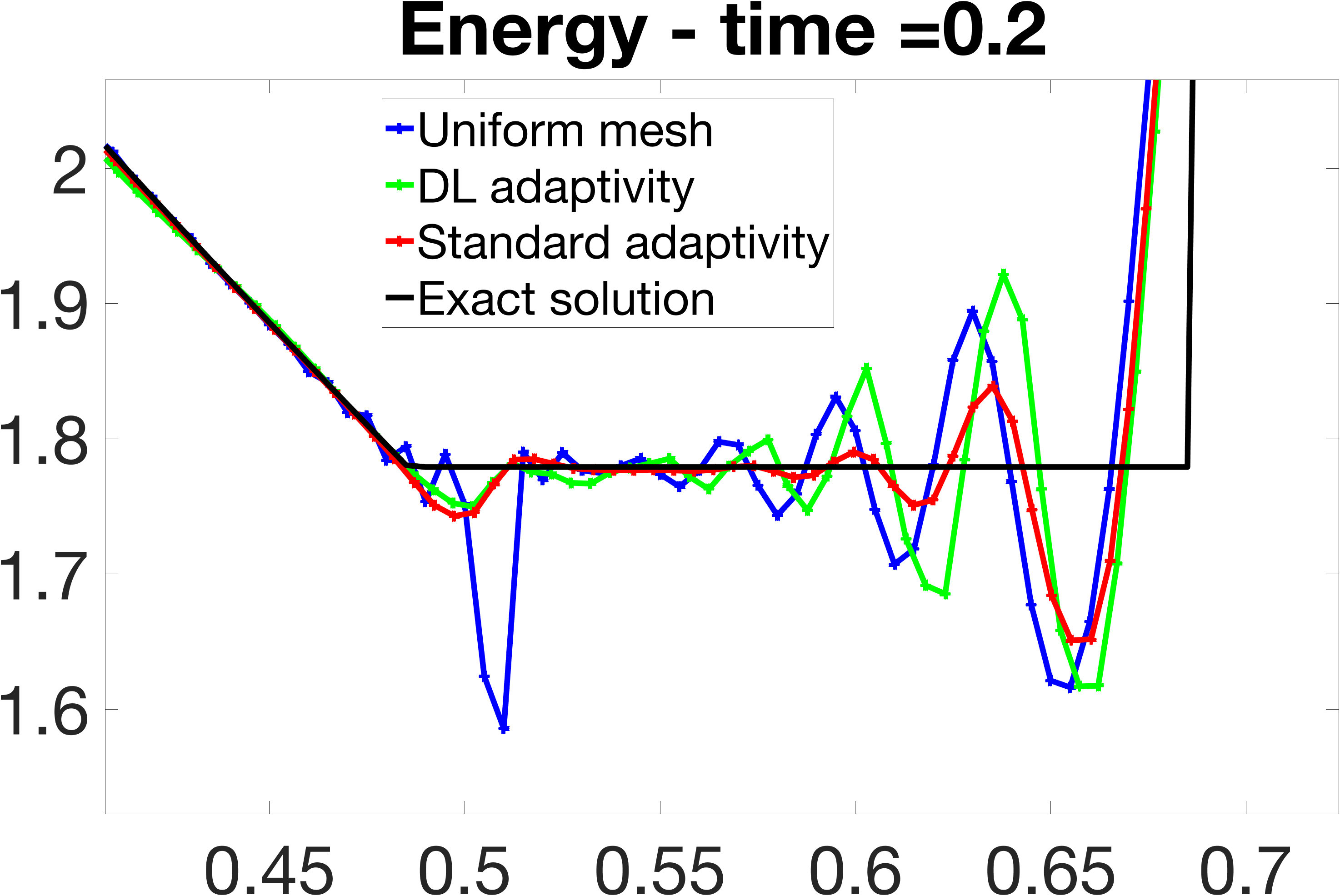} \\ 
\vspace{1cm}
\includegraphics[width=0.45\textwidth]{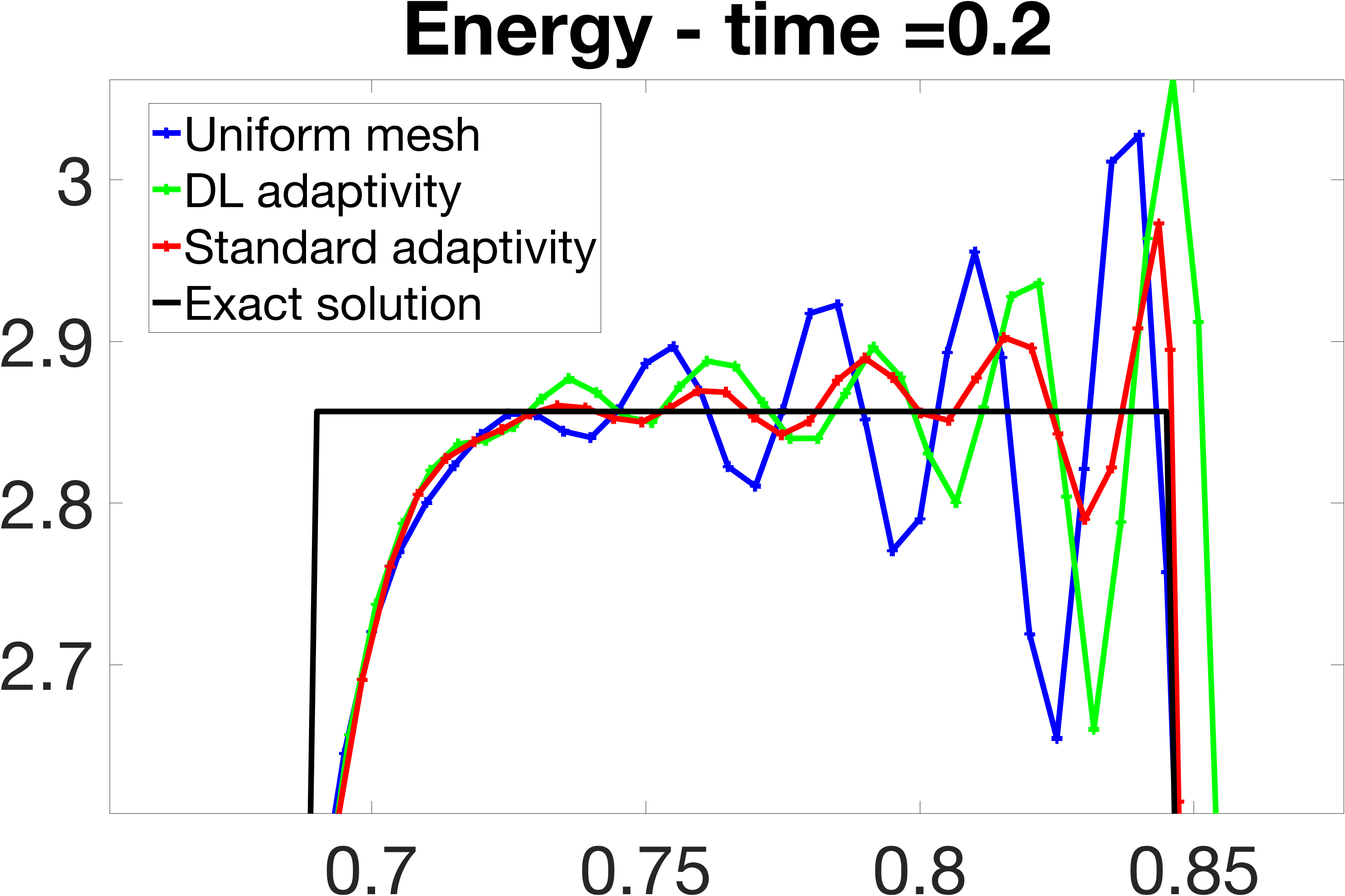}  
\caption{Energy profile of 1D Sod test case at $t=0.2$. Lax Wendroff scheme used for space discretization and time stepping. Entire profile shown at the top left, zoom-in in the $x$ range [0.4, 0.7] shown at the top right, and zoom-in in the $x$ range [0.7, 0.85] shown at the bottom.}
\label{euler_lw_energy}
\end{figure}

\begin{figure}[H]
\centering
\includegraphics[width=0.45\textwidth]{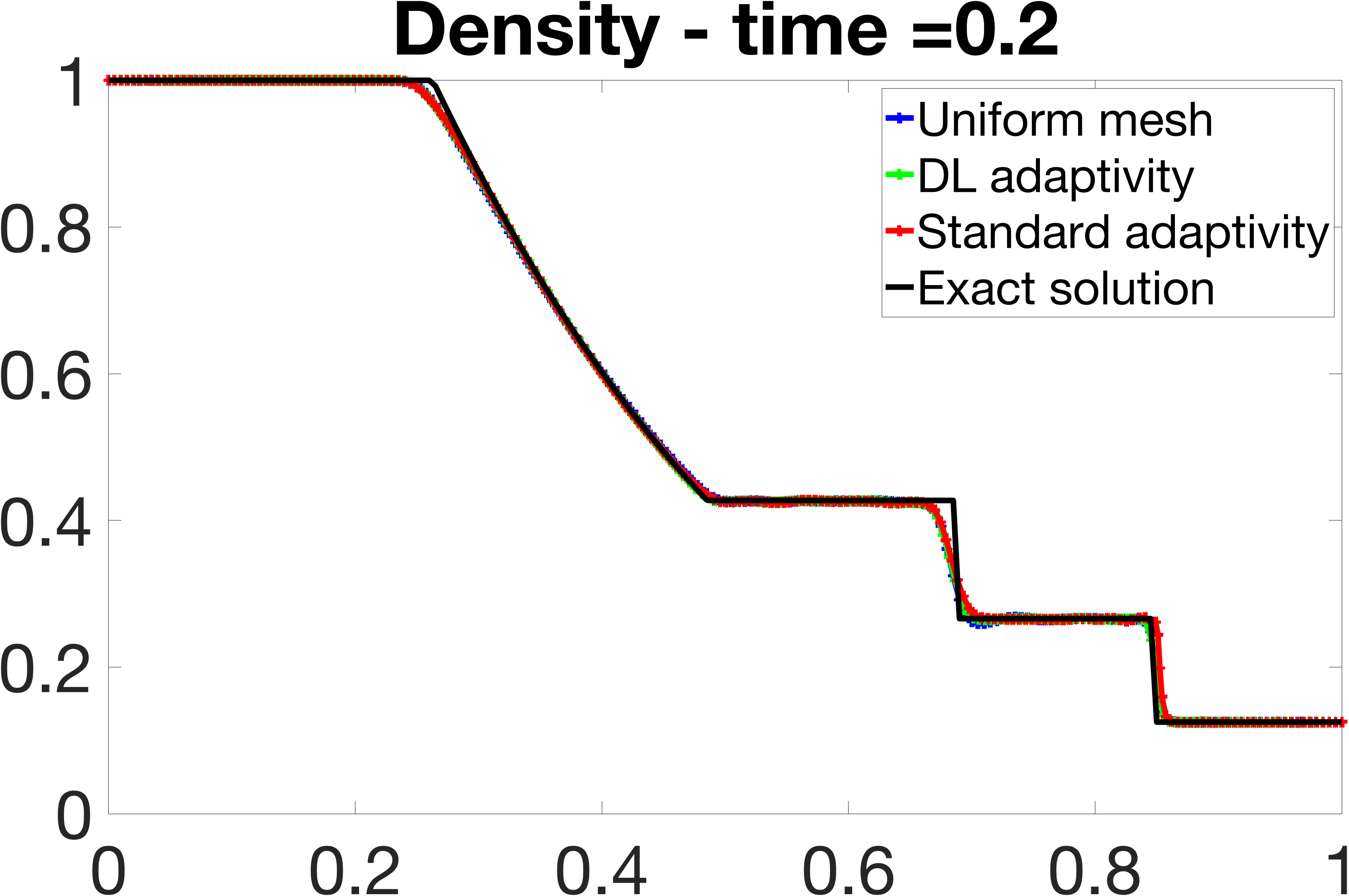}   
\hspace{1cm}
\includegraphics[width=0.45\textwidth]{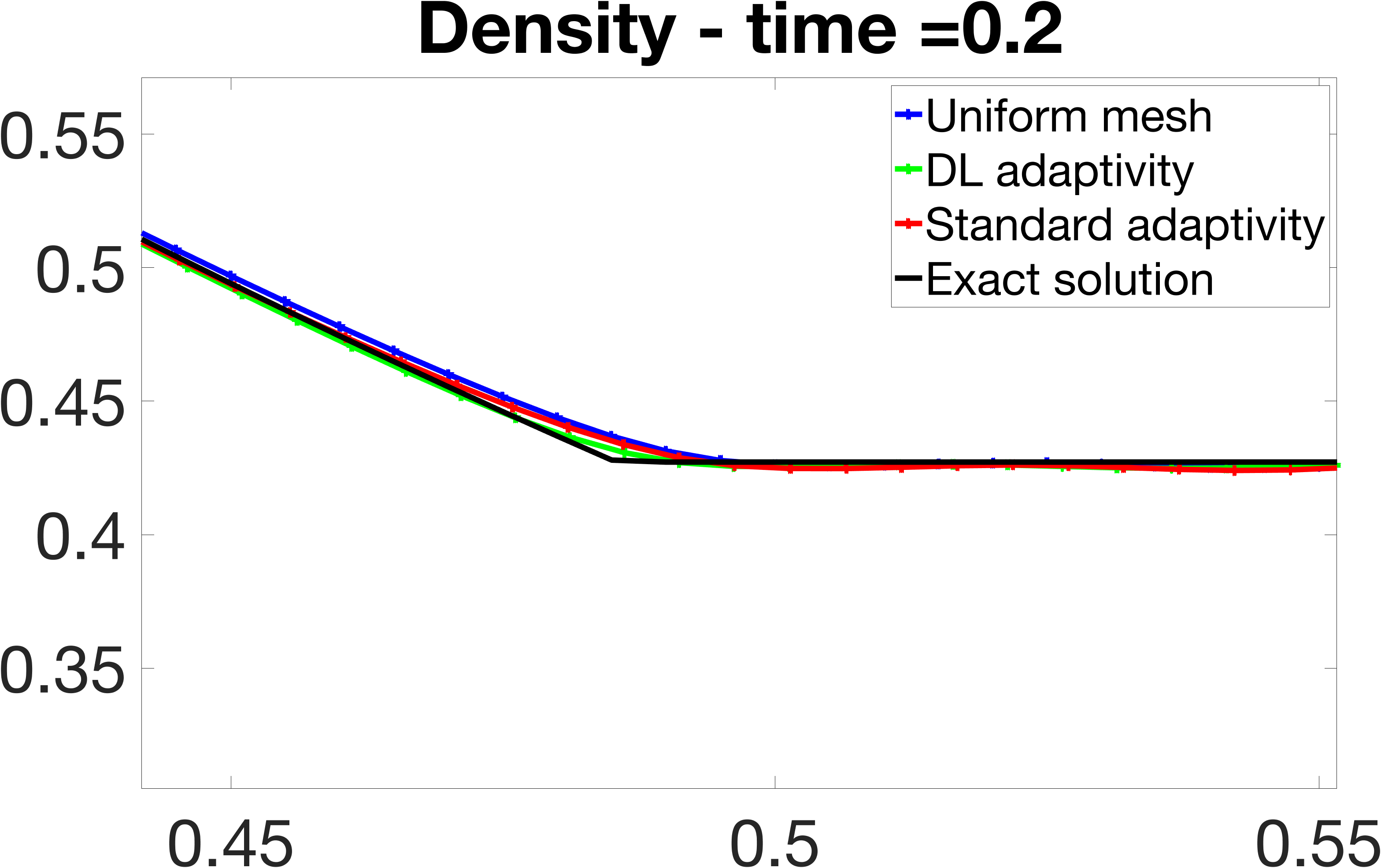}  \\
\vspace{1cm}
\includegraphics[width=0.45\textwidth]{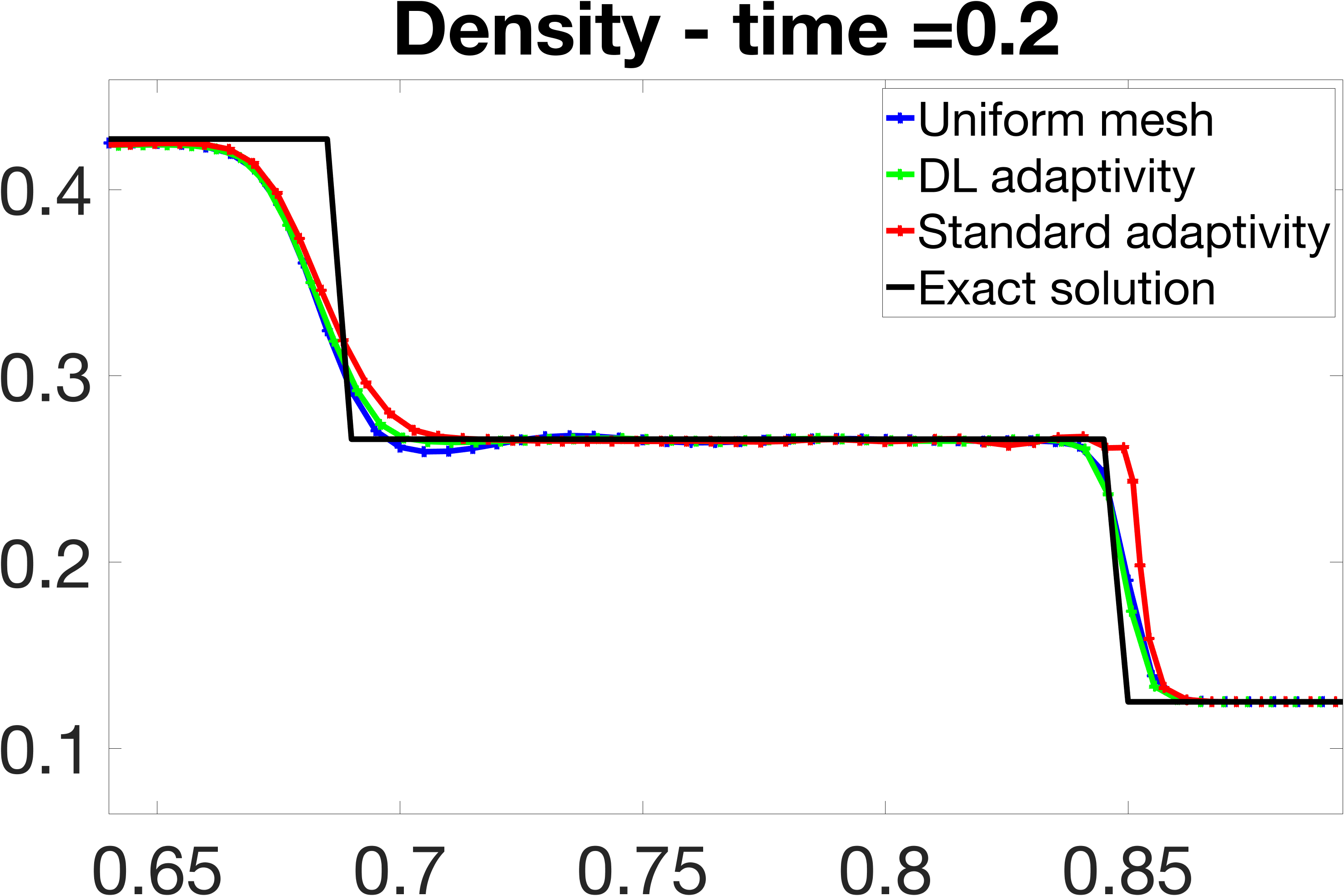}
\caption{Density profile of 1D Sod test case at $t=0.2$. WENO5 scheme used for space discretization and RK-3 used for time stepping. Entire profile shown at the top left, zoom-in in the $x$ range [0.4, 0.55] shown at the top right and zoom-in in the $x$ range [0.65, 085] shown at the bottom.}
\label{euler1_weno5_density}
\end{figure}

\begin{figure}[H]
\centering
\includegraphics[width=0.45\textwidth]{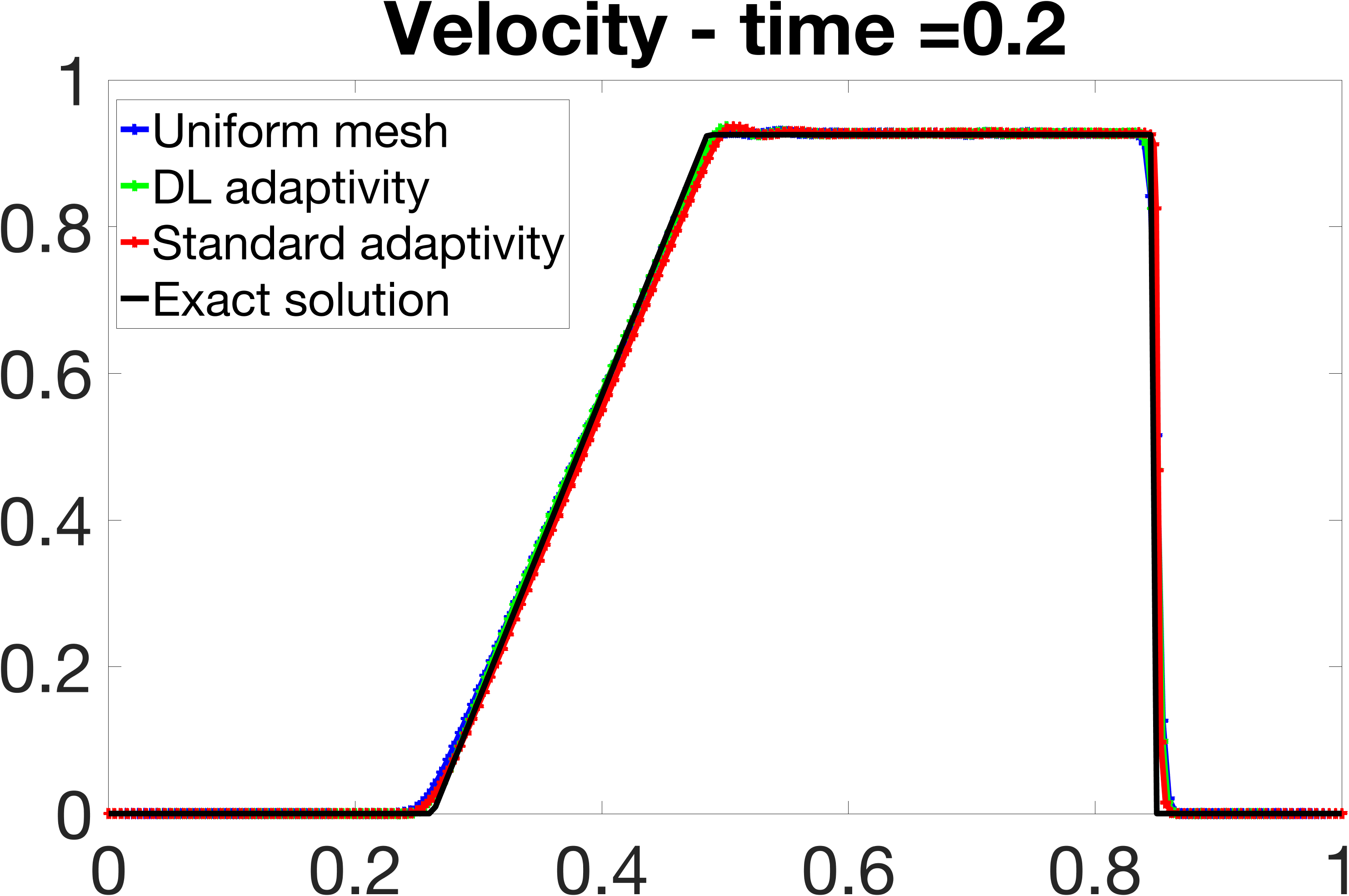}  
\hspace{1cm}
\includegraphics[width=0.47\textwidth]{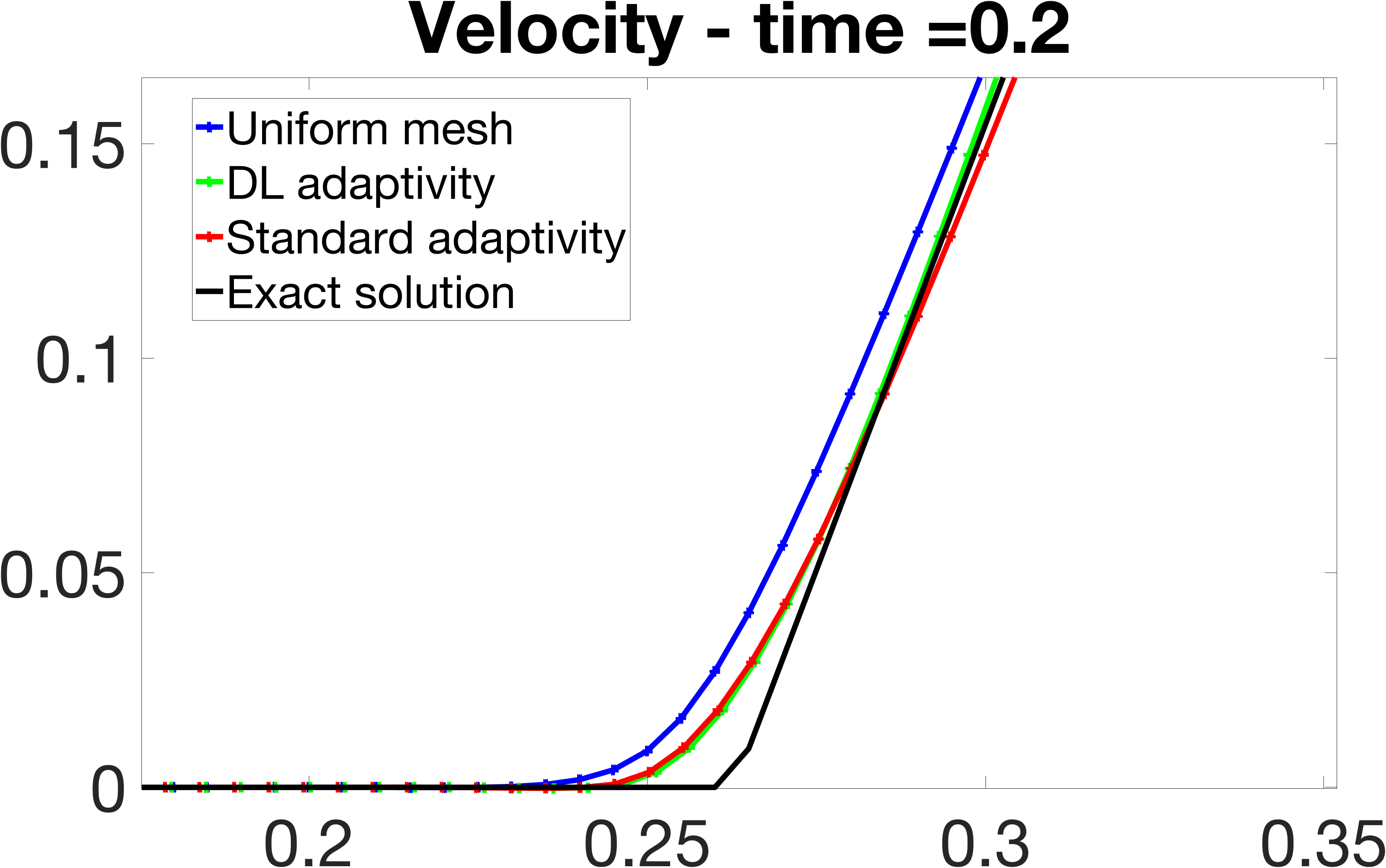}\\
\vspace{1cm}
\includegraphics[width=0.45\textwidth]{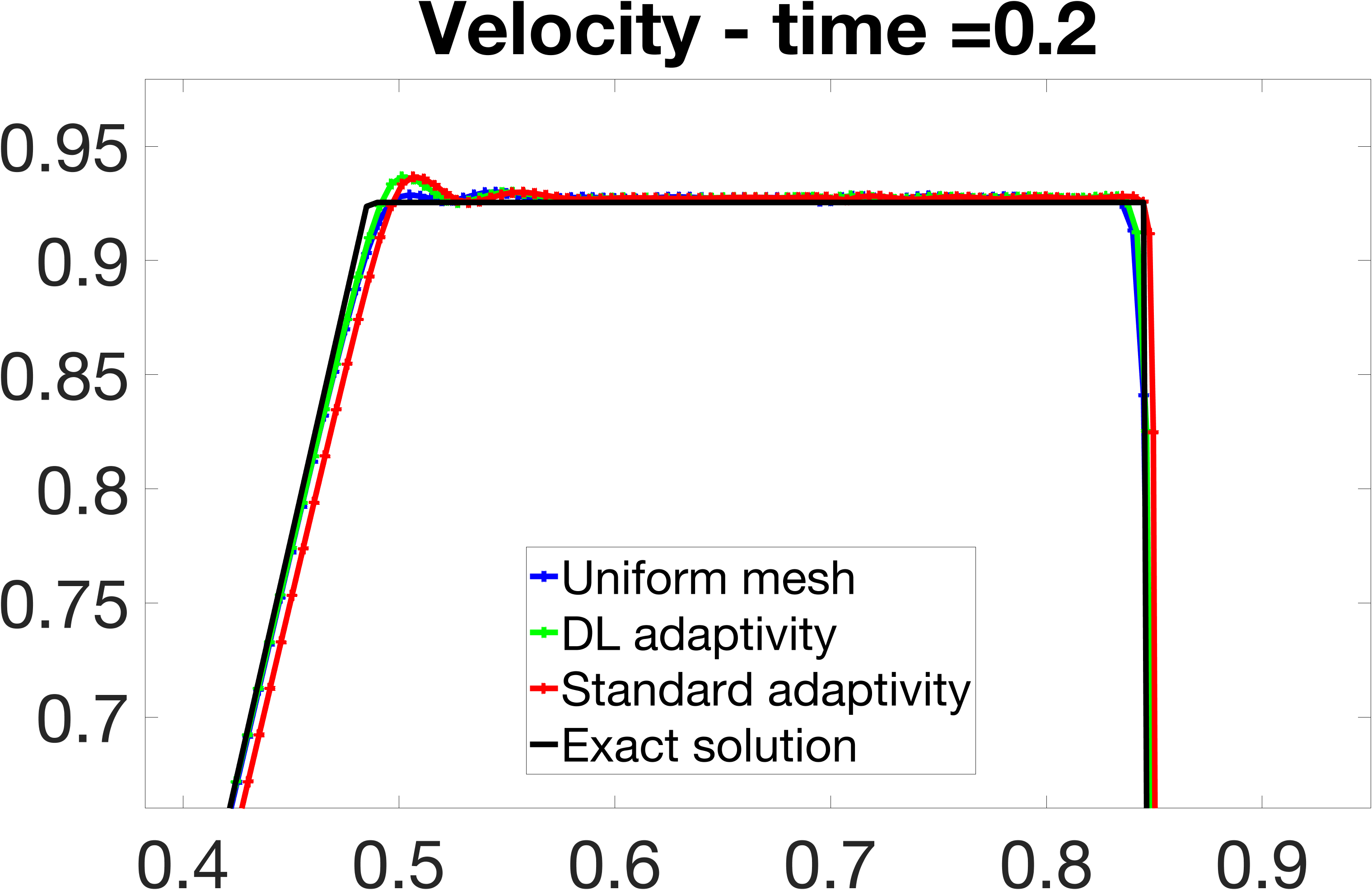}
\hspace{1cm}
\includegraphics[width=0.45\textwidth]{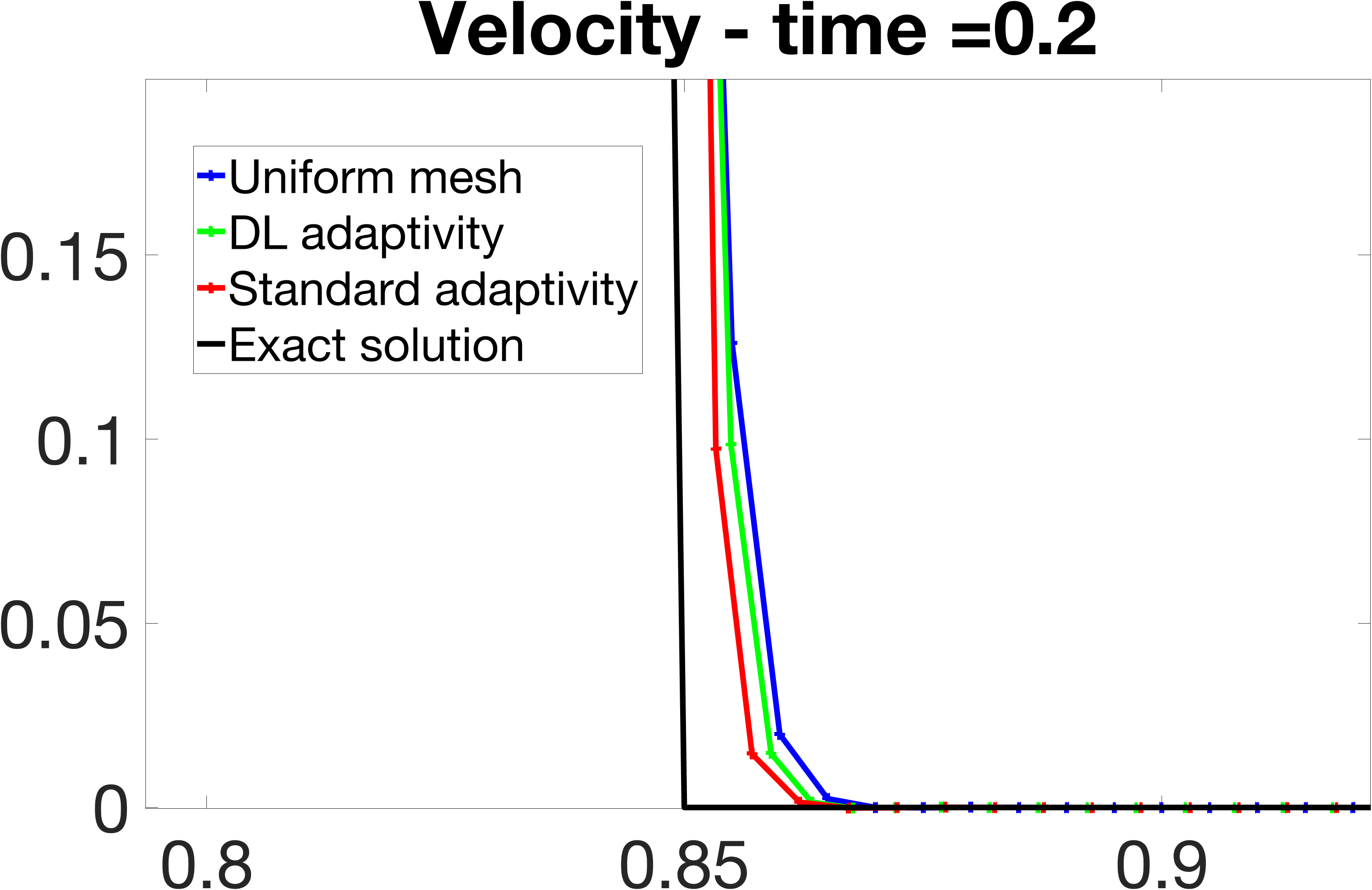}
\caption{Velocity profile 1D Sod test case at $t=0.2$. WENO5 scheme used for space discretization and RK-3 used for time stepping. Entire profile shown at the top left, zoom-in in the $x$ range [0.0, 0.35] shown at the top right, zoom-in in the $x$ range [0.4, 08] shown at the bottom left, and zoom-in in the $x$ range [0.8, 09] shown at the bottom right.}
\label{euler1_weno5_velocity}
\end{figure}

\begin{figure}[H]
\centering
\includegraphics[width=0.45\textwidth]{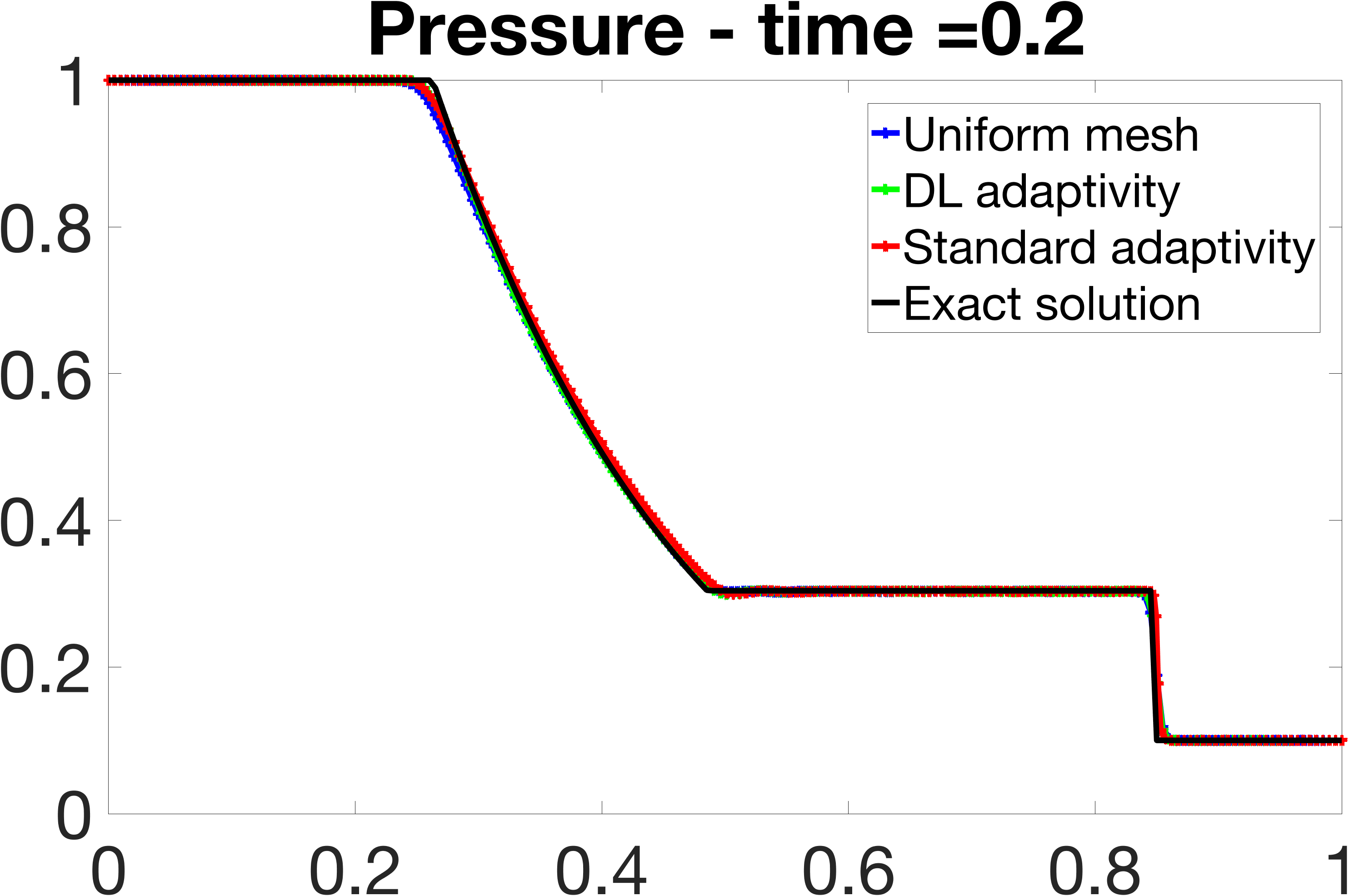}  
\hspace{1cm}
\includegraphics[width=0.45\textwidth]{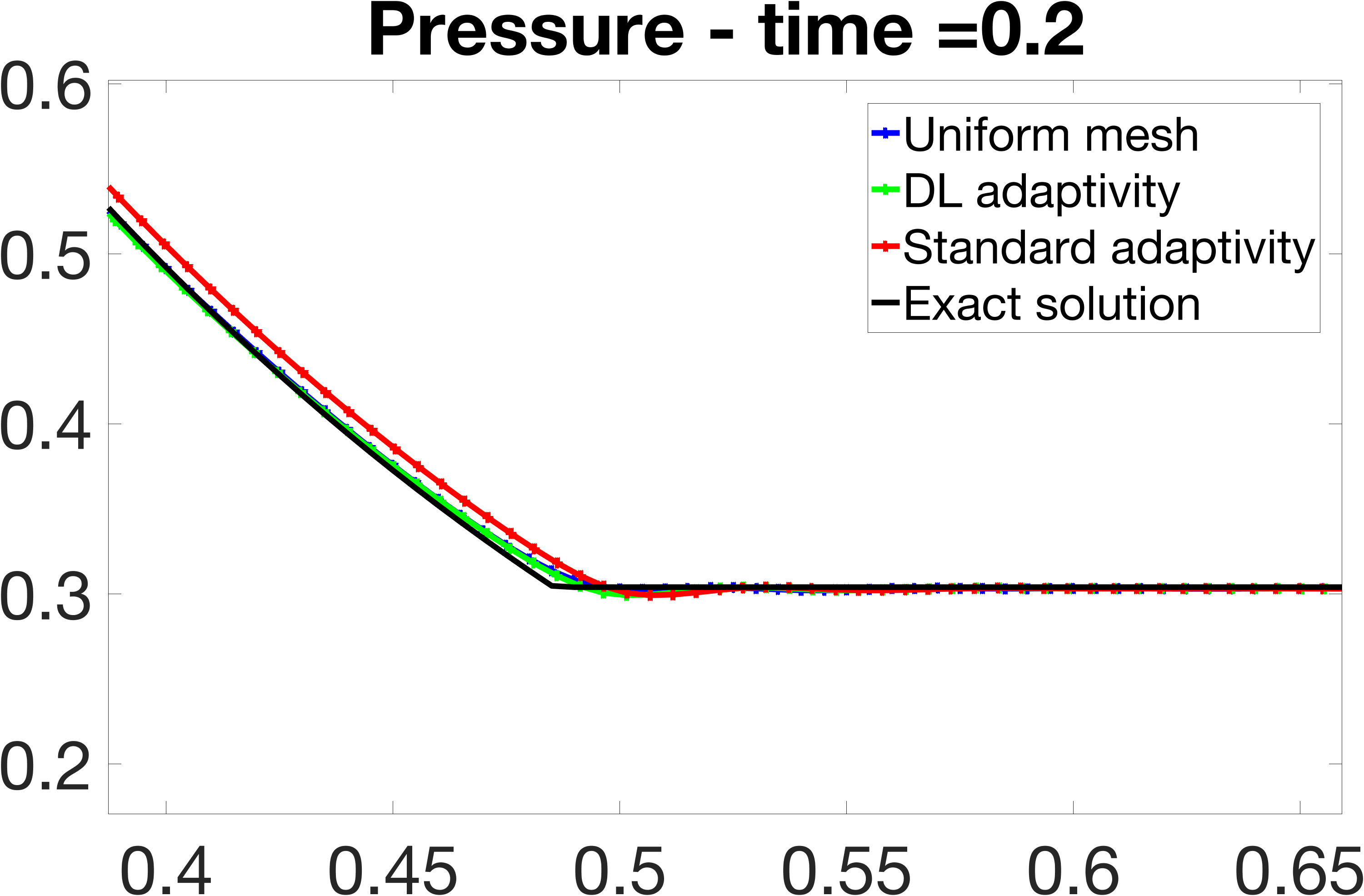}   \\
\vspace{1cm}
\includegraphics[width=0.45\textwidth]{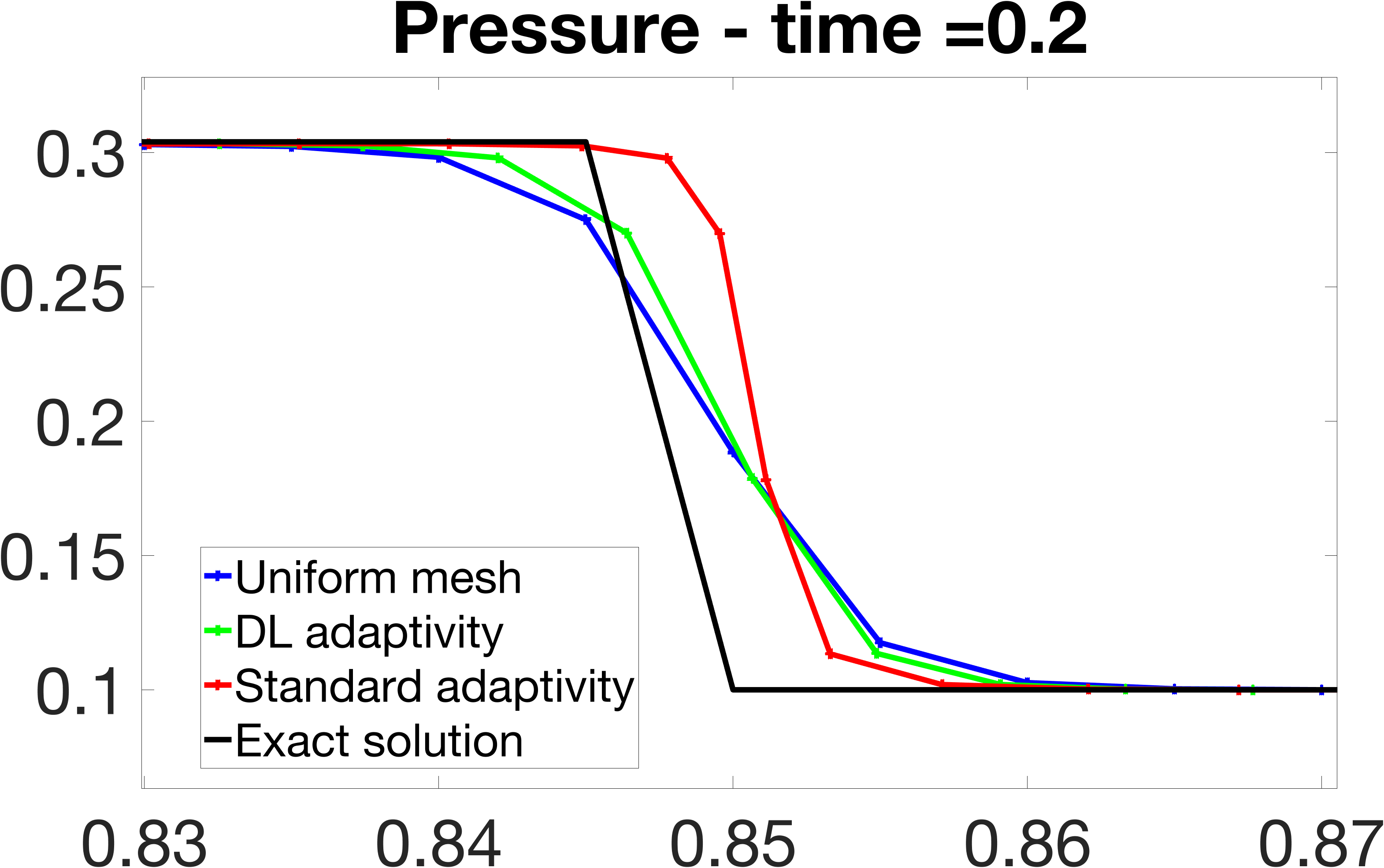} \\
\caption{Pressure profile of 1D Sod test case at $t=0.2$. WENO5 scheme used for space discretization and RK-3 used for time stepping. Entire profile shown at the top left, zoom-in in the $x$ range [0.4, 0.65] shown at the top right, and zoom-in in the $x$ range [0.83, 0.87] shown at the bottom. }
\label{euler1_weno5_pressure}
\end{figure}

\begin{figure}[H]
\centering
\includegraphics[width=0.45\textwidth]{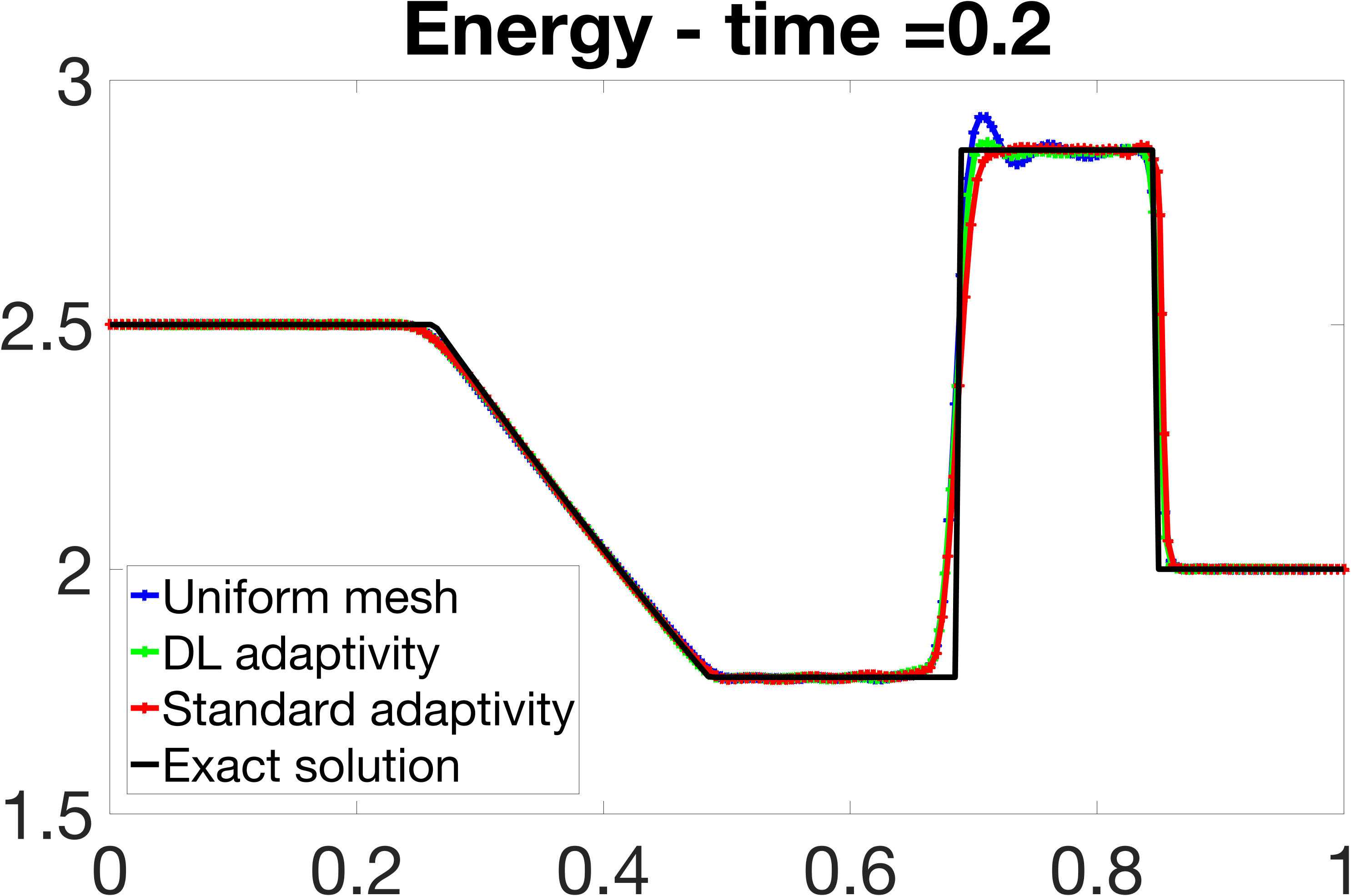}   
\hspace{1cm}
\includegraphics[width=0.45\textwidth]{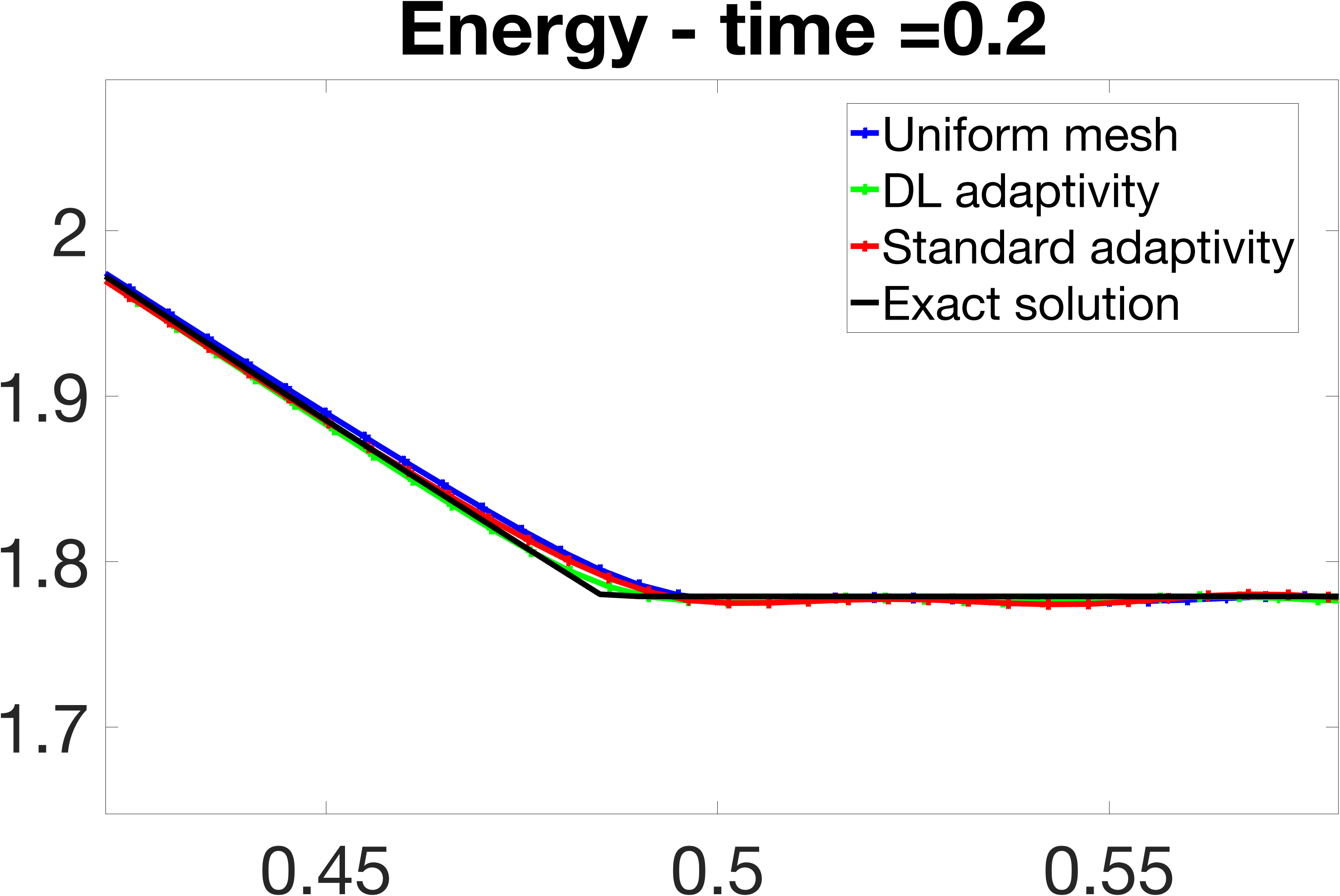} \\
\vspace{1cm}
\includegraphics[width=0.45\textwidth]{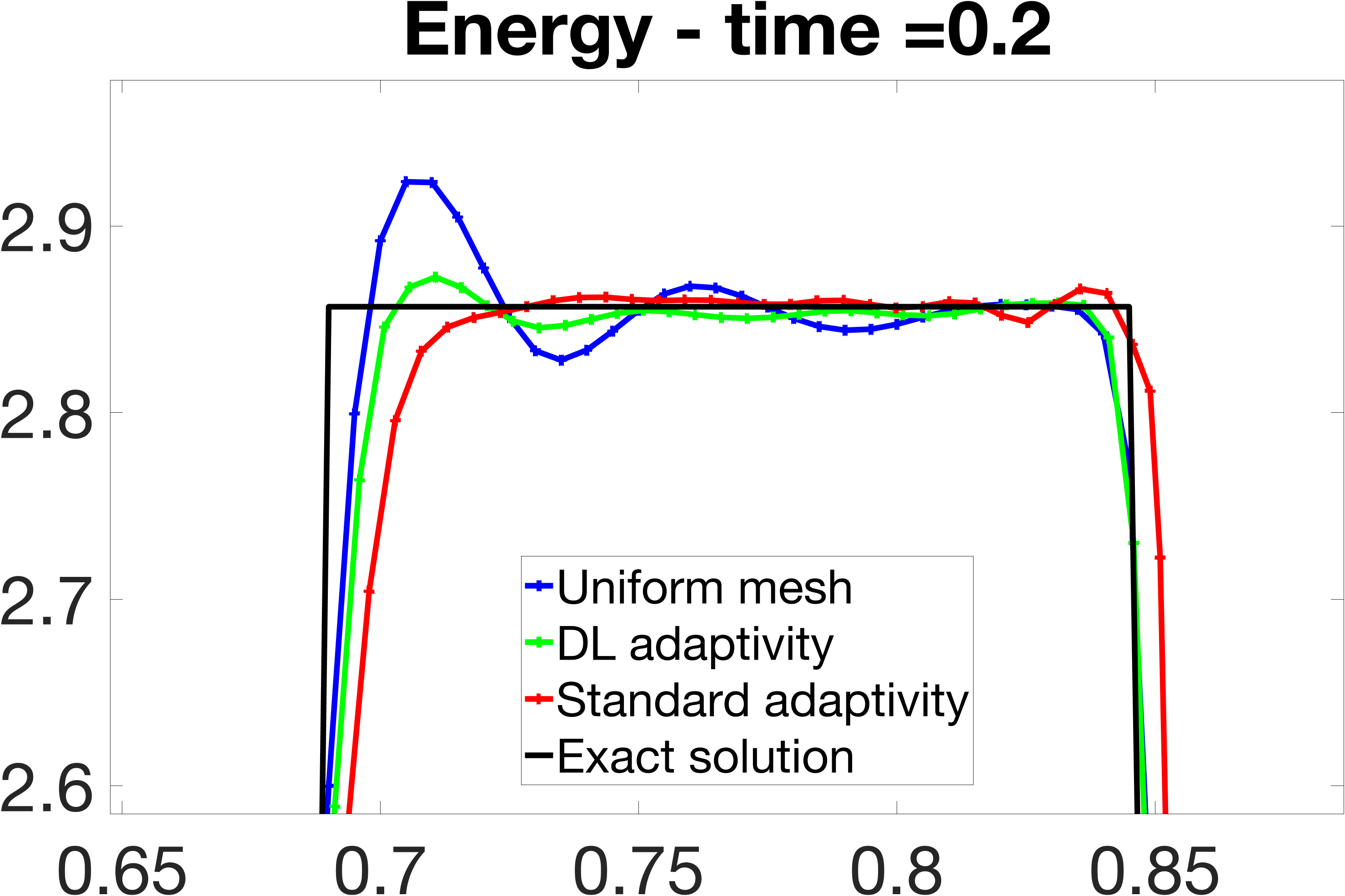} 
\hspace{1cm}
\includegraphics[width=0.45\textwidth]{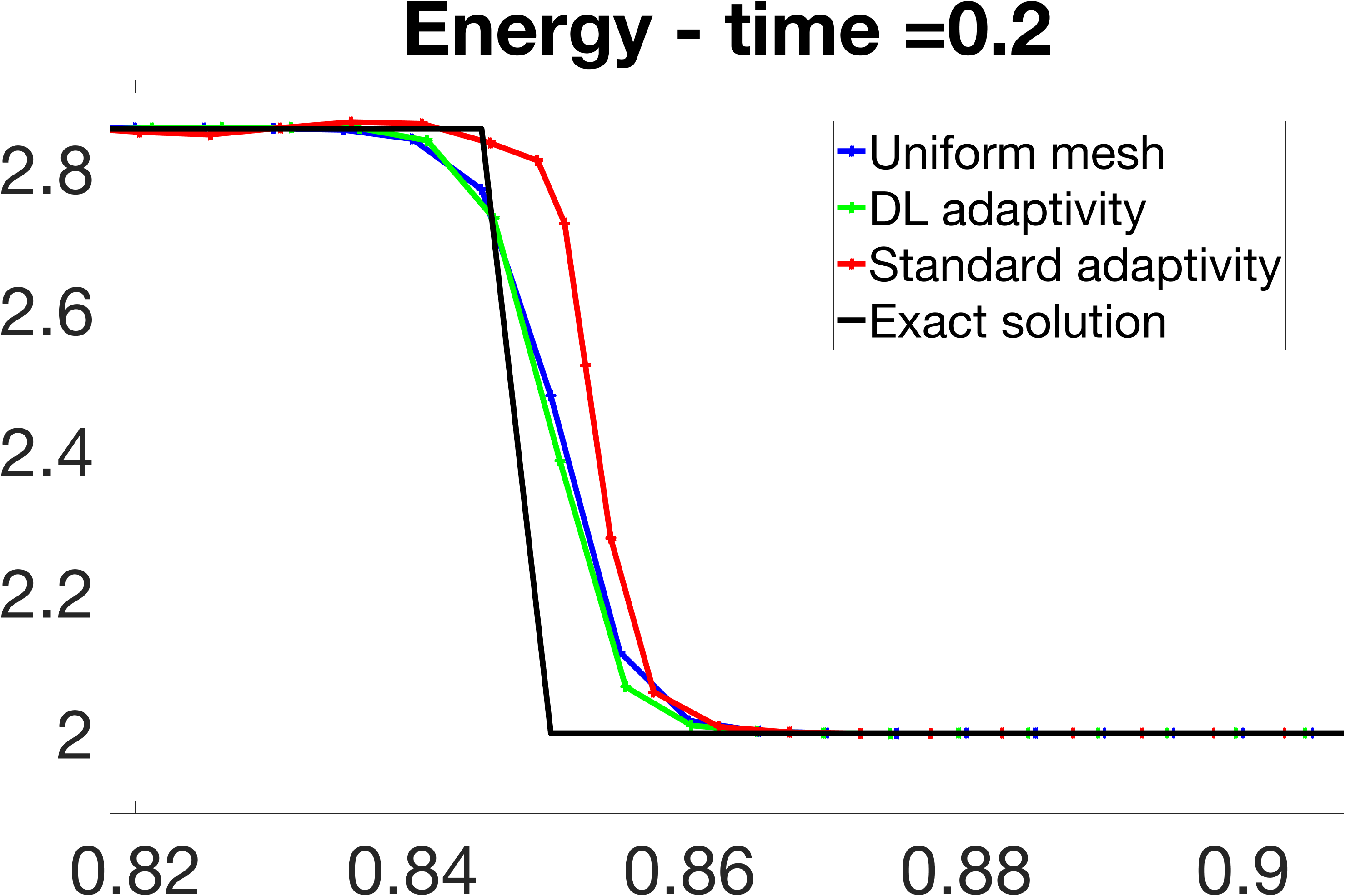} 
\caption{Energy profile of 1D Sod test case at $t=0.2$. WENO5 scheme used for space discretization and RK-3 used for time stepping. Entire profile shown at the top left, zoom-in in the $x$ range [0.4, 0.55] shown at the top right, zoom-in in the $x$ range [0.65, 0.85] shown at the bottom left, and zoom-in in the range [0.82, 0.9].}
\label{euler2_weno5_energy}
\end{figure}

\begin{figure}[H]
\centering
\includegraphics[width=0.50\textwidth]{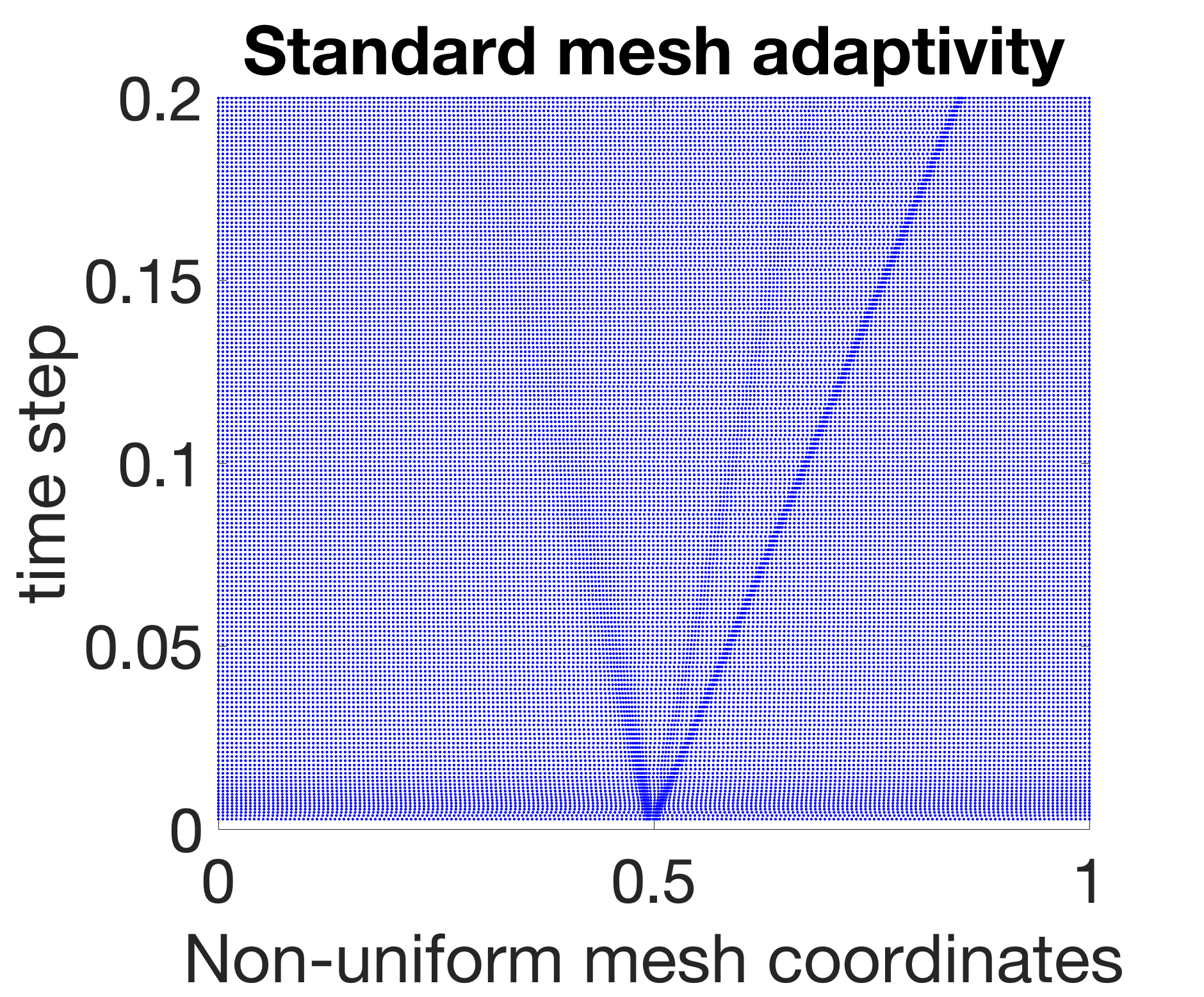}   
\includegraphics[width=0.49\textwidth]{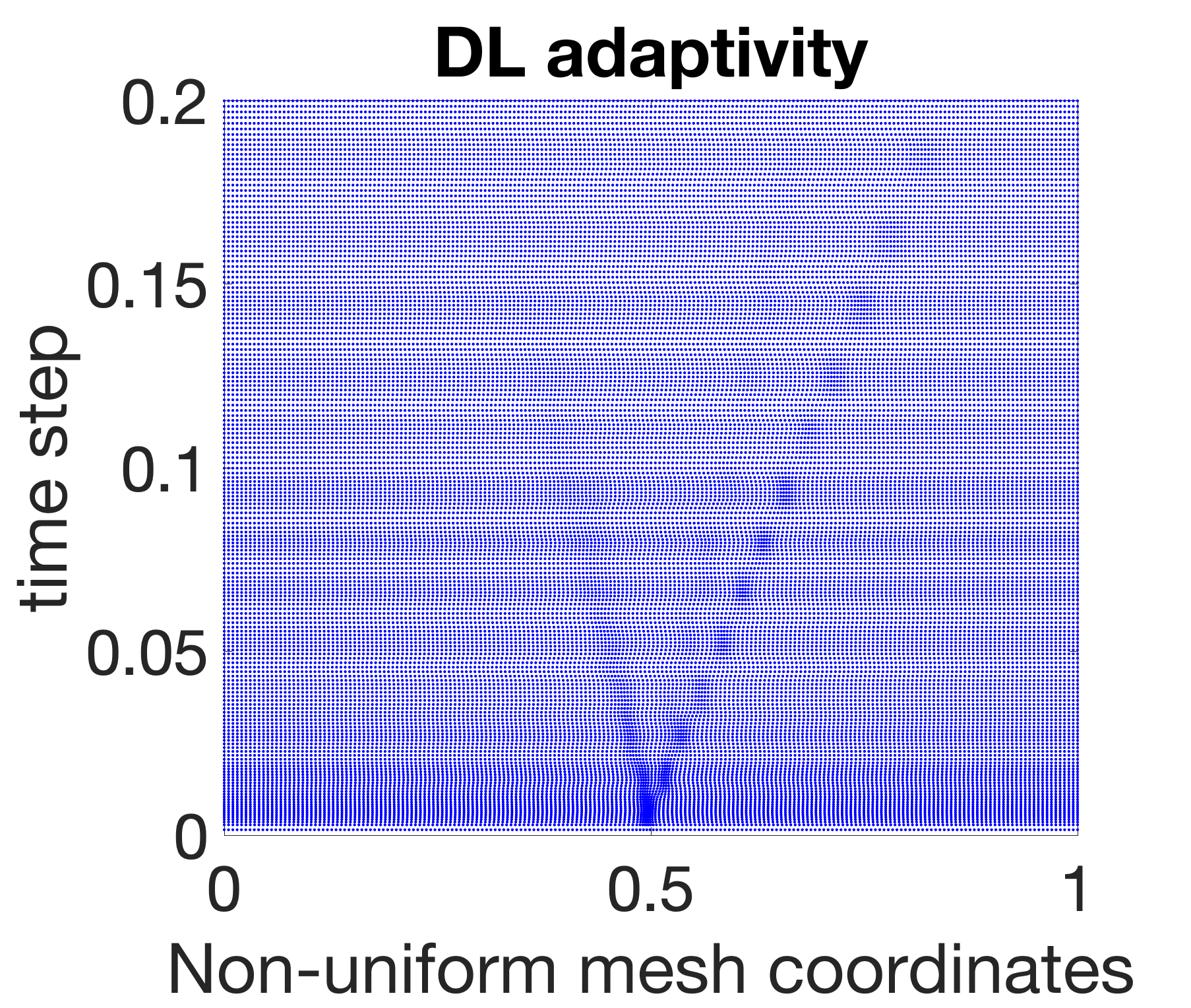}   
\caption{History of non-uniform adapted meshes for 1D Sod test case using standard adaptive zoning (left) and deep learning adaptive zoning (right) with the WENO5 scheme for space discretization and the RK-3 scheme used for time stepping.}
\label{mesh_history}
\end{figure}

The $L^1$ and $L^2$ norms of exact solution, numerical approximations, and numerical errors for density, velocity, pressure, and energy using Lax-Wendroff are provided in Table \ref{tab:sod_LW_norms}. Standard adaptive zoning and DL adaptive zoning are comparable in accuracy and both outperform the non-adaptive approach. Moreover, the computational time collected in Table \ref{tab:sod_lw_time} shows that DL is faster than the standard adaptive technique.  
The profiles of density, velocity, pressure, and energy using Lax-Wendroff are shown in Figures \ref{euler_lw_density}, \ref{euler_lw_velocity}, \ref{euler_lw_pressure}, and \ref{euler_lw_energy}. Both standard adaptive zoning and DL adaptive zoning dampen the numerical oscillations in the numerical reconstruction of each physical quantity with respect to the calculations performed on a uniform mesh. 
The $L^1$ and $L^2$ norms of exact solution, numerical approximations, and numerical errors for density, velocity, pressure, and energy using WENO5 for space discretization and RK-3 for time stepping are provided in Table \ref{tab:sod_norms}. The DL adaptivity allows the physics PDE solver to attain an accuracy comparable to what is obtained with the standard adaptivity on all the four physical quantities of interest. 
The profiles of density, velocity, pressure, and energy profiles using the WENO5 for space discretization and RK-3 for time stepping are shown in Figures \ref{euler1_weno5_density}, \ref{euler1_weno5_velocity}, \ref{euler1_weno5_pressure}, and \ref{euler2_weno5_energy}, respectively. Both standard adaptive zoning and DL adaptive zoning still dampen the numerical oscillations for each physical quantity with respect to the calculations performed on a uniform mesh. 
DL adaptive zoning still runs faster than standard adaptive zoning per time step as shown in Table \ref{tab:sod_weno5_time}. 
The evolution of the meshes across consecutive time steps for standard adaptive zoning and DL adaptive zoning are shown in Figure \ref{mesh_history}. The change of the mesh across consecutive time steps is smoother for the standard adaptive zoning, whereas the mesh updates produced by DL adaptivity are more irregular. Despite the irregularities of the adapted mesh computed with the DL approach, the numerical solution does not suffer any loss of resolution close to the shock.   

\begin{table}[ht]
    \centering
    \begin{tabular}{|c|c|c|c|c|}
    \hline
    \multicolumn{5}{|c|}{\textbf{Density}}\\
    \hline
    & \textbf{Exact solution}& \textbf{Uniform mesh} & \textbf{Standard adaptivity} & \textbf{DL adaptivity}\\
    \hline
    $L^2$-norm of the solution & 0.6503 & 0.6525 & 0.6522 & 0.6533\\  
    \hline
    $L^2$-norm of the relative error & - & 0.0461 & 0.0330 & 0.0357\\  
    \hline
    $L^1$-norm of the solution & 0.5625 & 0.5647 & 0.5662 & 0.5680\\ 
    \hline
    $L^1$-norm of the relative error & - & 0.0263 & 0.0136 & 0.0171 \\      
    \hline
    \multicolumn{5}{|c|}{\textbf{Velocity}}\\
    \hline
    & \textbf{Exact solution}& \textbf{Uniform mesh} & \textbf{Standard adaptivity} & \textbf{DL adaptivity}\\
    \hline
    $L^2$-norm of the solution & 0.6124 & 0.6183 & 0.6240 & 0.6257\\  
    \hline
    $L^2$-norm of the relative error & - & 0.2301 & 0.2100 & 0.2318\\  
    \hline
    $L^1$-norm of the solution & 0.4395 & 0.4420 & 0.4558 & 0.4575\\ 
    \hline
    $L^1$-norm of the relative error & - & 0.1068 & 0.0559 & 0.0722 \\      
    \hline
    \multicolumn{5}{|c|}{\textbf{Internal energy}}\\
    \hline
    & \textbf{Exact solution}& \textbf{Uniform mesh} & \textbf{Standard adaptivity} & \textbf{DL adaptivity}\\
    \hline
    $L^2$-norm of the solution & 2.2844 & 2.2817 & 2.3056 & 2.3024\\  
    \hline
    relative $L^2$-norm of the error & - & 0.0280 & 0.0295 & 0.0259\\  
    \hline
    $L^1$-norm of the solution & 2.2844 & 2.2508 & 2.2734 & 2.2712\\ 
    \hline
    relative $L^1$-norm of the error & - & 0.0280 & 0.0158 & 0.0206 \\      
    \hline
    \multicolumn{5}{|c|}{\textbf{Pressure}}\\
    \hline
    & \textbf{Exact solution}& \textbf{Uniform mesh} & \textbf{Standard adaptivity} & \textbf{DL adaptivity}\\
    \hline
    $L^2$-norm of the solution & 0.6210 & 0.6232 & 0.6233 & 0.6240\\  
    \hline
    relative $L^2$-norm of the error & - & 0.0616 & 0.0459 & 0.0509\\  
    \hline
    $L^1$-norm of the solution & 0.5210 & 0.5228 & 0.5261 & 0.5273\\ 
    \hline
    relative $L^1$-norm of the error & - & 0.0304 & 0.0149 & 0.0197 \\      
    \hline
    \end{tabular}
    \caption{1D Sod test case. Relative error for density, velocity, internal energy, and pressure in the $L^2$-norm and $L^1$-norm using the Lax-Wendroff scheme for space and time stepping.}
    \label{tab:sod_LW_norms}
\end{table}

\begin{table}[ht]
    \centering
    \begin{tabular}{|c|c|}
    \hline
    \textbf{Mesh Type} & \textbf{Wall-clock Time (s)}\\
    \hline
    Uniform mesh & 0.62\\
    \hline
    Standard adaptivity & 62.58\\
    \hline
    DL adaptivity & 748.25 - $(0.1476 \times 4,980)$ = 13.20\\
    \hline
    \end{tabular}
    \caption{1D Sod test case. Wall-clock computational time in seconds using the Lax-Wendroff scheme on a uniform mesh, non-uniform mesh computed with standard adaptive zoning, and non-uniform mesh computed with DL adaptive zoning.}
    \label{tab:sod_lw_time}
\end{table}

\begin{table}[ht]
    \centering
    \begin{tabular}{|c|c|c|c|c|}
    \hline
    \multicolumn{5}{|c|}{\textbf{Density}}\\
    \hline
    & \textbf{Exact solution}& \textbf{Uniform mesh} & \textbf{Standard adaptivity} & \textbf{DL adaptivity}\\
    \hline
    $L^2$-norm of the solution & 0.6503 & 0.6517 & 0.6517 & 0.6518\\  
    \hline
    $L^2$-norm of the relative error & - & 0.0198 & 0.0156 & 0.0150\\  
    \hline
    $L^1$-norm of the solution & 0.5625 & 0.5647 & 0.5640 & 0.5642\\ 
    \hline
    $L^1$-norm of the relative error & - & 0.0086 & 0.0073 & 0.0075 \\      
    \hline
    \multicolumn{5}{|c|}{\textbf{Velocity}}\\
    \hline
    & \textbf{Exact solution}& \textbf{Uniform mesh} & \textbf{Standard adaptivity} & \textbf{DL adaptivity}\\
    \hline
    $L^2$-norm of the solution & 0.6124 & 0.6118 & 0.6085 & 0.6094\\  
    \hline
    $L^2$-norm of the relative error & - & 0.1103 & 0.0380 & 0.0388\\  
    \hline
    $L^1$-norm of the solution & 0.4395 & 0.4421 & 0.4363 & 0.4374\\ 
    \hline
    $L^1$-norm of the relative error & - & 0.0271 & 0.0149 & 0.0153 \\      
    \hline
    \multicolumn{5}{|c|}{\textbf{Internal energy}}\\
    \hline
    & \textbf{Exact solution}& \textbf{Uniform mesh} & \textbf{Standard adaptivity} & \textbf{DL adaptivity}\\
    \hline
    $L^2$-norm of the solution & 2.2844 & 2.2913 & 2.2834 & 2.2837\\  
    \hline
    relative $L^2$-norm of the error & - & 0.0344 & 0.0295 & 0.0259\\  
    \hline
    $L^1$-norm of the solution & 2.2568 & 2.2613 & 2.2540 & 2.2543\\ 
    \hline
    relative $L^1$-norm of the error & - & 0.0085 & 0.0069 & 0.0066 \\      
    \hline
    \multicolumn{5}{|c|}{\textbf{Pressure}}\\
    \hline
    & \textbf{Exact solution}& \textbf{Uniform mesh} & \textbf{Standard adaptivity} & \textbf{DL adaptivity}\\
    \hline
    $L^2$-norm of the solution & 0.6210 & 0.6224 & 0.6224 & 0.6224\\  
    \hline
    relative $L^2$-norm of the error & - & 0.0231 & 0.0131 & 0.0132\\  
    \hline
    $L^1$-norm of the solution & 0.5210 & 0.5235 & 0.5224 & 0.5225\\ 
    \hline
    relative $L^1$-norm of the error & - & 0.0084 & 0.0068 & 0.0068 \\      
    \hline
    \end{tabular}
    \caption{1D Sod test case. Relative error for density, velocity, internal energy, and pressure in the $L^2$-norm and $L^1$-norm using the WENO5 scheme for space discretization and the RK-3 scheme for time stepping.}
    \label{tab:sod_norms}
\end{table}

\begin{table}[ht]
    \centering
    \begin{tabular}{|c|c|}
    \hline
    \textbf{Mesh Type} & \textbf{Wall-clock Time (s)}\\
    \hline
    Uniform mesh & 0.30\\
    \hline
    Standard adaptivity & 12.38\\
    \hline
    DL adaptivity & 276.05 - $(0.1476 \times 1,816)$ = 8.01\\
    \hline
    \end{tabular}
    \caption{1D Sod test case. Wall-clock computational time in seconds using the WENO5 scheme for space discretization and the RK-3 scheme used for time stepping on a uniform mesh, non-uniform mesh computed with standard adaptive zoning, and non-uniform mesh computed with DL adaptive zoning.}
    \label{tab:sod_weno5_time}
\end{table}

\clearpage

\subsubsection{Taylor–von Neumann–Sedov blast wave test case}
We present here the numerical results to assess the performance of DL adaptive zoning compared against standard adaptive zoning to model the propagation of a shock wave for the Taylor–von Neumann–Sedov problem, where the 1D Euler equations in the spatial domain [0, 0.6] are combined with the following initial conditions 
\begin{equation*}
\rho(x, t=0) = 1.0, \quad u(x, t=0) = 1.0, \quad e(x, t=0)  = 2.8049 \mathrm{e}-4
\end{equation*}
and non-homogeneous Dirichlet boundary conditions 
\begin{equation*}
\rho(0,t) = 1.0, \quad u(0,t) = 1.0, \quad e(0,t) = 2.8049 \mathrm{e}-4, \quad t\in [0,T),
\end{equation*}
\begin{equation*}
\rho(1,t) = 1.0, \quad u(1,t) = 1.0, \quad e(1,t) = 2.8049 \mathrm{e}-4, \quad t\in [0,T).
\end{equation*}
The initial and boundary values for the pressure are obtained from the state equations described in Section \ref{background}. 
We use WENO5 to discretize the domain [0,0.6] with 200 finite volume cells. The time integration is stopped at $t=1.0$.
The parameters used for the standard adaptivity are described in Table \ref{parameters_mmpde1}. The DL model is the same as for the square wave propagation test case. 

\begin{figure}[H]
\centering
\includegraphics[width=0.45\textwidth]{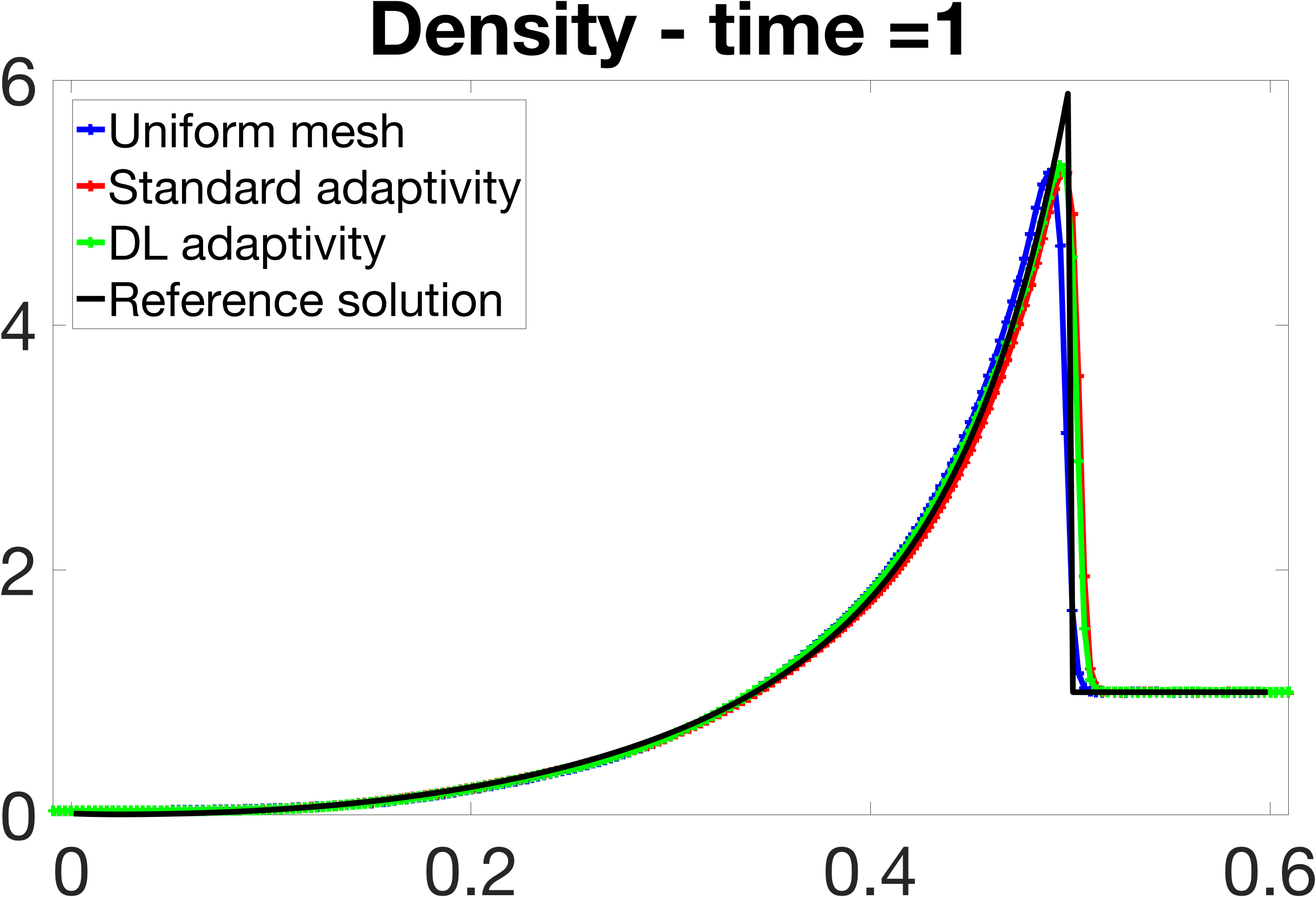}   
\hspace{1cm}
\includegraphics[width=0.45\textwidth]{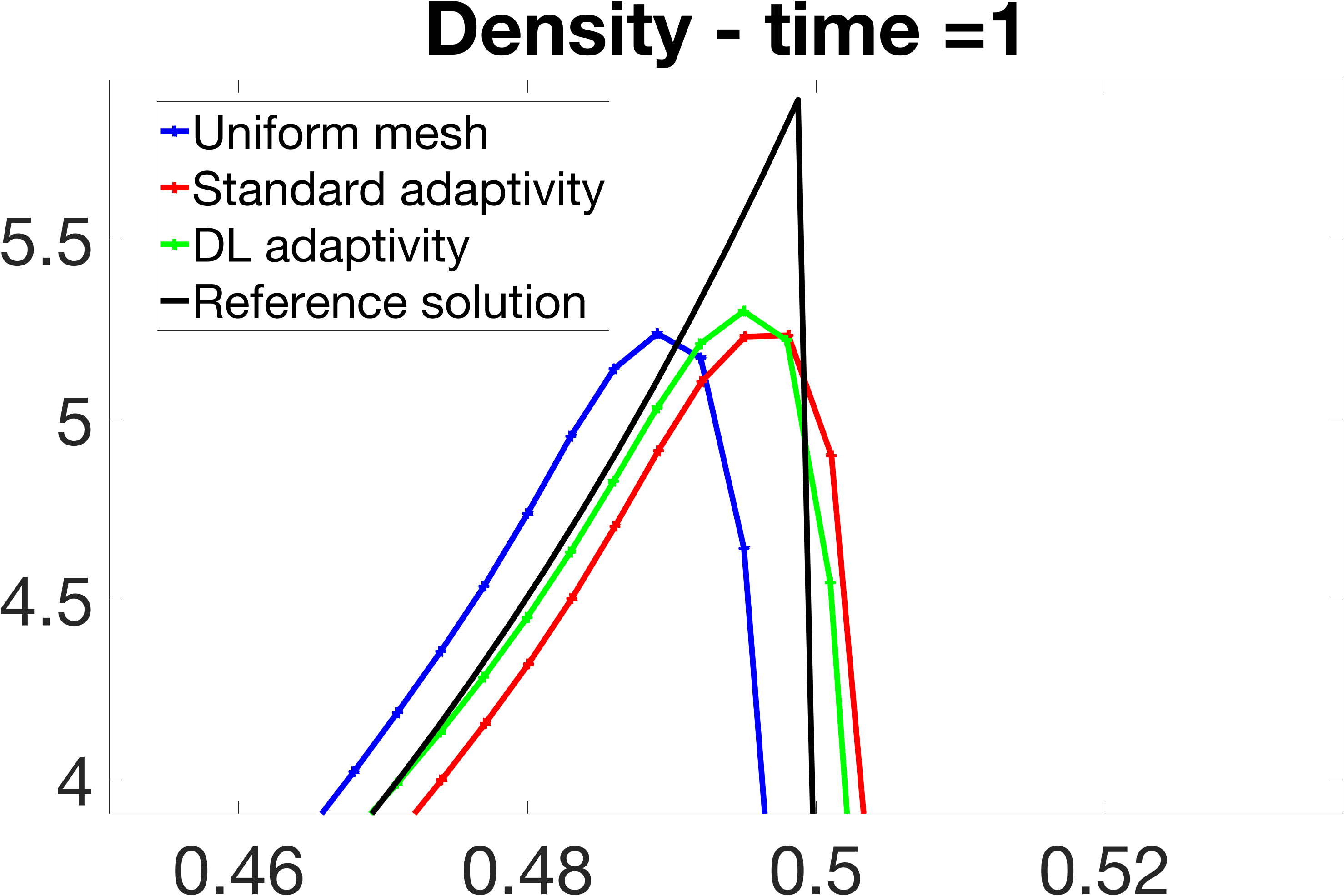} \\
\vspace{1cm}
\includegraphics[width=0.45\textwidth]{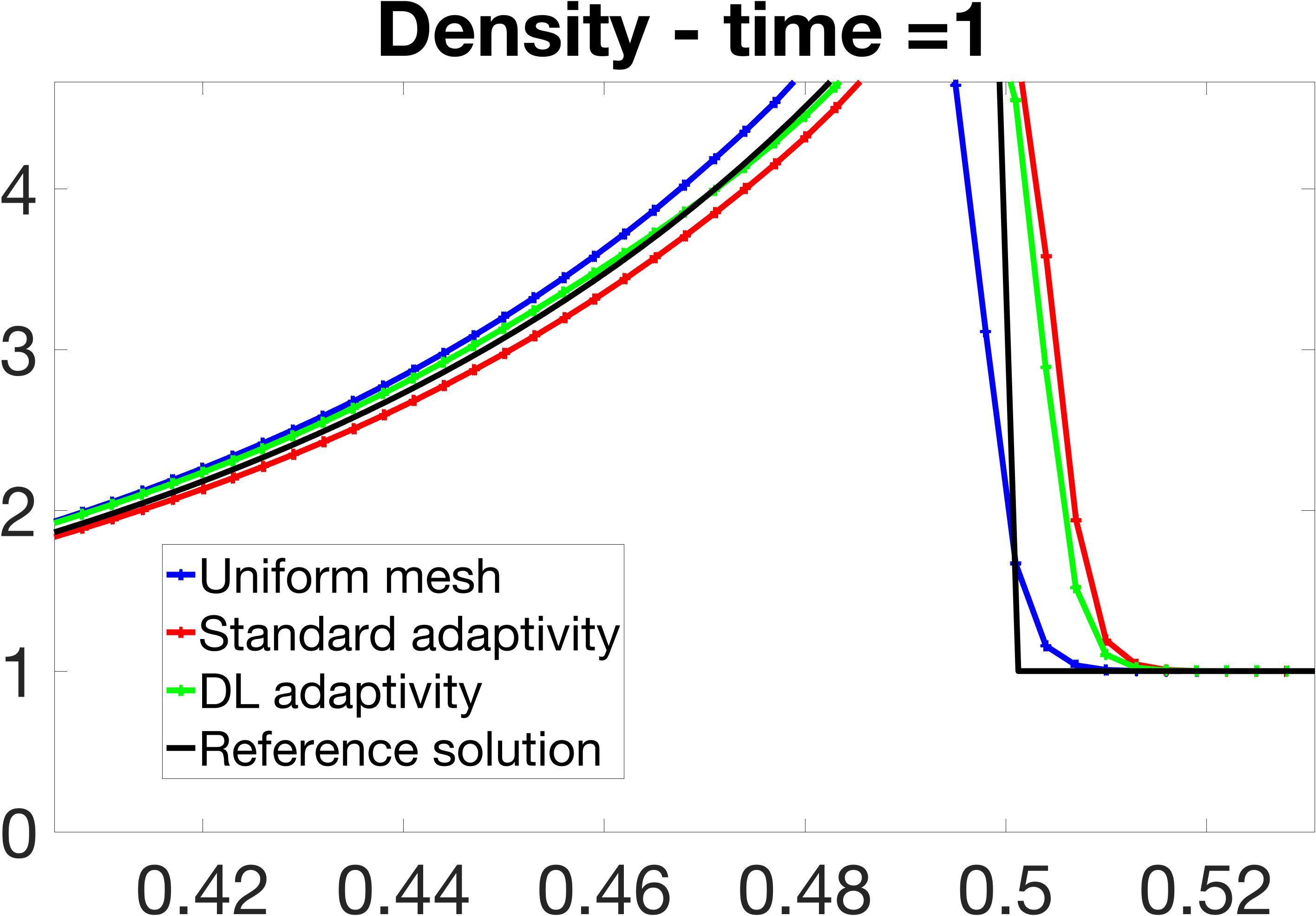} 
\caption{Density profile of 1D Taylor–von Neumann–Sedov test case at $t=1.0$. WENO5 scheme used for space discretization and RK-3 scheme used for time stepping. Entire profile shown at the top left, zoom-in in the $x$ range [0.4, 0.55] shown at the top right, zoom-in in the $x$ range [0.65, 0.85] shown at the bottom, and zoom-in in the $x$ range [0.82, 0.9].}
\label{sedov_weno5_density}
\end{figure}

\begin{figure}[H]
\centering
\includegraphics[width=0.45\textwidth]{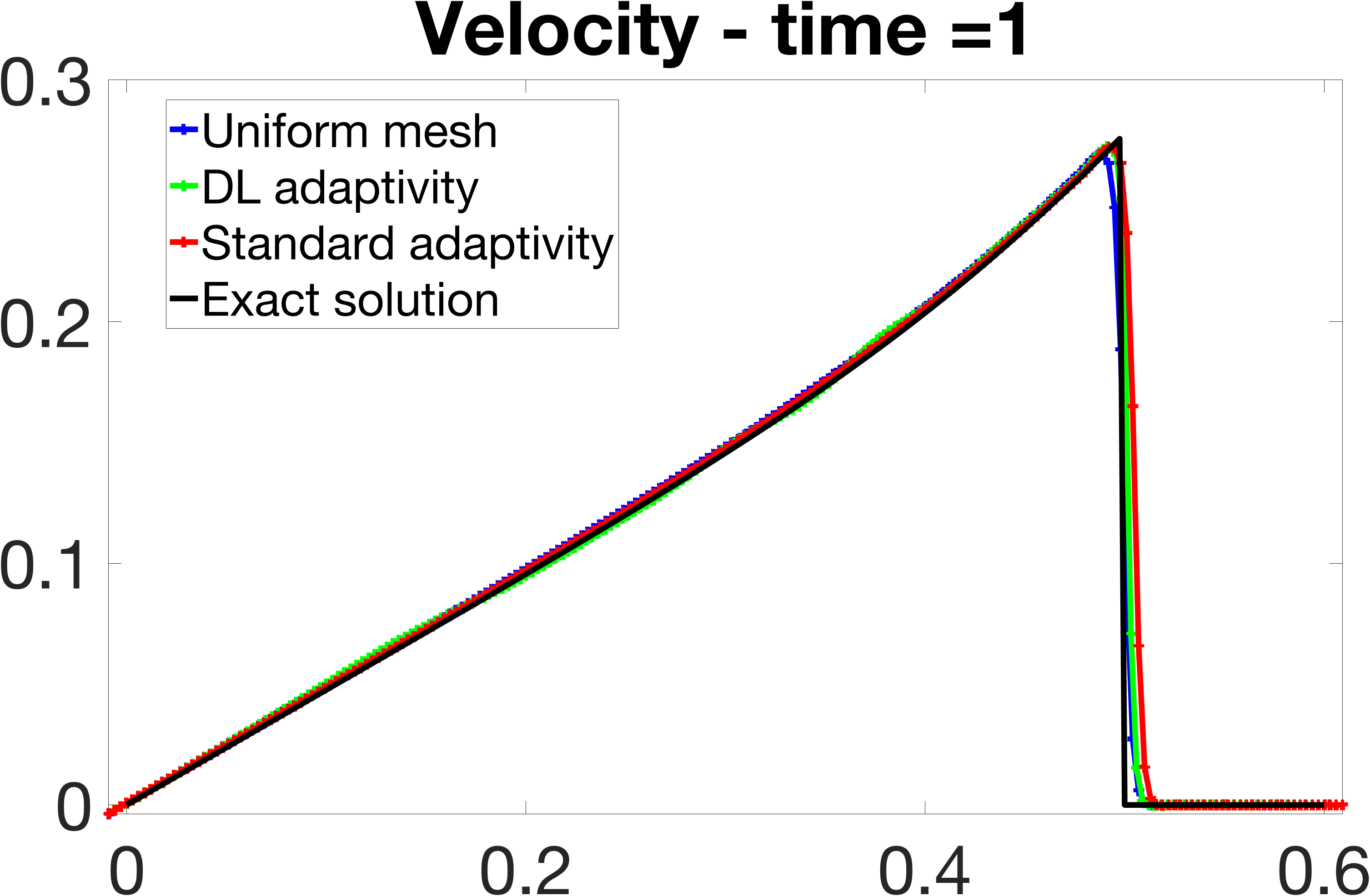}   
\hspace{1cm}
\includegraphics[width=0.45\textwidth]{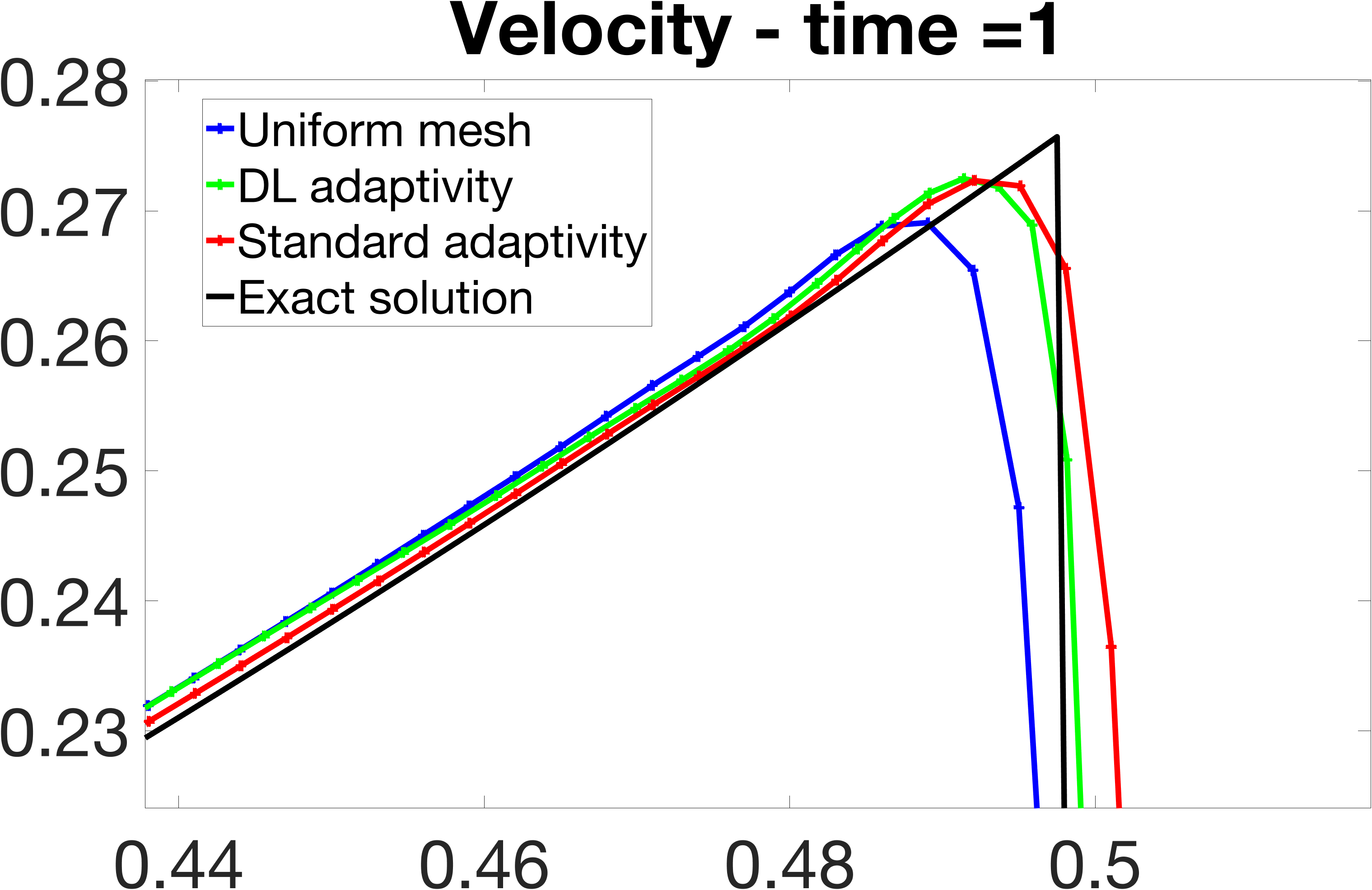} \\
\vspace{1cm}
\includegraphics[width=0.45\textwidth]{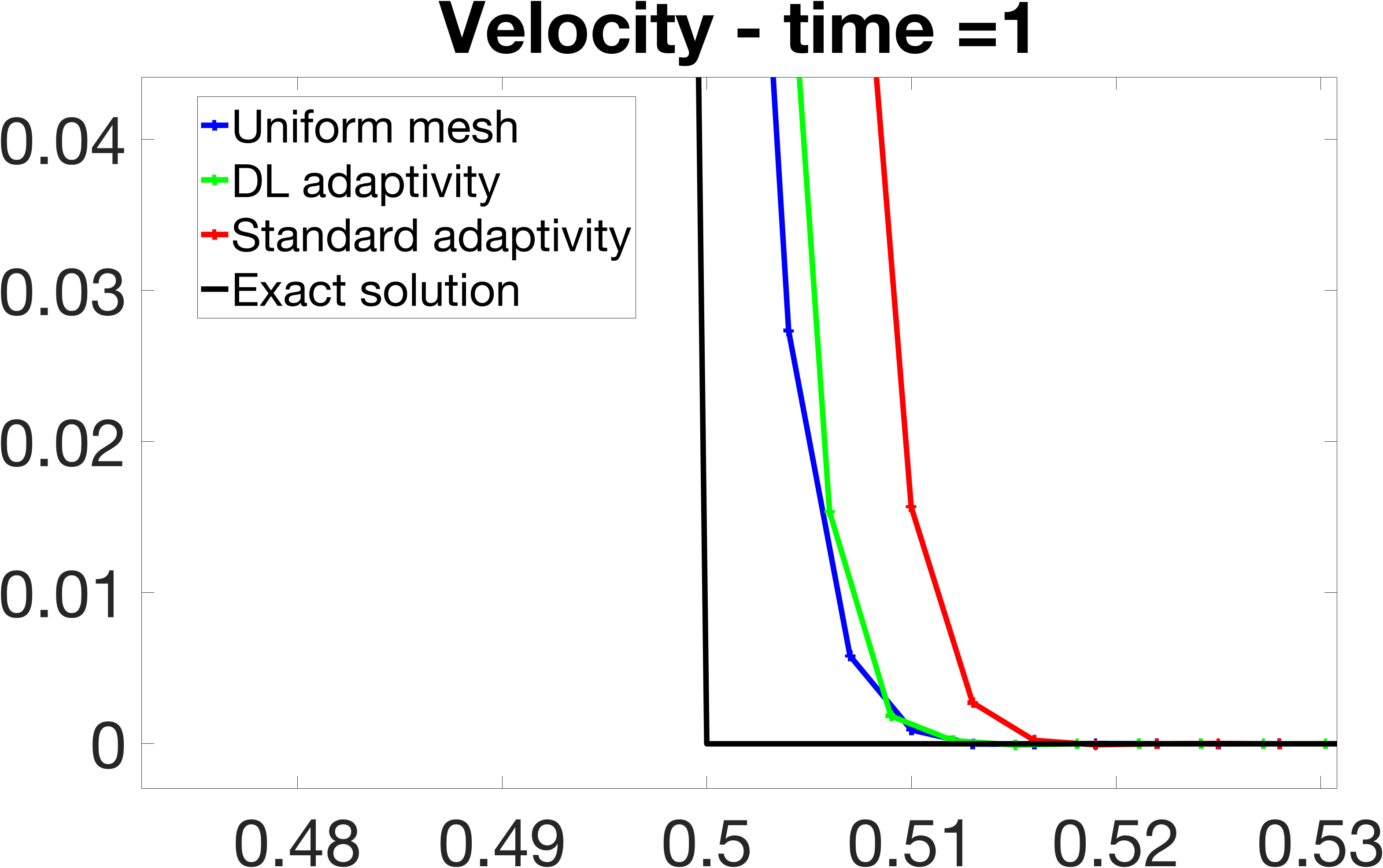} 
\caption{Velocity profile of 1D Taylor–von Neumann–Sedov test case at $t=1.0$. WENO5 scheme used for space discretization and RK-3 scheme used for time stepping. Entire profile shown at the top left, zoom-in in the $x$ range [0.4, 0.55] shown at the top right, zoom-in in the $x$ range [0.48, 0.53] shown at the bottom.}
\label{sedov_weno5_velocity}
\end{figure}

\begin{figure}[H]
\centering
\includegraphics[width=0.45\textwidth]{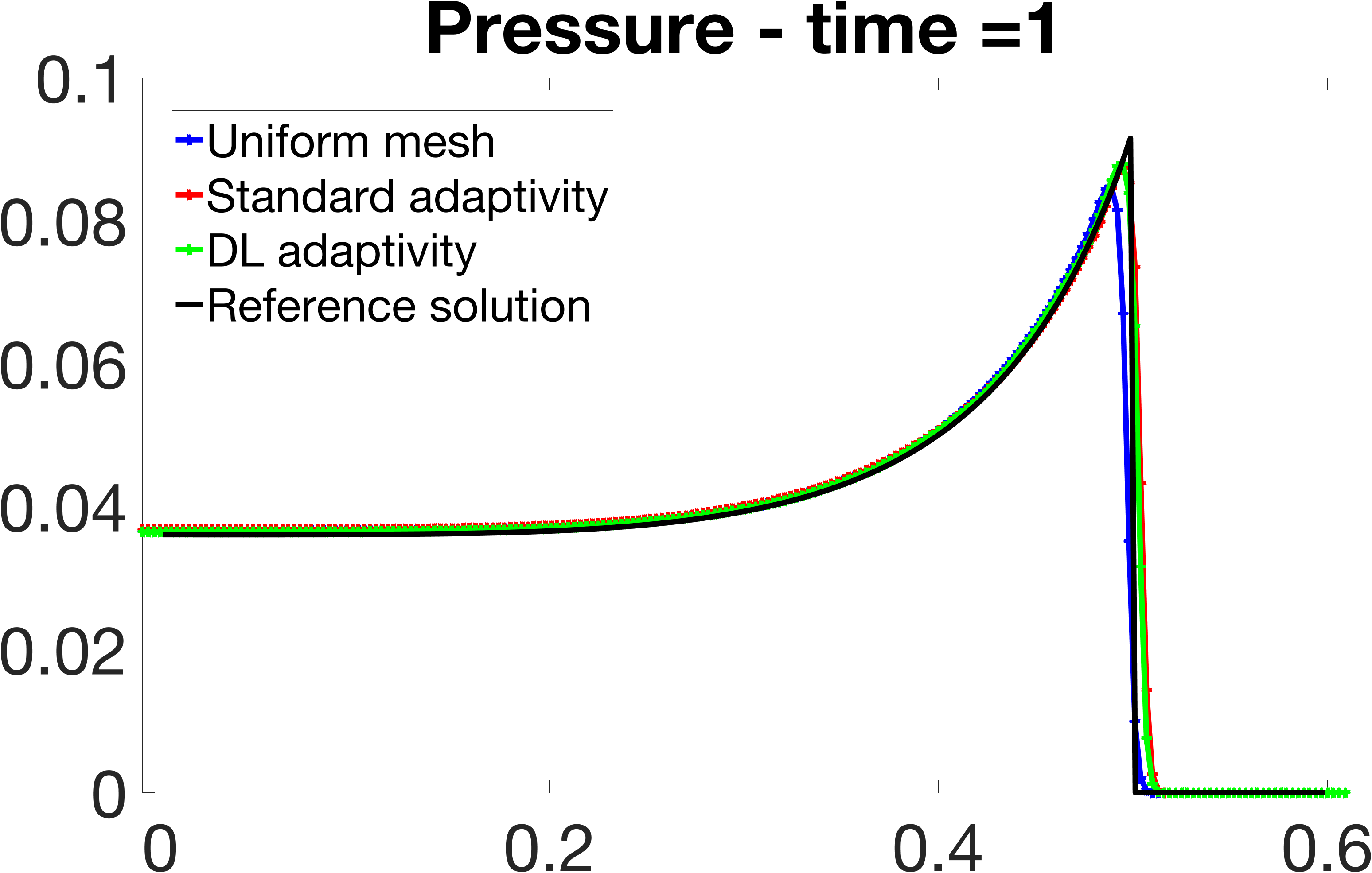}  
\hspace{1cm}
\includegraphics[width=0.45\textwidth]{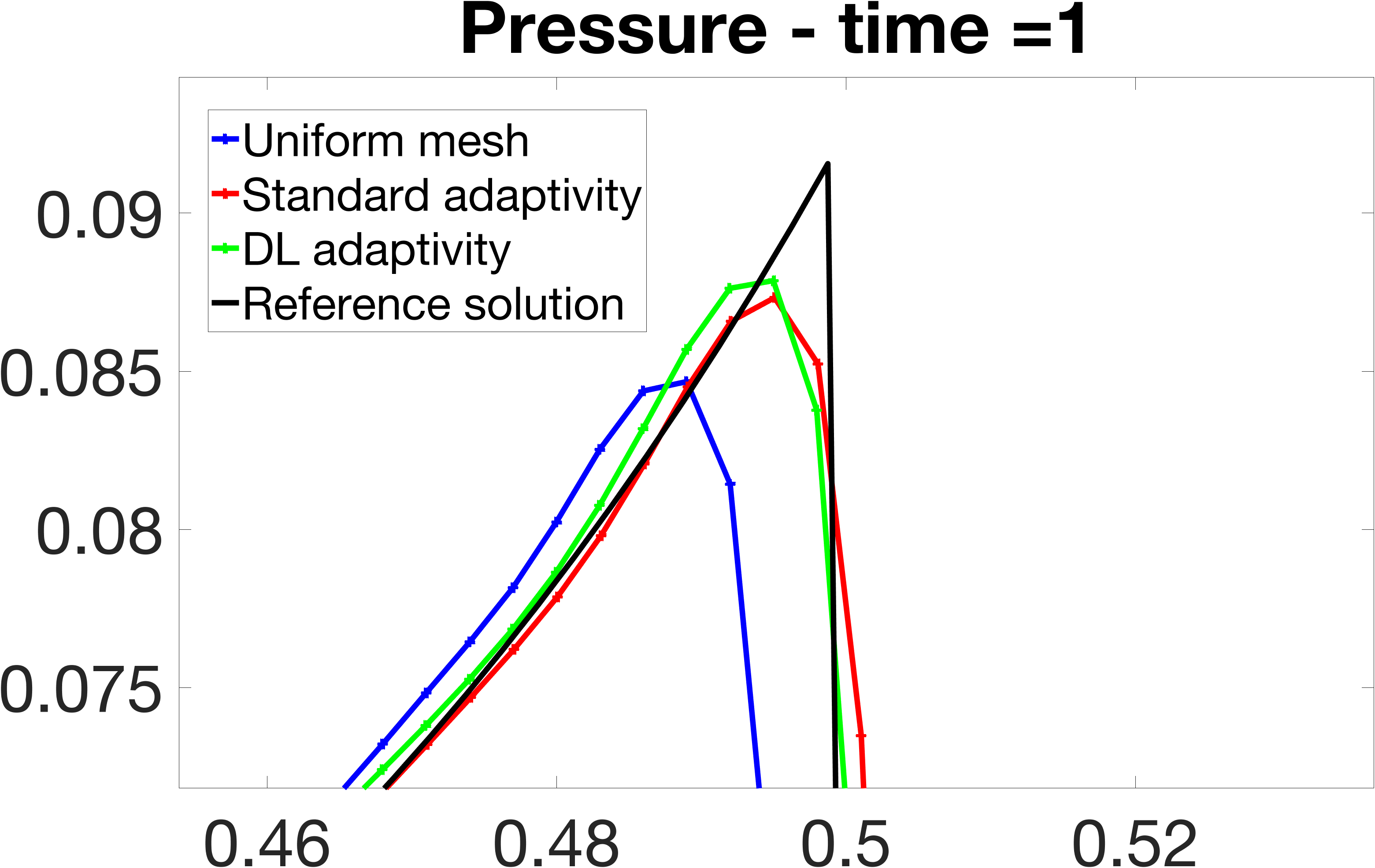}  
\caption{Pressure profile of 1D Taylor–von Neumann–Sedov test case at $t=1.0$. WENO5 scheme used for space discretization and RK-3 scheme used for time stepping. Entire profile shown at left and zoom-in in the $x$ range [0.4, 0.5] shown at right.}
\label{sedov_weno5_pressure}
\end{figure}

\begin{figure}[H]
\centering
\includegraphics[width=0.45\textwidth]{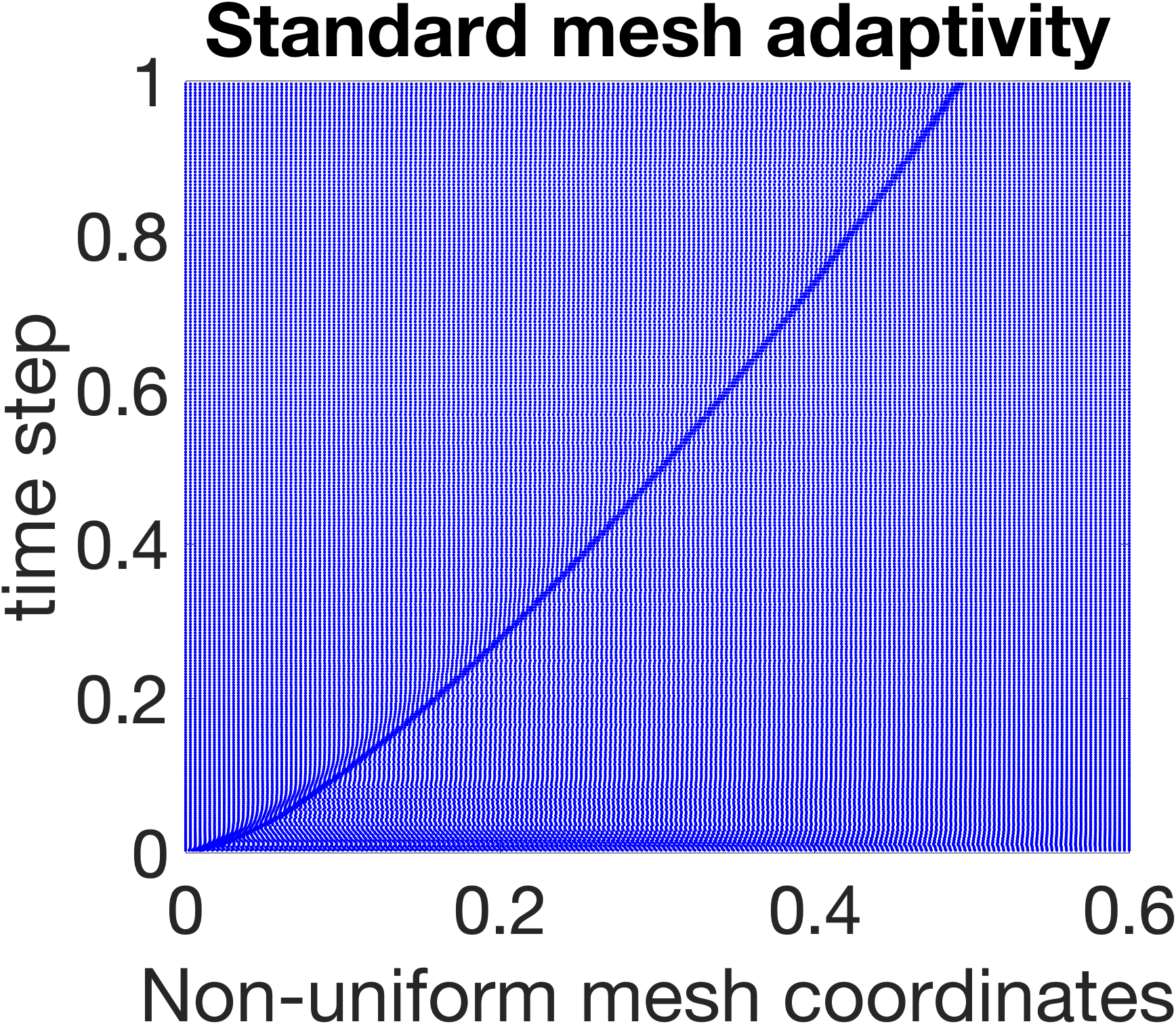}  
\hspace{1cm} 
\includegraphics[width=0.44\textwidth]{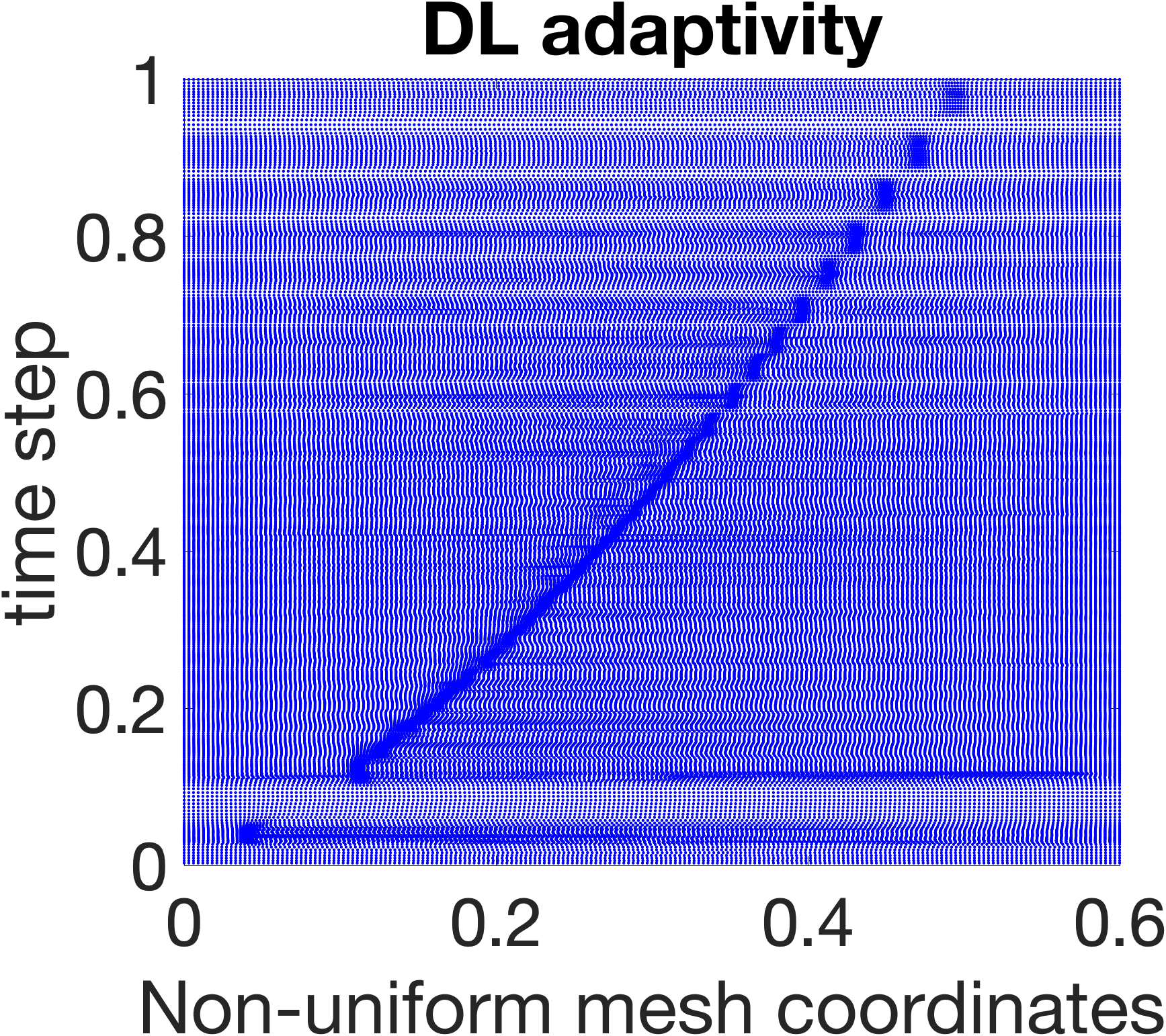}   
\caption{History of non-uniform adapted meshes for 1D Taylor–von Neumann–Sedov test case using standard adaptive zoning (left) and deep learning adaptive zoning (right). }
\label{mesh_history_sedov}
\end{figure}

The $L^1$ and $L^2$ norms of exact, numerical approximations, and numerical error for density and velocity using WENO5 for space discretization and RK-3 for time stepping are provided in Table \ref{tab:sedov_norms}. The DL adaptivity allows the physics PDE solver to attain an accuracy comparable to the one obtained with the standard adaptivity and better than the uniform mesh. 
The numerical reconstruction of density, velocity, and pressure profiles using the WENO5 space discretization and RK-3 for time stepping  are shown in Figures \ref{sedov_weno5_density}, \ref{sedov_weno5_velocity}, and \ref{sedov_weno5_pressure} respectively. Both standard adaptive zoning and DL adaptive zoning still dampen the numerical oscillations in the numerical reconstruction of each physical quantity with respect to the calculations performed on a uniform mesh. 
DL adaptive zoning still runs faster than standard adaptive zoning for each time step, with standard adaptive zoning taking one wall-clock second and DL adaptive zoning taking half a second per time step. 
The evolution of the meshes across consecutive time steps for standard adaptive zoning and DL adaptive zoning are shown in Figure \ref{mesh_history_sedov}. The change of the mesh across consecutive time steps is smoother for the standard adaptive zoning, whereas the mesh updates produced by the DL adaptivity are more irregular. 

\begin{table}
    \centering
    \begin{tabular}{|c|c|c|c|c|}
    \hline
    \multicolumn{5}{|c|}{\textbf{Density}}\\
    \hline
    & \textbf{Exact solution}& \textbf{Uniform mesh} & \textbf{Standard adaptivity} & \textbf{DL adaptivity}\\
    \hline
    $L^2$-norm of the solution & 1.226 & 1.2154 & 1.2134 & 1.2140\\  
    \hline
    $L^2$-norm of the relative error & - & 0.1543 & 0.1867 & 0.1549\\  
    \hline
    $L^1$-norm of the solution & 0.5988 & 0.6015 & 0.6007 & 0.6012\\ 
    \hline
    $L^1$-norm of the relative error & - & 0.0584 & 0.0661 & 0.0557 \\      
    \hline
    \multicolumn{5}{|c|}{\textbf{Velocity}}\\
    \hline
    & \textbf{Exact solution}& \textbf{Uniform mesh} & \textbf{Standard adaptivity} & \textbf{DL adaptivity}\\
    \hline
    $L^2$-norm of the solution & 0.1039 & 0.1048 & 0.1043 & 0.1046\\  
    \hline
    $L^2$-norm of the relative error & - & 0.0706 & 0.0680 & 0.0697\\  
    \hline
    $L^1$-norm of the solution & 0.0625 & 0.0637 & 0.0633 & 0.0636\\ 
    \hline
    $L^1$-norm of the relative error & - & 0.0271 & 0.0271 & 0.0260 \\ 
    \hline
    \multicolumn{5}{|c|}{\textbf{Pressure}}\\
    \hline
    & \textbf{Exact solution}& \textbf{Uniform mesh} & \textbf{Standard adaptivity} & \textbf{DL adaptivity}\\
    \hline
    $L^2$-norm of the solution & 0.0324 & 0.1048 & 0.1043 & 0.1046\\  
    \hline
    $L^2$-norm of the relative error & - & 0.0706 & 0.0680 & 0.0697\\  
    \hline
    $L^1$-norm of the solution & 0.0219 & 0.0220 & 0.0216 & 0.0220\\ 
    \hline
    $L^1$-norm of the relative error & - & 0.0263 & 0.0249 & 0.0250 \\ 
    \hline    
    \end{tabular}
    \caption{1D Taylor–von Neumann–Sedov test case. Relative error for density and velocity in the $L^2$-norm and $L^1$-norm using the WENO5 scheme used for space discretization and the RK-3 scheme used for time stepping.}
    \label{tab:sedov_norms}
\end{table}

\begin{table}
    \centering
    \begin{tabular}{|c|c|}
    \hline
    \textbf{Mesh Type} & \textbf{Wall-clock Time (s)}\\
    \hline
    Uniform mesh & 1.92\\
    \hline
    Standard adaptivity & 32.84\\
    \hline
    DL adaptivity & $1127.75 - (0.1476 \times 7,612) = 4.22$\\
    \hline
    \end{tabular}
    \caption{1D Taylor–von Neumann–Sedov test case. Wall-clock computational time in seconds using the WENO5 scheme for space discretization and the RK-3 scheme for time stepping on a uniform mesh, non-uniform mesh computed with standard adaptive zoning, and non-uniform mesh computed with DL adaptive zoning. The time spent to perform \texttt{Python} calls is subtracted from total computational time for the DL approach.}
    \label{tab:sedov_weno5_time}
\end{table}

\subsection{Woodward-Colella blast waves test case}
We now present the numerical results the performance of DL adaptive zoning compared against standard adaptive zoning with WENO5 space discretization scheme for the 1D Woodward-Colella blast waves test case \cite{woodward} described by the 1D Euler equations in the domain [0,1] with the following initial conditions 
\begin{equation*}
    \rho(x, t=0) = 
    \begin{cases}
     1.0\quad \text{if} \quad x<=0.5\\
     \frac{1}{8}\quad \text{if} \quad x>0.5\\
    \end{cases}, \quad 
    u(x, t=0) = 0.0, \quad 
   p(x, t=0) = 
    \begin{cases}
     1.0\quad \text{if} \quad x<=0.5\\
     \frac{1}{10}\quad \text{if} \quad x>0.5\\
    \end{cases}\\    
\end{equation*}
and the initial value of the energy is obtained from the state equations described in Section \ref{background}. 
The parameters used for the standard adaptivity and generate the training dataset for DL are described in Table \ref{parameters_mmpde2}.
The $L^1$ and $L^2$ norms of exact, numerical approximations, and numerical error for density, velocity, pressure, and energy using WENO5 for space discretization and RK-3 for time stepping are provided in Table \ref{tab:woodward_norms}. The numerical solution computed by the physics PDE solver with DL adaptivity is comparable to the one obtained with the standard adaptivity, and both adaptive methods attain better accuracy than the uniform mesh. 

\begin{figure}[H]
\centering
\includegraphics[width=0.45\textwidth]{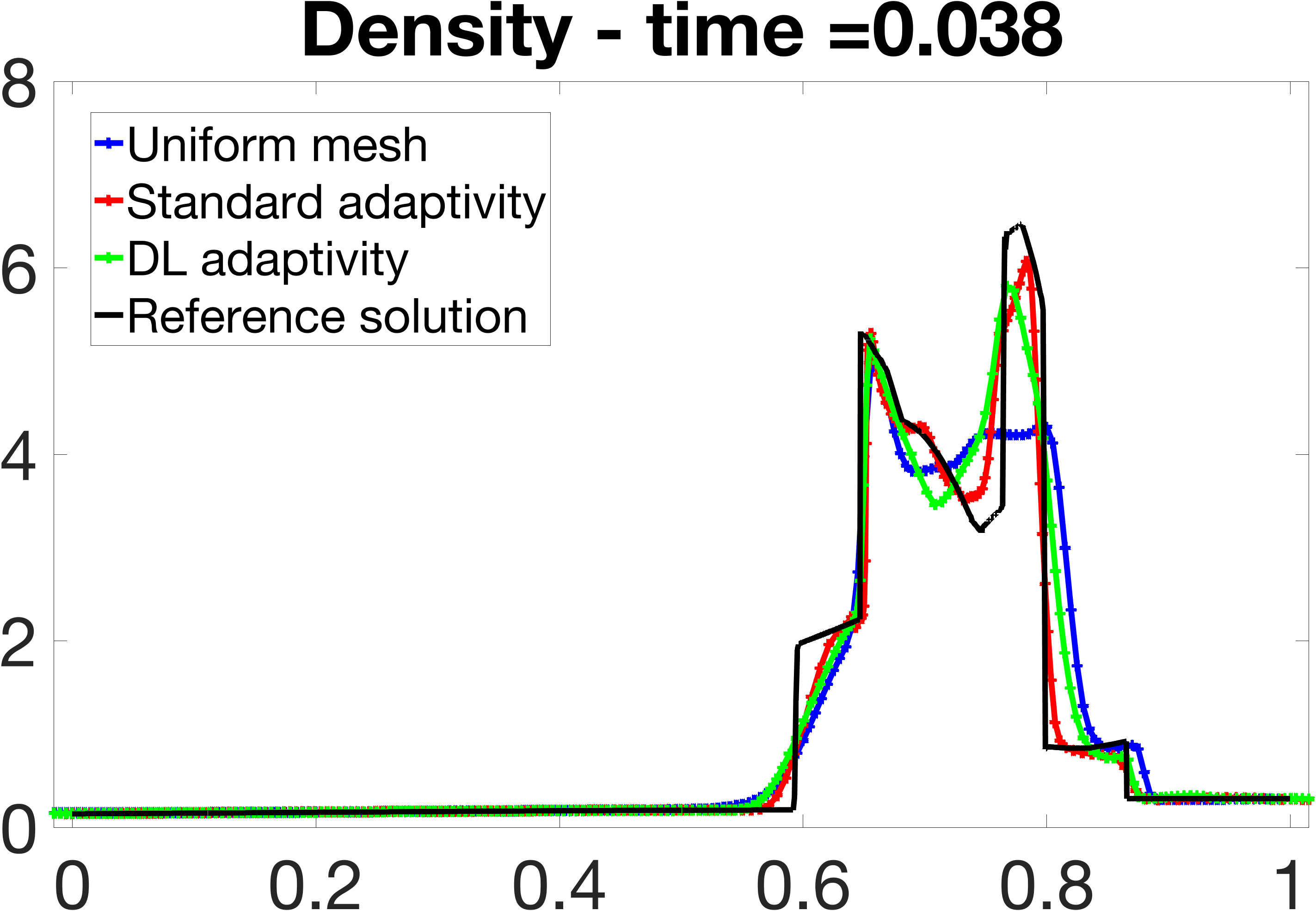}  \\
\vspace{1cm}
\includegraphics[width=0.45\textwidth]{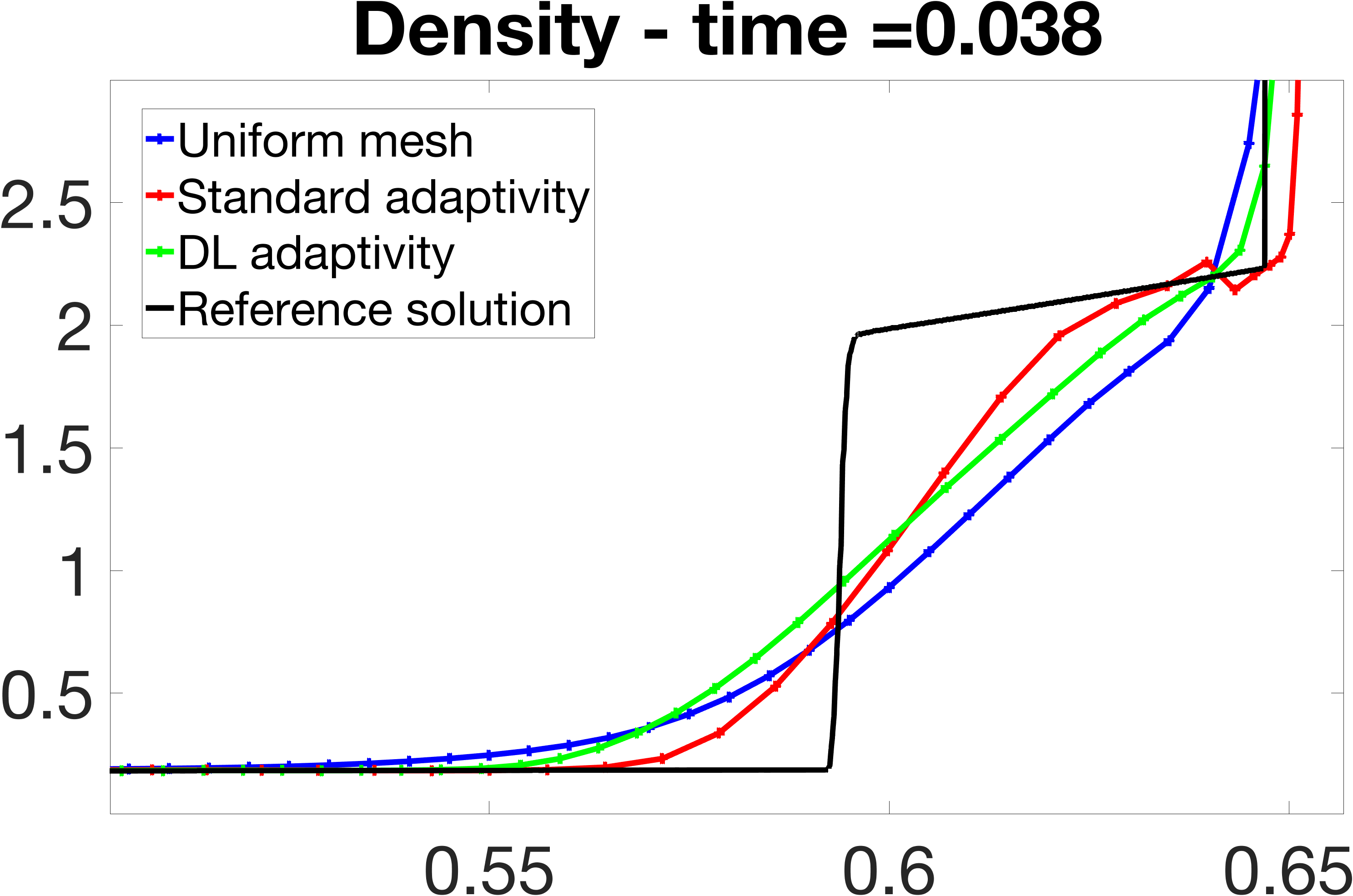} \hspace{1cm}
\includegraphics[width=0.45\textwidth]{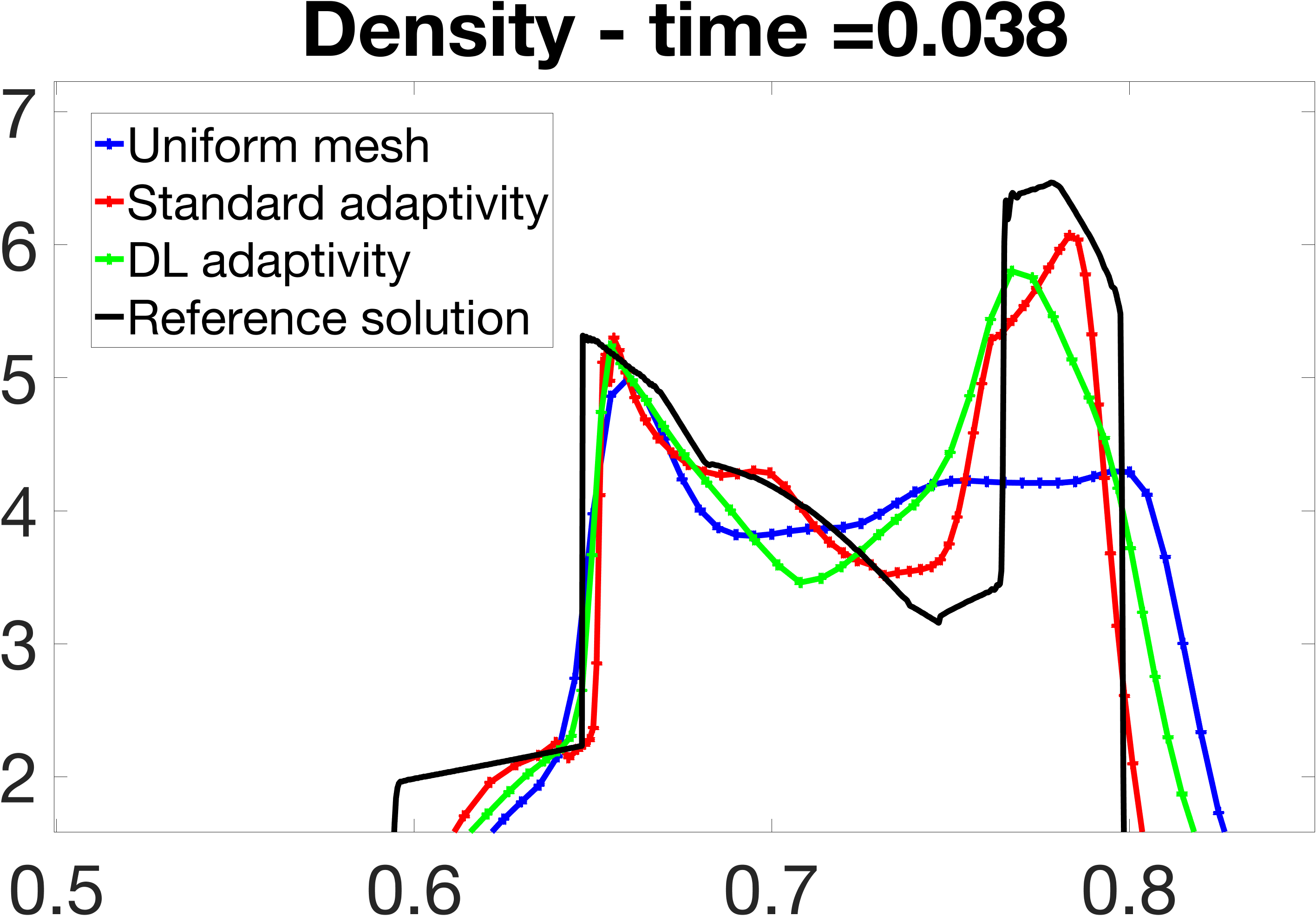}
\caption{Density profile of 1D Woodward-Colella test case at $t=0.038$ using the WENO5 scheme for space discretization and the RK-3 scheme used for time stepping. Entire profile shown at the top, zoom-in in the $x$ range [0.5, 0.65] shown at the bottom left, and zoom-in in the $x$ range [0.6, 0.8].}
\label{woodward_weno5_density}
\end{figure}

\begin{figure}[H]
\centering
\includegraphics[width=0.45\textwidth]{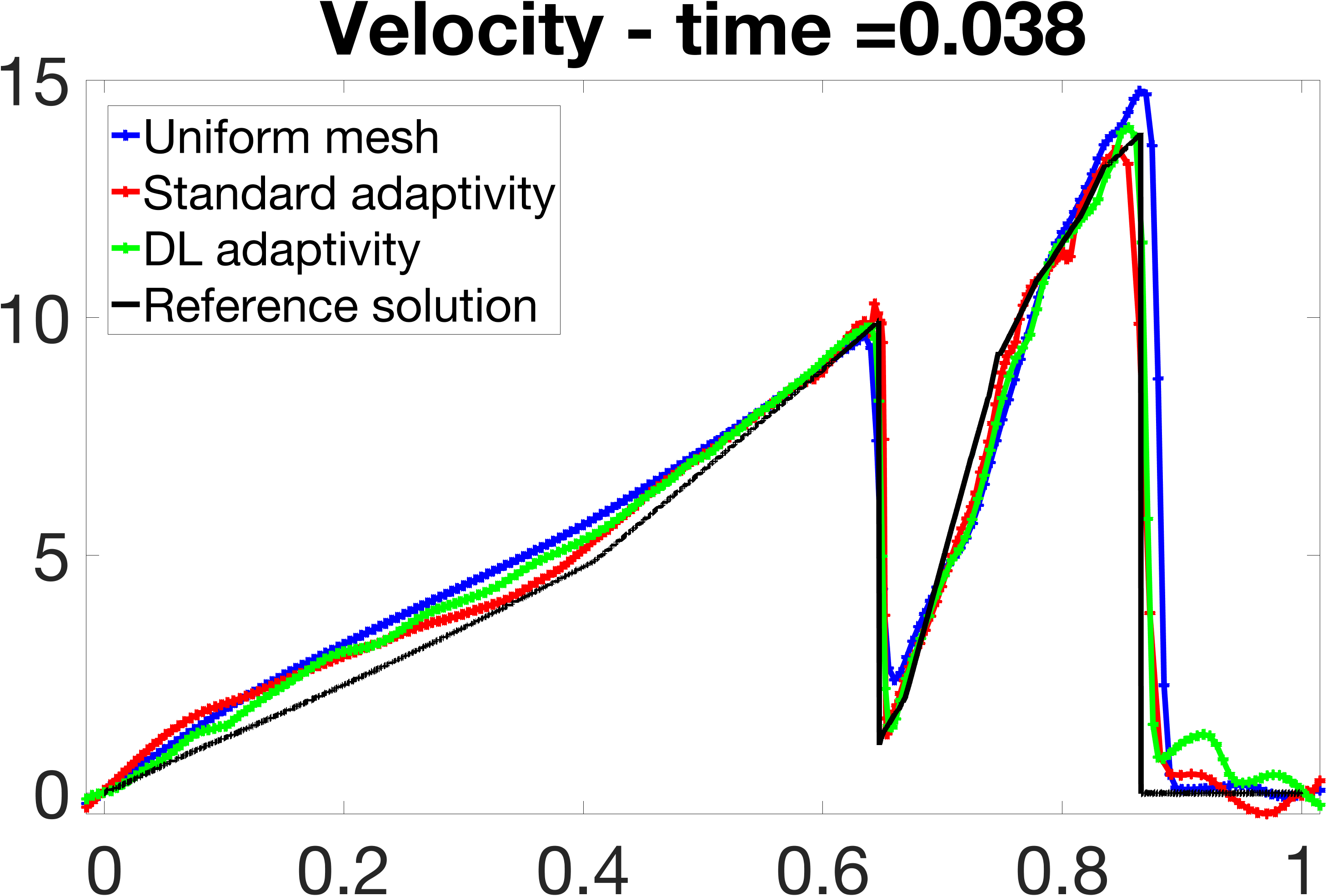}  \\
\caption{Velocity profile of 1D Woodward-Colella test case at $t=0.038$ using the WENO5 scheme for space discretization and the RK-3 scheme for time stepping.}
\label{woodward_weno5_velocity}
\end{figure}

\begin{figure}[H]
\centering
\includegraphics[width=0.45\textwidth]{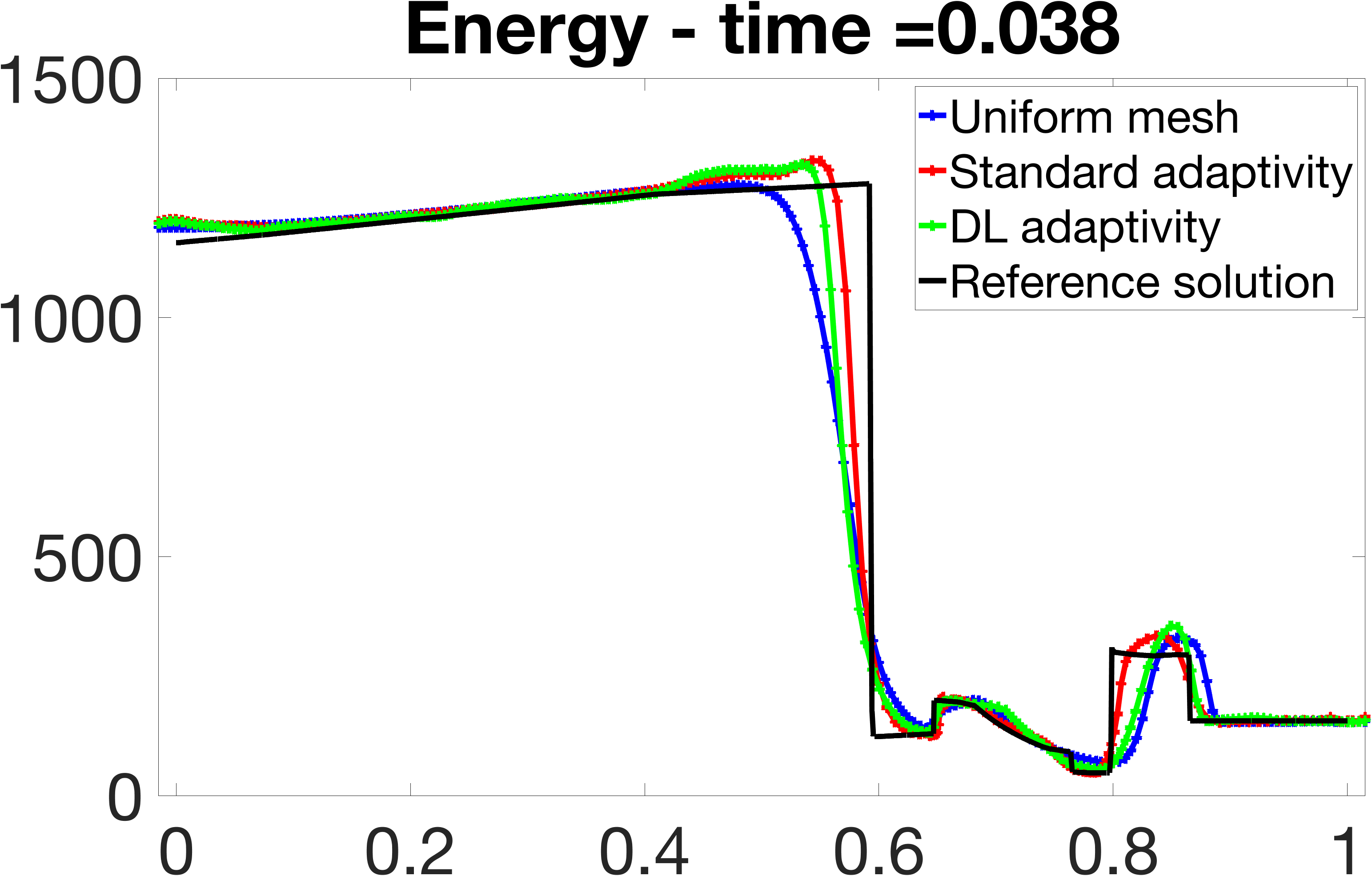}  \\
\caption{Energy profile of 1D Woodward-Colella test case at $t=0.038$ using the WENO5 scheme for space discretization and the RK-3 scheme for time stepping.}
\label{woodward_weno5_energy}
\end{figure}

\begin{figure}[H]
\centering
\includegraphics[width=0.45\textwidth]{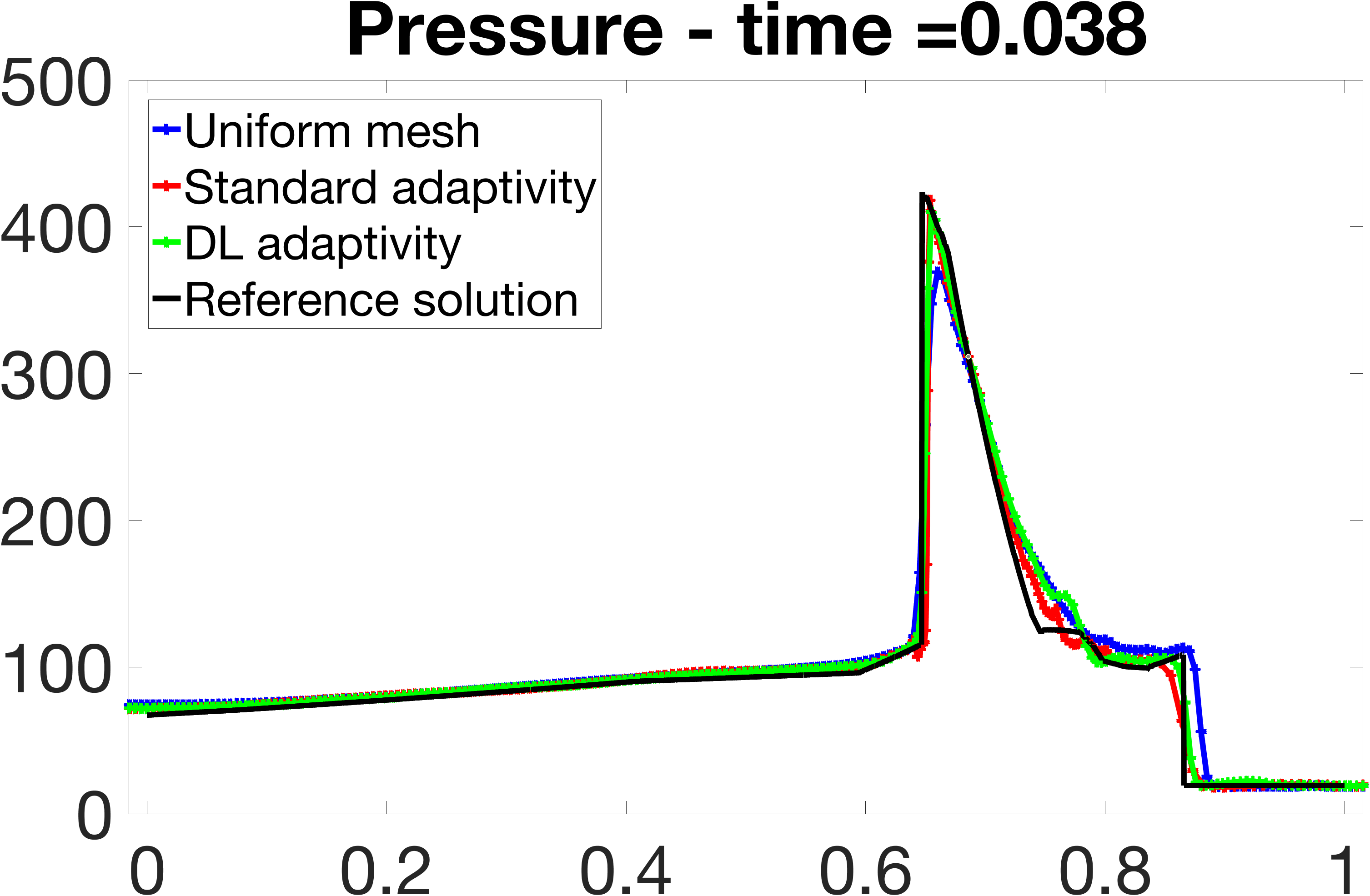}  \hspace{1cm}
\includegraphics[width=0.45\textwidth]{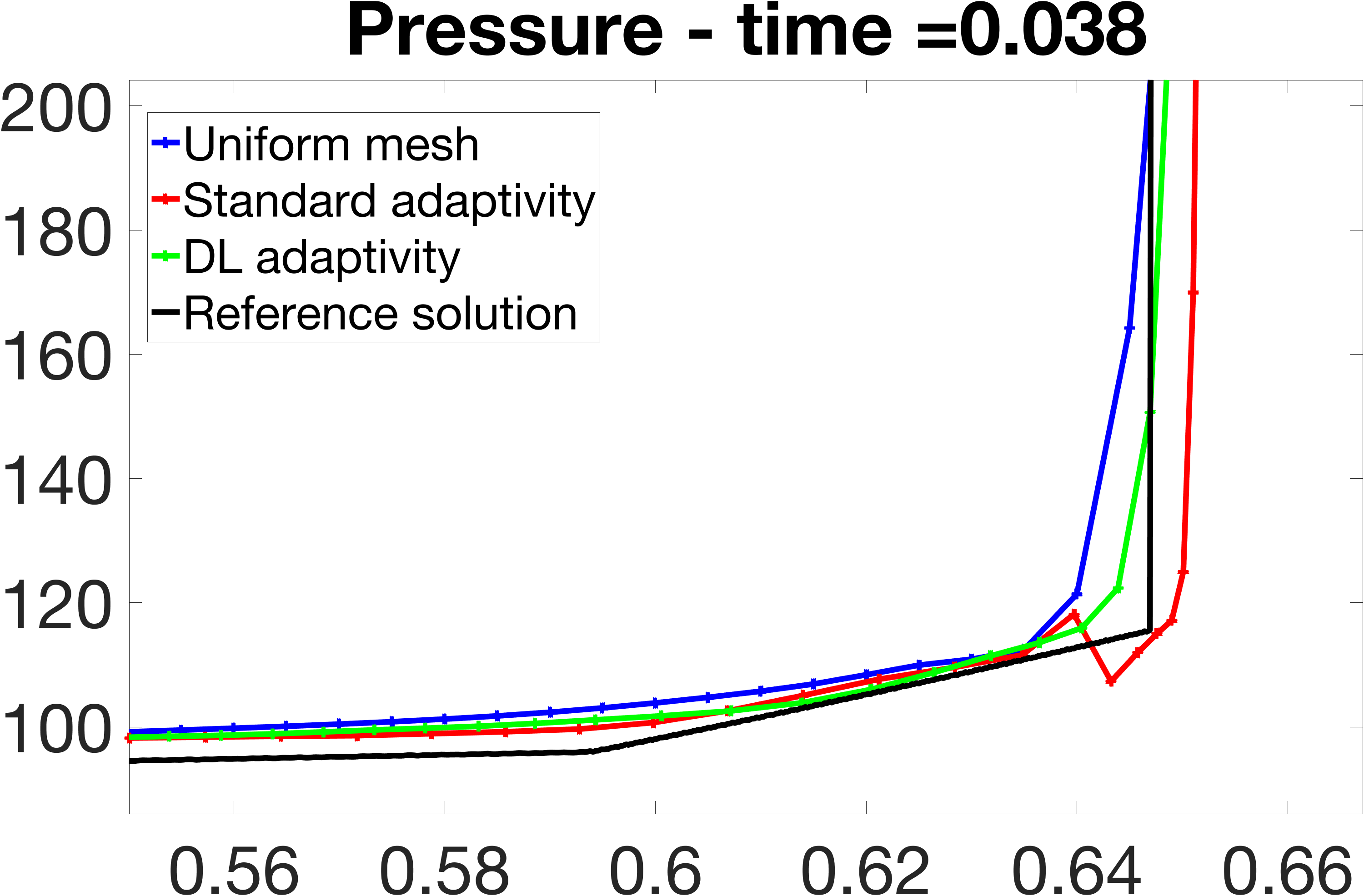} \\
\vspace{1cm}
\includegraphics[width=0.45\textwidth]{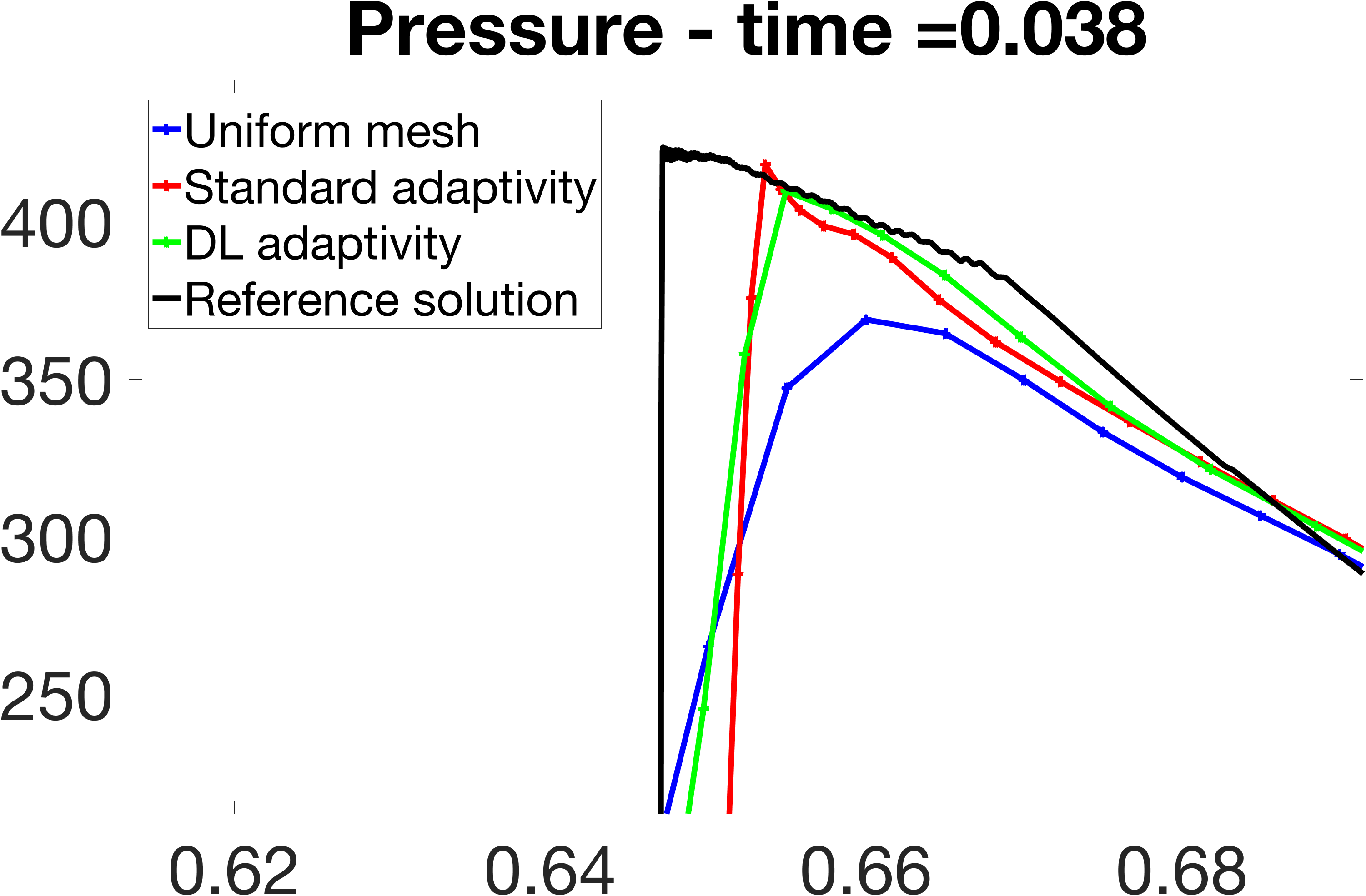} \hspace{1cm}
\includegraphics[width=0.45\textwidth]{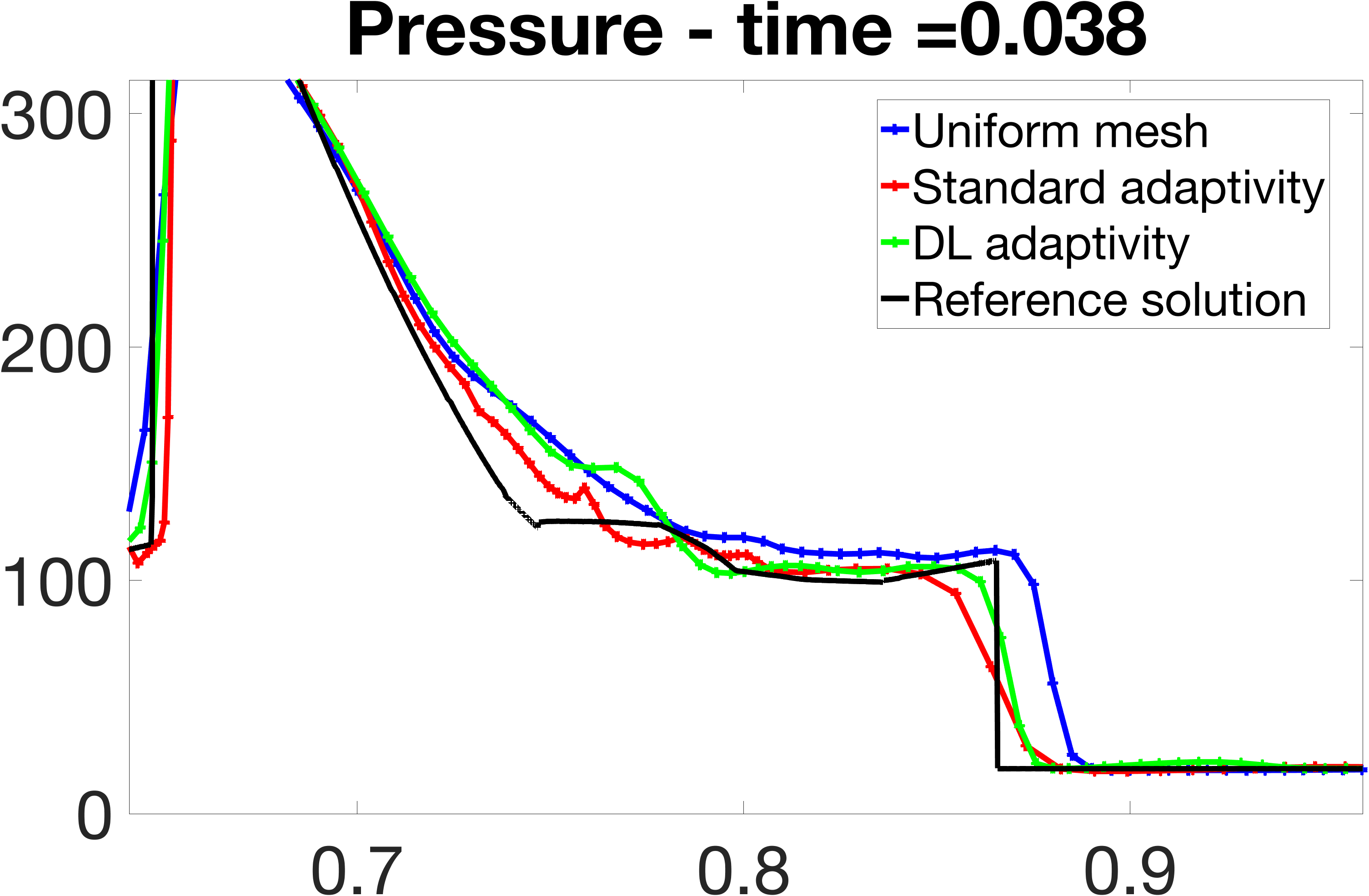} 
\caption{Pressure profile of 1D Woodward-Colella test case at $t=0.038$ using the WENO5 scheme used for space discretization and the RK-3 scheme for time stepping. Entire profile shown at the top left, zoom-in in the $x$ range [0.55, 0.65] shown at the top right, zoom-in in the $x$ range [0.62, 0.68] shown at the bottom left, and zoom-in in the $x$ range [0.7, 0.9] shown at the bottom right.}
\label{woodward_weno5_pressure}
\end{figure}

The profiles of density, velocity, pressure, and energy profiles using WENO5 for space discretization and RK-3 for time steppin g are shown in Figures \ref{woodward_weno5_density}, \ref{woodward_weno5_velocity}, \ref{woodward_weno5_energy}, and \ref{woodward_weno5_pressure}, respectively. The density profile reconstructed with DL adaptivity does not capture the peaks as accurately as the standard adaptivity, whereas the reconstruction of velocity, energy, and pressure obtained with DL adaptivity and standard adaptivity are comparable and are overall close to the reference solution than the reconstruction obtained with the uniform mesh. 

\begin{table}
\begin{center}
\begin{tabular}{ |c| c| }
\hline
\multicolumn{2}{|c|}{\textbf{Parabolic MMPDE }}\\
\hline
Definition of monitor function & Equation \eqref{monitor_function}\\
\hline
 $\alpha_1$ & 0.0 \\ 
 \hline
 $\alpha_2$ &600.0 \\ 
 \hline
 $\alpha_3$ & 0.0 \\  
 \hline
 $\alpha_4$ & 0.0 \\  
 \hline
 type of smoothing for monitor function & Gaussian\\
 \hline
 size of sliding window for smoothing & $10\%$ of domain\\
 \hline
 number of time steps & 1,000 \\
 \hline
\end{tabular}
\caption{Parameters of the monitor function and fixed point solver to perform standard adaptivity and generate the training dataset for DL using the parabolic MMPDE in Woodward-Colella blast waves test case.}
\label{parameters_mmpde2}
\end{center}
\end{table}

\begin{table}
    \centering
    \begin{tabular}{|c|c|c|c|c|}
    \hline
    \multicolumn{5}{|c|}{\textbf{Density}}\\
    \hline
    & \textbf{Reference solution}& \textbf{Uniform mesh} & \textbf{Standard adaptivity} & \textbf{DL adaptivity}\\
    \hline
    $L^2$-norm of the solution & 1.9018 & 1.7815 & 1.8150 & 1.8522\\  
    \hline
    $L^2$-norm of the relative error & - & 0.3311 & 0.1956 & 0.2365\\  
    \hline
    $L^1$-norm of the solution & 1.0000 & 1.0039 & 0.9711 & 1.0124\\ 
    \hline
    $L^1$-norm of the relative error & - & 0.2522 & 0.1202 & 0.1820 \\      
    \hline
    \multicolumn{5}{|c|}{\textbf{Velocity}}\\
    \hline
    & \textbf{Reference solution}& \textbf{Uniform mesh} & \textbf{Standard adaptivity} & \textbf{DL adaptivity}\\
    \hline
    $L^2$-norm of the solution & 6.0089 & 6.4313 & 6.1148 & 6.1330\\  
    \hline
    $L^2$-norm of the relative error & - & 0.2844 & 0.1483 & 0.1802\\  
    \hline
    $L^1$-norm of the solution & 4.5984 & 5.1780 & 4.8878 & 4.9377\\ 
    \hline
    $L^1$-norm of the relative error & - & 0.1670 & 0.1047 &  0.1196 \\      
    \hline
    \multicolumn{5}{|c|}{\textbf{Internal energy}}\\
    \hline
    & \textbf{Reference solution}& \textbf{Uniform mesh} & \textbf{Standard adaptivity} & \textbf{DL adaptivity}\\
    \hline
    $L^2$-norm of the solution & 950.3885 & 927.4720 & 949.0341 & 939.2098\\  
    \hline
    relative $L^2$-norm of the error & - & 0.1487 & 0.1030 & 0.1430\\  
    \hline
    $L^1$-norm of the solution & 793.0439 & 770.7242 & 790.0884 & 778.4181\\ 
    \hline
    relative $L^1$-norm of the error & - & 0.0655 & 0.0395 & 0.0555 \\      
    \hline
    \multicolumn{5}{|c|}{\textbf{Pressure}}\\
    \hline
    & \textbf{Reference solution}& \textbf{Uniform mesh} & \textbf{Standard adaptivity} & \textbf{DL adaptivity}\\
    \hline
    $L^2$-norm of the solution & 120.7597 & 121.0085 & 118.9435 & 121.8832\\  
    \hline
    relative $L^2$-norm of the error & - & 0.1475 & 0.1808 & 0.1115\\  
    \hline
    $L^1$-norm of the solution & 97.6510 & 101.9923 & 98.5290 & 100.6408\\ 
    \hline
    relative $L^1$-norm of the error & - & 0.0812 & 0.0584 & 0.0569 \\      
    \hline
    \end{tabular}
    \caption{1D Woodward-Colella test case. Relative error for density, velocity, internal energy, and pressure in the $L^2$-norm and $L^1$-norm using the WENO5 scheme for space discretization and the RK-3 scheme for time stepping.}
    \label{tab:woodward_norms}
\end{table}

\begin{table}
    \centering
    \begin{tabular}{|c|c|}
    \hline
    \textbf{Mesh Type} & \textbf{Wall-clock Time (s)}\\
    \hline
    Uniform mesh & 1.21 \\
    \hline
    Standard adaptivity & 350.48\\
    \hline
    DL adaptivity & $2,056.26 - (0.1476\times 13,528) = 59.53$\\
    \hline
    \end{tabular}
    \caption{1D Woodward-Colella test case. Wall-clock computational time in seconds using the WENO5 scheme for space discretization and the RK-3 scheme for time stepping on a uniform mesh, non-uniform mesh computed with standard adaptive zoning, and non-uniform mesh computed with DL adaptive zoning.}
    \label{tab:woodward_weno5_time}
\end{table}

\begin{figure}[H]
\centering
\includegraphics[width=0.5\textwidth]{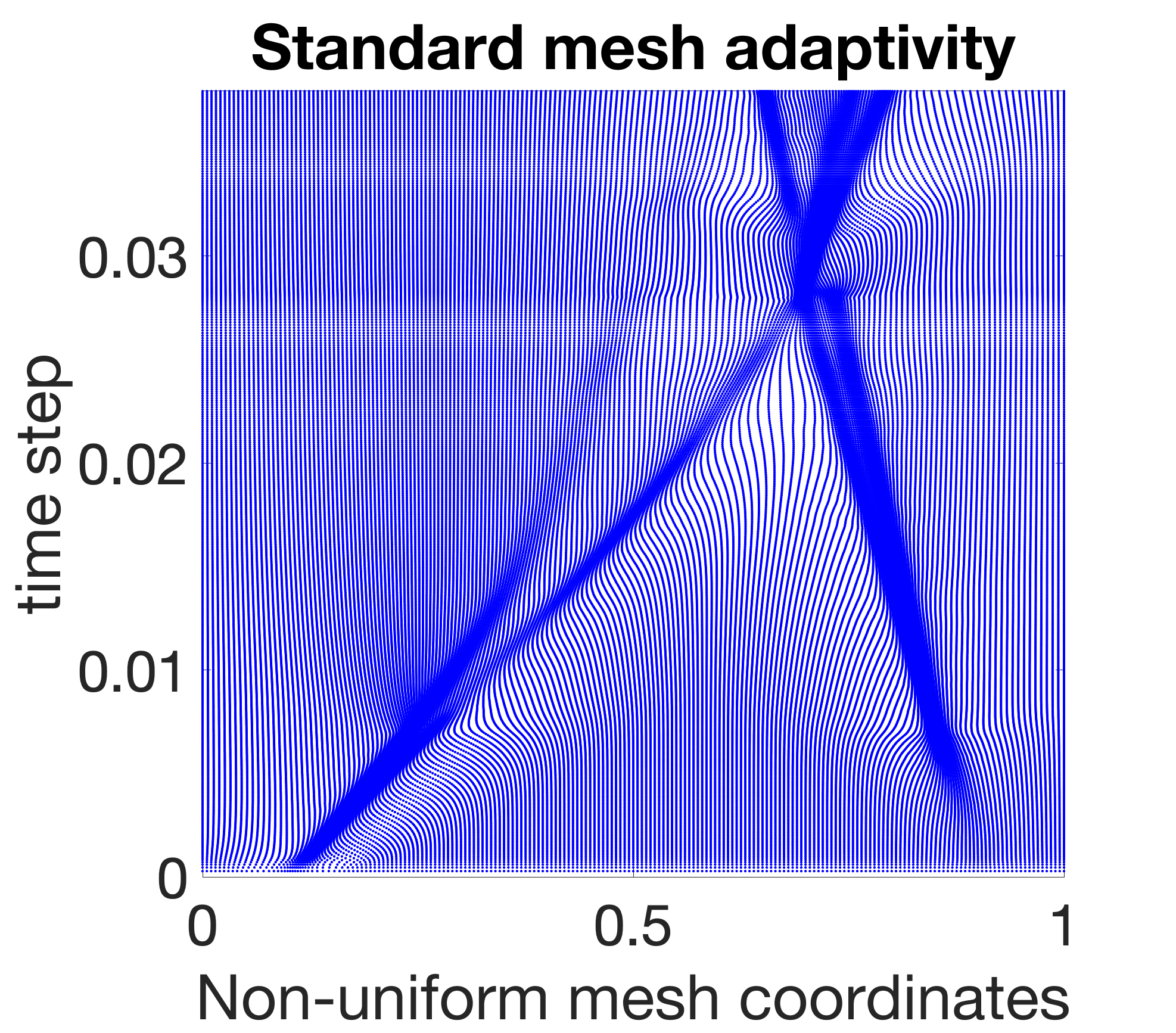}  
\includegraphics[width=0.49\textwidth]{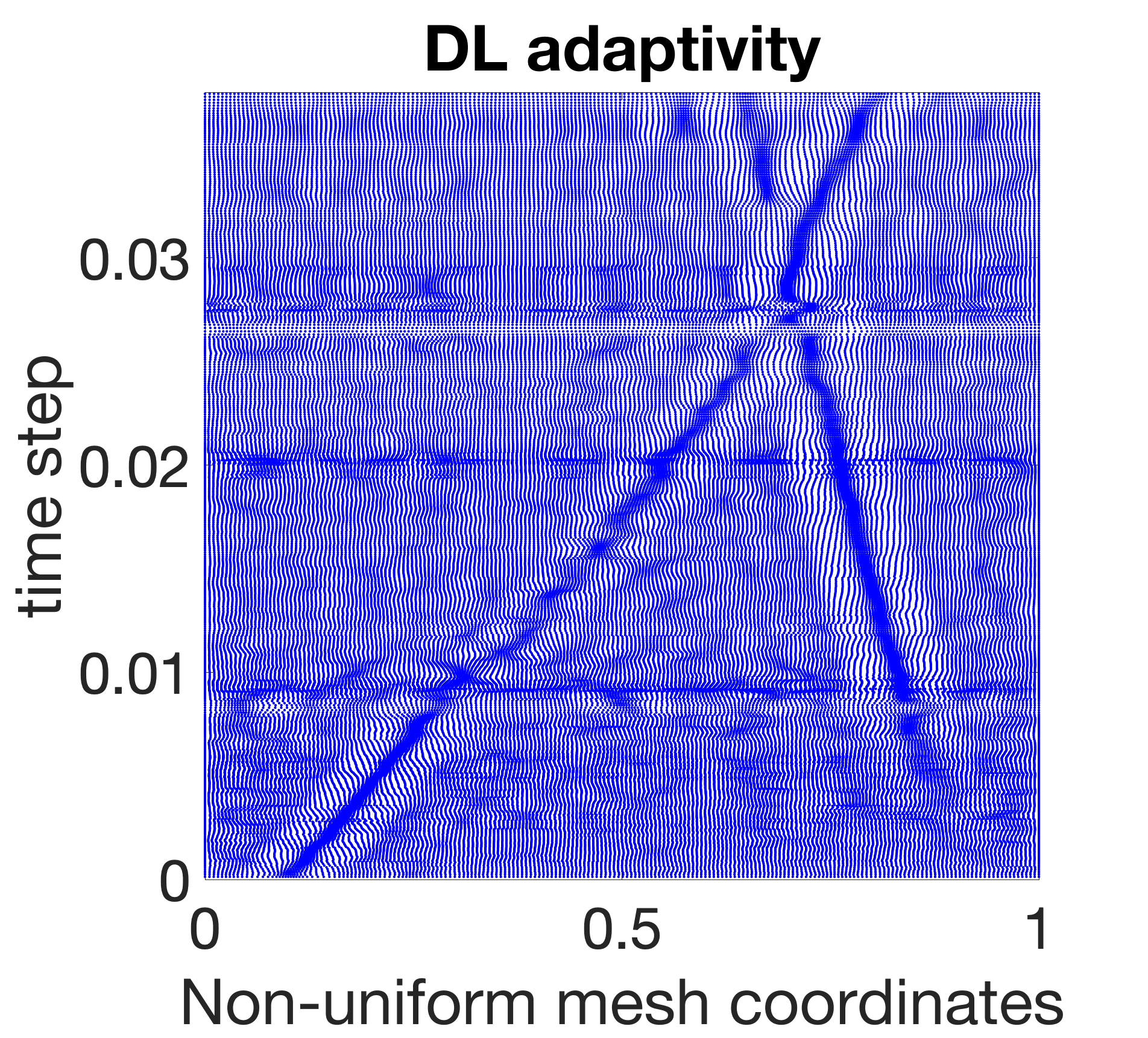}   
\caption{History of non-uniform adapted meshes for 1D Woodward-Colella test case using standard adaptive zoning (left) and deep learning adaptive zoning (right). }
\label{mesh_history_woodward}
\end{figure}

\section{Discussion and future developments}
\label{conclusion_section}
We presented a DL approach to perform adaptive zoning that changes the location of grid points in a numerical discretization scheme to increase the concentration of nodes close to discontinuities in the quantity of interest. The training dataset consists of artificial staircase shock profiles where the number and the location of the jumps changes arbitrarily across the data samples.  
This particular choice of the training dataset forces the DL model to focus only on identifying regions with discontinuities, whereas all regions with smoother gradients are disregarded. 
Numerical results show that this gives an advantage to the DL model over standard adaptive zoning in accurately reconstructing the trend of the solution close to the shock. Moreover, the DL approach reduces the computational time to perform the adaptive zoning per time step with respect to standard adaptive zoning techniques. 

Future work will extend the approach to multi-dimensional problems, as well as include long short term memory neural networks (LSTM) to mimic the evolutionary nature of the auxiliary parabolic MMPDE. 

\section*{Acknowledgements}
Massimiliano Lupo Pasini thanks Dr. Vladimir Protopopescu for his valuable feedback in the preparation of this manuscript.

This work was supported in part by the Office of Science of the Department of Energy and by the Laboratory Directed Research and Development (LDRD) Program of Oak Ridge National Laboratory. 

\FloatBarrier

\bibliographystyle{unsrt}
\bibliography{bibliography}

\begin{thebibliography}{10}

\bibitem{brackbill82}
J.~U. Brackbill and J.~S. Saltzman.
\newblock Adaptive zoning for singular problems in two dimensions.
\newblock {\em Journal of Computational Physics}, 46:342--368, September 1982.

\bibitem{brackbill93}
J.~U. Brackbill.
\newblock An adaptive grid with directional control.
\newblock {\em Journal of Computational Physics}, 108:38--50, August 1993.

\bibitem{ceniceros_efficient_2001}
H.~D. Ceniceros and T.~Y. Hou.
\newblock An efficient dynamically adaptive mesh for potentially singular
  solutions.
\newblock {\em Journal of Computational Physics}, 172(2):609--639, September
  2001.

\bibitem{lapenta2003}
G.~Lapenta.
\newblock Variational grid adaptation based on the minimization of local
  truncation error: time-independent problems.
\newblock {\em Journal of Computational Physics}, 193:159--179, August 2003.

\bibitem{budd_adaptivity_2009}
C.~J. Budd, W.~Huang, and R.~D. Russell.
\newblock Adaptivity with moving grids.
\newblock {\em Acta Numerica}, 18:111--241, May 2009.

\bibitem{huang2010}
W.~Huang and R.~D. Russell.
\newblock {\em {Adaptive Moving Mesh Methods}}.
\newblock Springer-Verlag Berlin Heidelberg, Berlin, Heidelberg, 2010.

\bibitem{huang_variational_2011}
Weizhang Huang and Robert~D. Russell.
\newblock Variational mesh adaptation methods.
\newblock In {\em Adaptive {Moving} {Mesh} {Methods}}, volume 174, pages
  281--378. Springer New York, New York, NY, 2011.
\newblock Series Title: Applied Mathematical Sciences.

\bibitem{delzanno2011}
G.~L. Delzanno and J.~M. Finn.
\newblock The fluid dynamic approach to equidistribution methods for grid
  adaptation.
\newblock {\em Computer Physics Communications}, 182:330--346, October 2011.

\bibitem{browne2014}
P.~A. Browne, C.~J. Budd, C.~Piccolo, and M.~Cullen.
\newblock Fast three dimensional r-adaptive mesh redistribution.
\newblock {\em Journal of Computational Physics}, 275:174--196, 2014.

\bibitem{remski_balanced_2014}
J.~Remski, J.~Zhang, and Q.~Du.
\newblock On balanced moving mesh methods.
\newblock {\em Journal of Computational and Applied Mathematics}, 265:255--263,
  August 2014.

\bibitem{hu2015}
F.~Hu, R.~Wang, X.~Chen, and H.~Feng.
\newblock An adaptive mesh method for 1d hyperbolic conservation laws.
\newblock {\em Applied Numerical Mathematics}, 91:11--25, 2015.

\bibitem{huang_comparative_2015}
W.~Huang, L.~Kamenski, and R.~D. Russell.
\newblock A comparative numerical study of meshing functionals for variational
  mesh adaptation.
\newblock {\em Journal of Mathematical Study}, 48(2):168--186, June 2015.
\newblock arXiv: 1503.04709.

\bibitem{budd2015}
C.~J. Budd, R.~D. Russell, and E.~Walsh.
\newblock The geometry of r-adaptive meshes generated using optimal transport
  methods.
\newblock {\em Journal of Computational Physics}, 265:113--137, 2015.

\bibitem{pathak2016}
H.~S. Pathak and R.~K. Shukla.
\newblock Adaptive finite-volume {WENO} schemes on dynamically redistributed
  grids for compressible euler equations.
\newblock {\em Journal of Computational Physics}, 319:200--230, May 2016.

\bibitem{heaton_ian_2018}
J.~Heaton, I.~Goodfellow, Y.~Bengio, and A.~Courville.
\newblock Deep learning.
\newblock {\em Genetic Programming and Evolvable Machines}, 19(1-2):305--307,
  June 2018.

\bibitem{fritzen_--fly_2019}
F.~Fritzen, M.~Fernández, and F.~Larsson.
\newblock On-the-fly adaptivity for nonlinear twoscale simulations using
  artificial neural networks and reduced order modeling.
\newblock {\em Frontiers in Materials}, 6:75, May 2019.

\bibitem{yang2012}
X.~Yang, W.~Huang, and J.~Qiu.
\newblock A moving mesh {WENO} method for one-dimensional conservation laws.
\newblock {\em Computer Physics Communications}, 34(4):A2317--A2343, 2012.

\bibitem{CrnjaricZic2007EfficientIO}
N.~Crnjaric-Zic, S.~Macesic, and Bojan Crnkovic.
\newblock Efficient implementation of {WENO} schemes to nonuniform meshes.
\newblock {\em Annali dell'Universita' di Ferrara}, 53:199--215, 2007.

\bibitem{CAPDEVILLE20082977}
G.~Capdeville.
\newblock A central {WENO} scheme for solving hyperbolic conservation laws on
  non-uniform meshes.
\newblock {\em Journal of Computational Physics}, 227(5):2977--3014, 2008.

\bibitem{CORALIC201495}
Vedran Coralic and Tim Colonius.
\newblock Finite-volume {WENO} scheme for viscous compressible multicomponent
  flows.
\newblock {\em Journal of Computational Physics}, 274:95--121, 2014.

\bibitem{huang2018}
W.~F. Huang, Y.~X. Ren, and X.~Jiang.
\newblock A simple algorithm to improve the performance of the {WENO} scheme on
  non-uniform grids.
\newblock {\em Acta Mechanica Sinica}, 34(1):37--47, October 2018.

\bibitem{cao_study_1999}
W.~Cao, W.~Huang, and R.~D. Russell.
\newblock A study of monitor functions for two-dimensional adaptive mesh
  generation.
\newblock {\em SIAM Journal on Scientific Computing}, 20(6):1978--1994, January
  1999.

\bibitem{cao_approaches_2003}
W.~Cao, W.~Huang, and R.~D. Russell.
\newblock Approaches for generating moving adaptive meshes: location versus
  velocity.
\newblock {\em Applied Numerical Mathematics}, 47(2):121--138, November 2003.

\bibitem{chacon2006}
L.~Chac\'on and G.~Lapenta.
\newblock A fully implicit, nonlinear adaptive grid strategy.
\newblock {\em Journal of Computational Physics}, 212:703--717, June 2006.

\bibitem{lapenta2006}
G.~Lapenta and L.~Chac\'on.
\newblock Cost-effectiveness of fully implicit moving mesh adaptation: a
  practical investigation in {1D}.
\newblock {\em Journal of Computational Physics}, 219:86--103, August 2006.

\bibitem{chacon2011}
L.~Chac\'on, G.~L. Delzanno, and J.~M. Finn.
\newblock Robust, multidimensional mesh-motion based on {M}onge–{K}antorovich
  equidistribution.
\newblock {\em Journal of Computational Physics}, 230:87--103, September 2011.

\bibitem{sod1978}
G.~A. Sod.
\newblock A survey of several finite difference methods for systems of
  nonlinear hyperbolic conservation laws.
\newblock {\em Journal of Computational Physics}, 27:1--31, November 1978.

\bibitem{sedov1946}
L.~I. Sedov.
\newblock Propagation of strong shock waves.
\newblock {\em Journal of Applied Mathematics and Mechanics}, 10:241--250,
  November 1946.

\bibitem{woodward}
The numerical simulation of two-dimensional fluid flow with strong shocks.
\newblock {\em Journal of Computational Physics}, 54:115--173, 1984.

\bibitem{toro99}
E.~F. Toro.
\newblock {\em { Riemann Solvers and Numerical Methods for Fluid Dynamics: A
  Practical Introduction}}.
\newblock Springer-Verlag Berlin Heidelberg, Berlin, Heidelberg, 1999.

\bibitem{hirsch}
C.~Hirsch.
\newblock {\em { Numerical Computation of Internal and External Flows, Vol. II:
  Computational Methods for Inviscid and Viscous Flows}}.
\newblock Wiley, Vrije Universiteit Brussel, Brussels, Belgium, 1991.

\bibitem{shu}
X.-D. Liu, S.~Osher, and T.~Chan.
\newblock Weighted essentially non-oscillatory schemes.
\newblock {\em Journal of Computational Physics}, pages 200--212, 1994.

\bibitem{runge}
C.~Runge.
\newblock Ueber die numerische aufl\"{o}sung von differentialgleichungen.
\newblock {\em Mathematische Annalen}, 46:167--178, 1895.

\bibitem{kutta}
M.~Kutta.
\newblock Beitrag zur n\"{a}herungsweisen integration totaler
  differentialgleichungen".
\newblock {\em Zeitschrift für Mathematik und Physik}, 46:435--453, 1901.

\bibitem{Savitzky64}
A.~Savitzky and M.~J.~E Golay.
\newblock Smoothing and differentiation of data by simplified least squares
  procedures.
\newblock {\em Analytical Chemistry}, 36:1627--39, 1964.

\bibitem{kecs82}
W.~Kecs.
\newblock {\em {The Convolution Product - And Some Applications}}.
\newblock D.~Reidel Publishing Company, Dordrecht, Holland / Boston, USA /
  London, England, 1982.

\bibitem{Enders2004}
W.~Enders.
\newblock Stationary time-series models.
\newblock {\em Applied Econometric Time Series (Second edition)}, pages
  48--107, 2004.

\bibitem{mlp}
F.~Rosenblatt.
\newblock Principles of neurodynamics: perceptrons and the theory of brain
  mechanisms.
\newblock {\em Nature}, 542(7639):75--79, 2017.

\bibitem{goodfellow}
I.~Goodfellow, Y.~Bengio, and A.~Courville.
\newblock {\em Deep Learning}.
\newblock The MIT Press, Cambridge, Massachusetts, 2016.

\bibitem{resnet}
K.~He, X.~Zhang, S.~Ren, and Jian Sun.
\newblock Deep residual learning for image recognition.
\newblock {\em arXiv:1512.03385 [cs.CV]}, December 2015.

\bibitem{adam}
D.~P. Kingma and J.~Ba.
\newblock Adam: a method for stochastic optimization, 2017.

\bibitem{MATLAB:2019}
MATLAB.
\newblock {\em (R2019b)}.
\newblock The MathWorks Inc., Natick, Massachusetts, 2019.

\bibitem{cfl}
R.~Courant, K.~Friedrichs, and H.~Lewy.
\newblock {On the partial difference equations of mathematical physics. AEC
  Research and Development Report.}
\newblock Technical report, NYO-7689, New York: AEC Computing and Applied
  Mathematics Centre – Courant Institute of Mathematical Sciences, 1928.

\end{thebibliography}

\newpage

\section{Appendix}
\subsection*{Numerical experiments using shock profiles as input data and adapted mesh as output data}

\begin{figure}[H]
    \centering
    \includegraphics[width=0.45\textwidth]{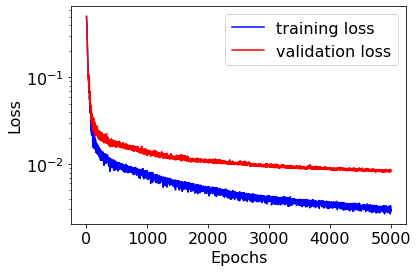}
    \includegraphics[width=0.45\textwidth]{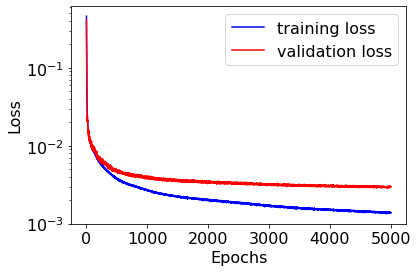}
    \includegraphics[width=0.45\textwidth]{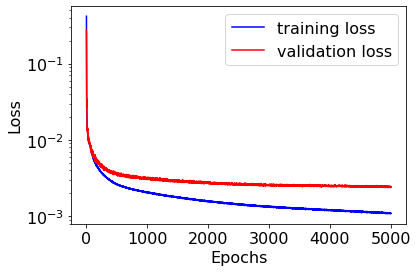}    
    \includegraphics[width=0.45\textwidth]{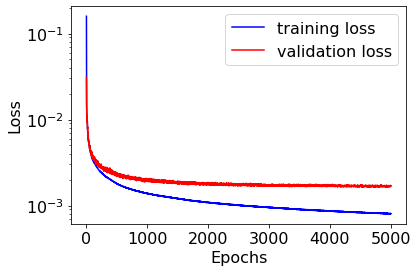}
    \caption{Training input data: shock profile. Training output data: adapted mesh. Training and validation loss function for training of DL models on dataset generated with 500 toy problems (top-left), 5,000 toy problems (top-right), 10,000 toy problems (bottom-left), and 50,000 toy problems (bottom-right).}
    \label{loss_functions_appendix}
\end{figure}




\begin{figure}[H]
    \includegraphics[width=0.32\textwidth]{images/shock1.png}
    \includegraphics[width=0.32\textwidth]{images/shock2.png}
    \includegraphics[width=0.32\textwidth]{images/shock3.png}
    \includegraphics[width=0.32\textwidth]{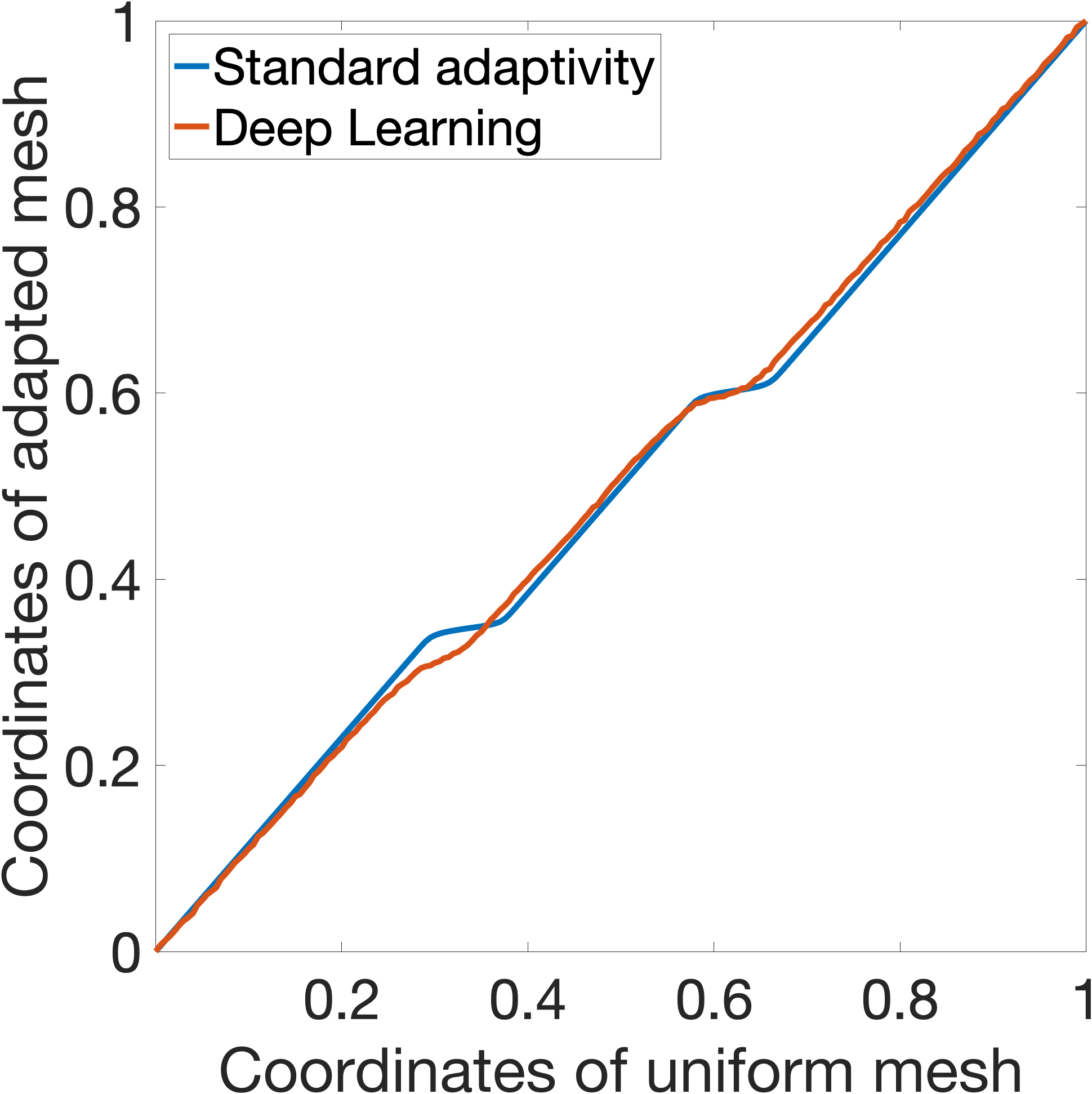}
    \includegraphics[width=0.32\textwidth]{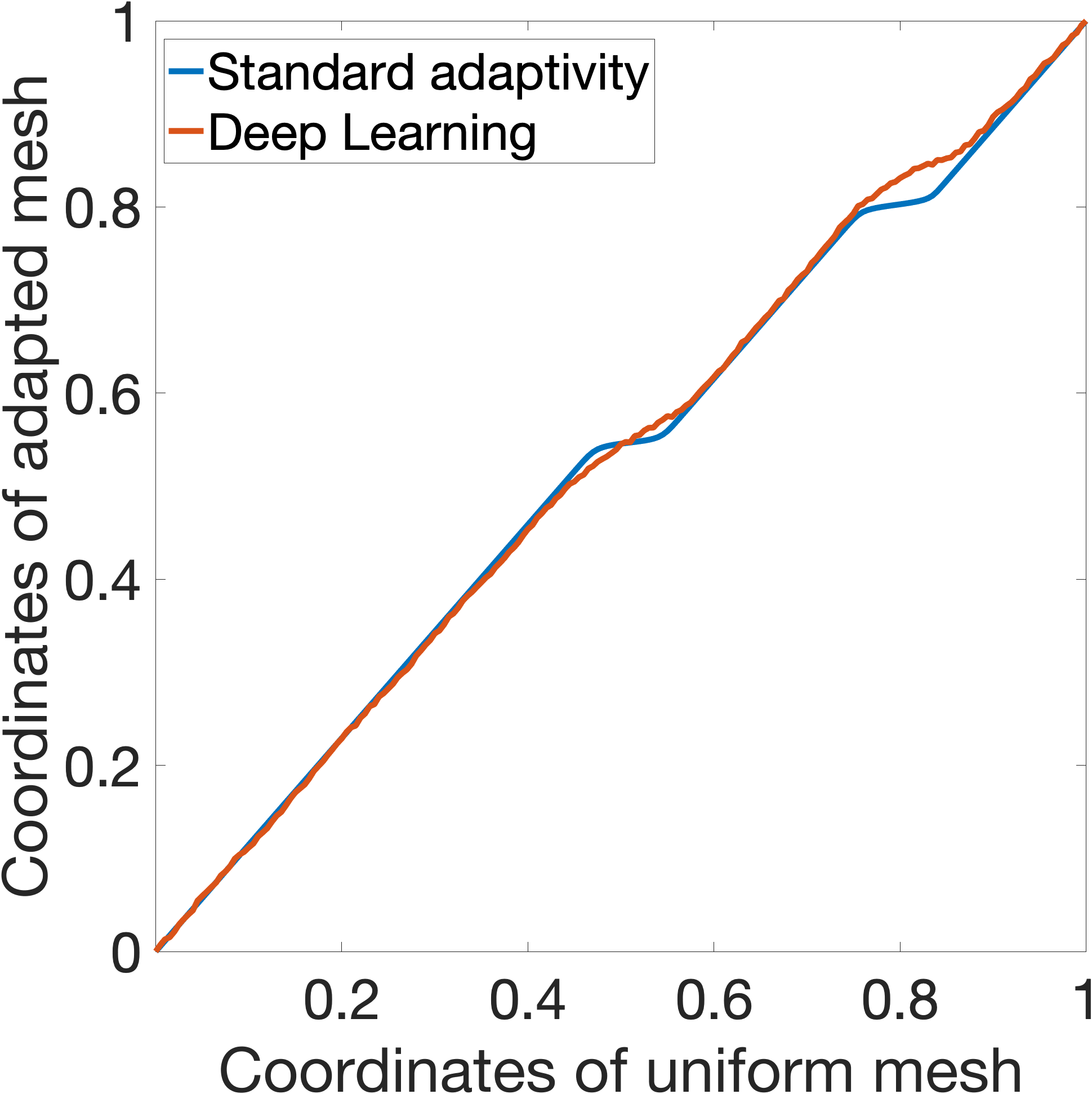}
    \includegraphics[width=0.32\textwidth]{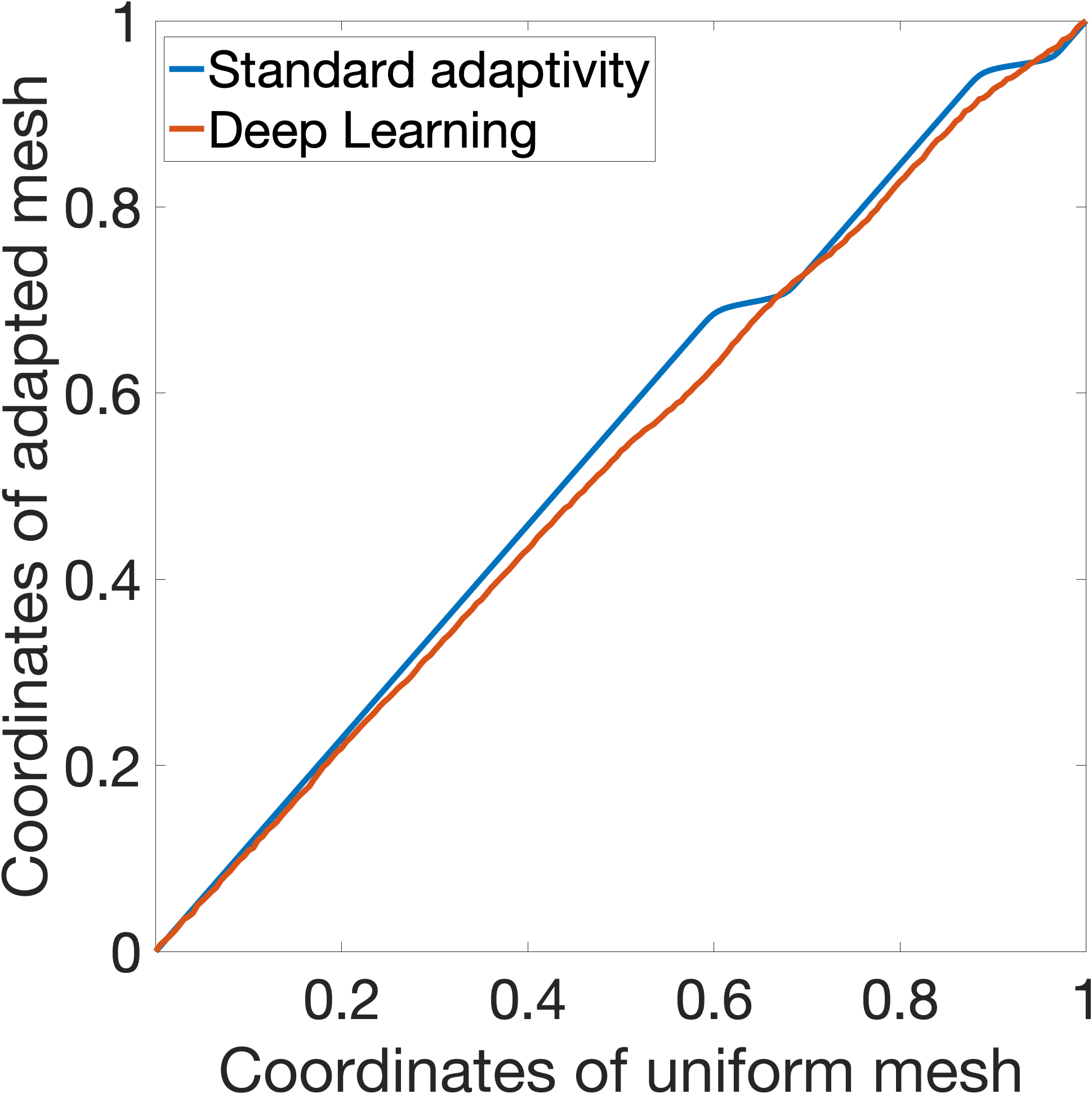}
    \includegraphics[width=0.32\textwidth]{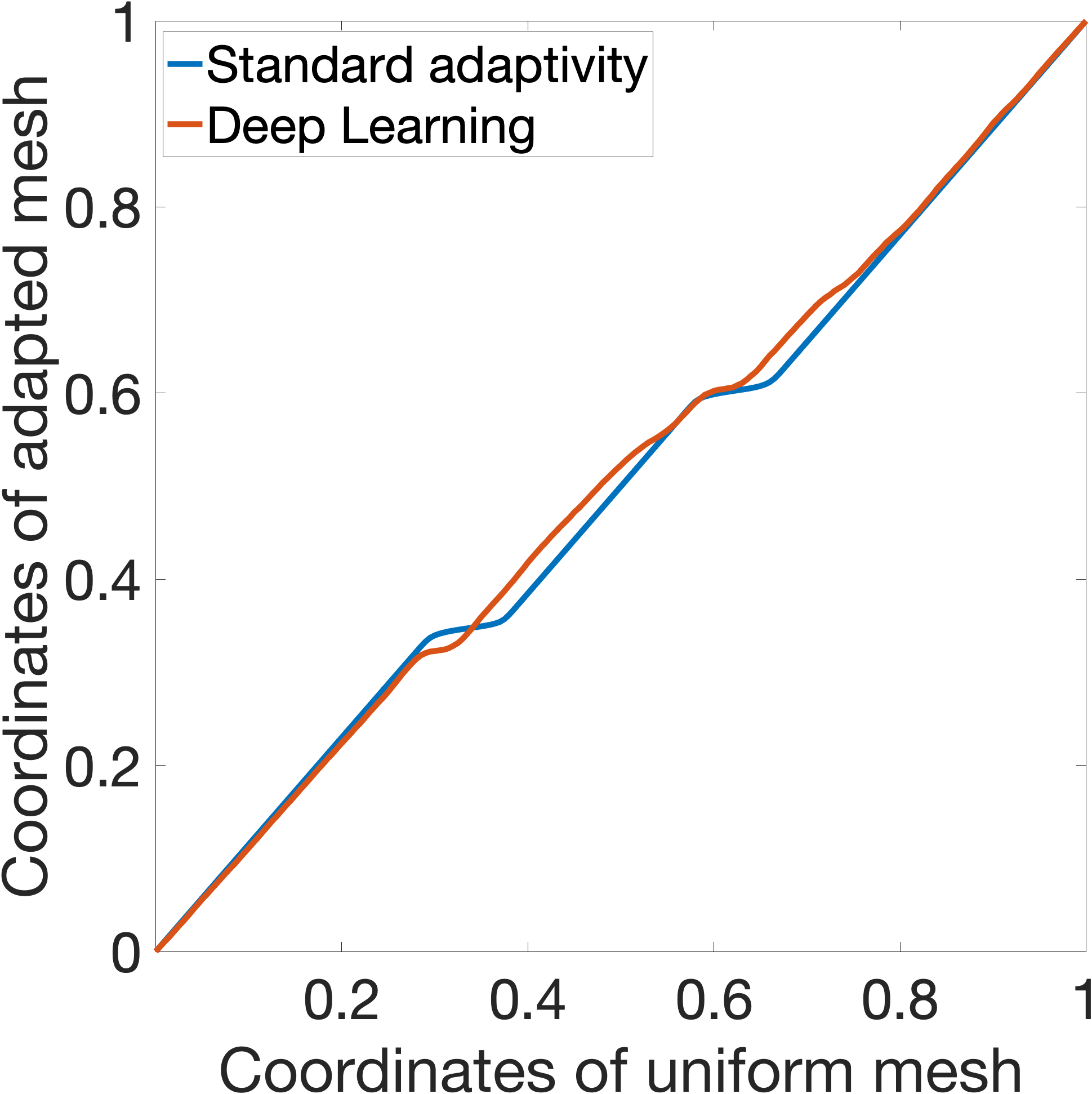}    
    \includegraphics[width=0.32\textwidth]{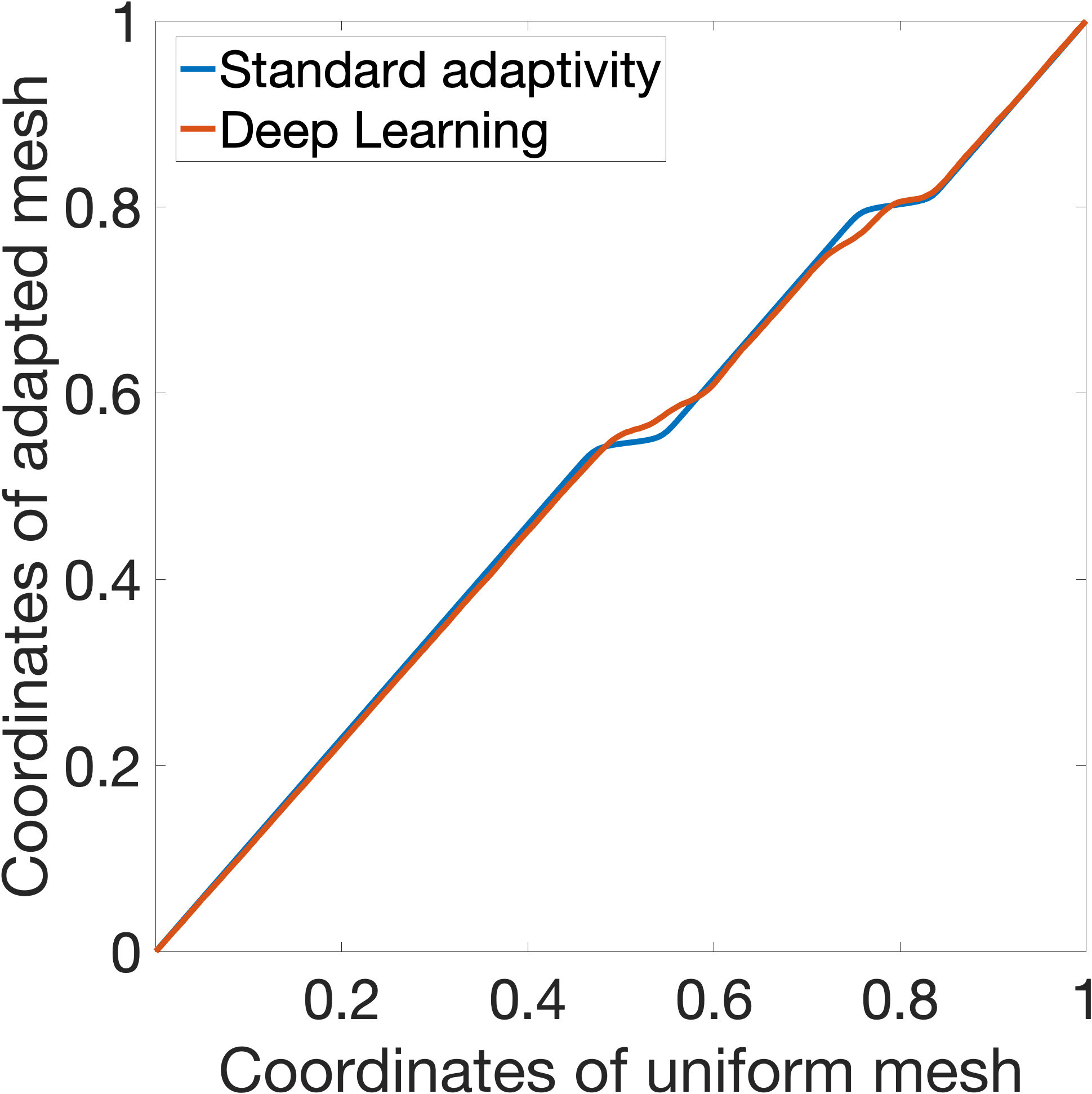}
    \includegraphics[width=0.32\textwidth]{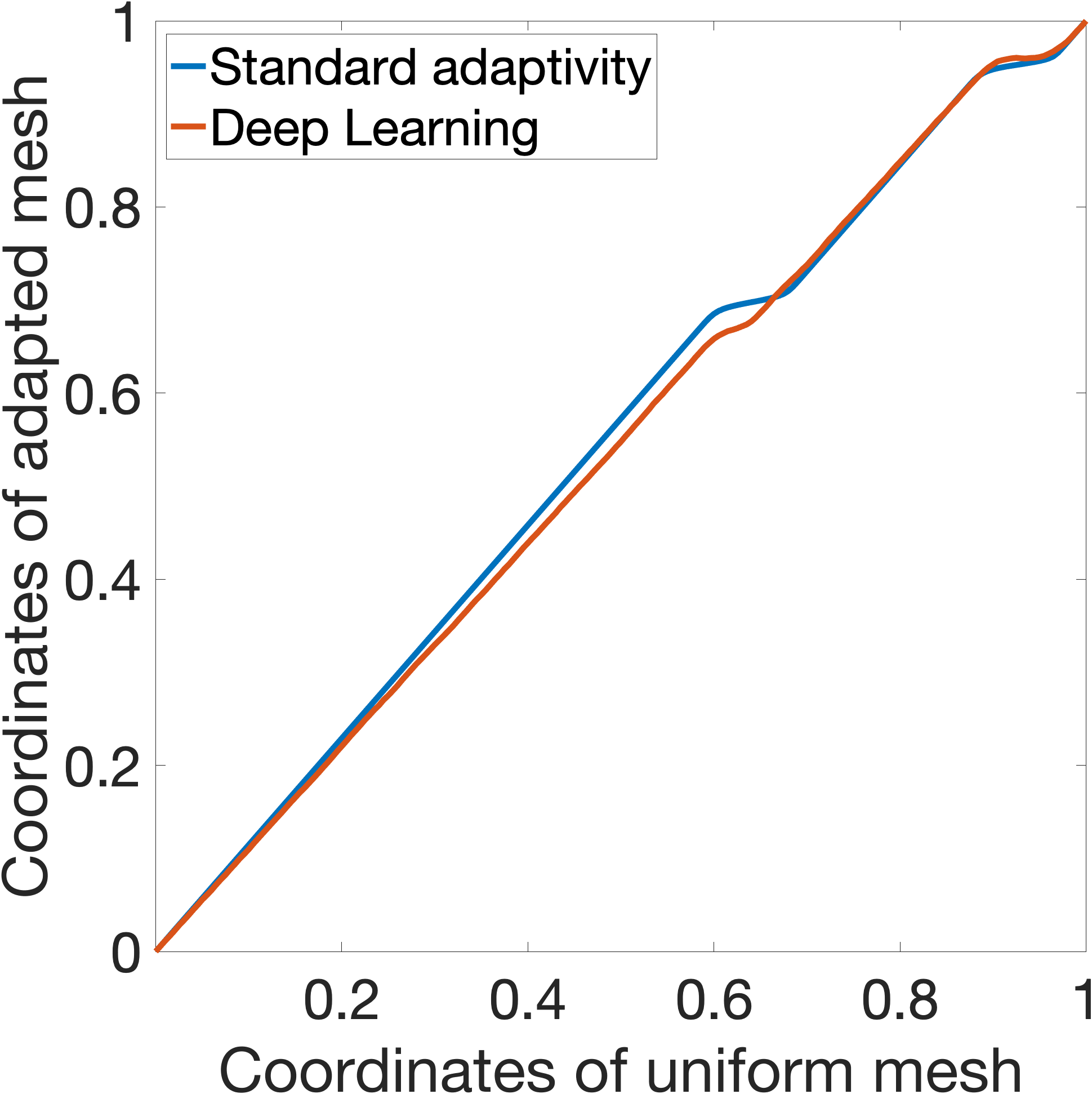}    
    \includegraphics[width=0.32\textwidth]{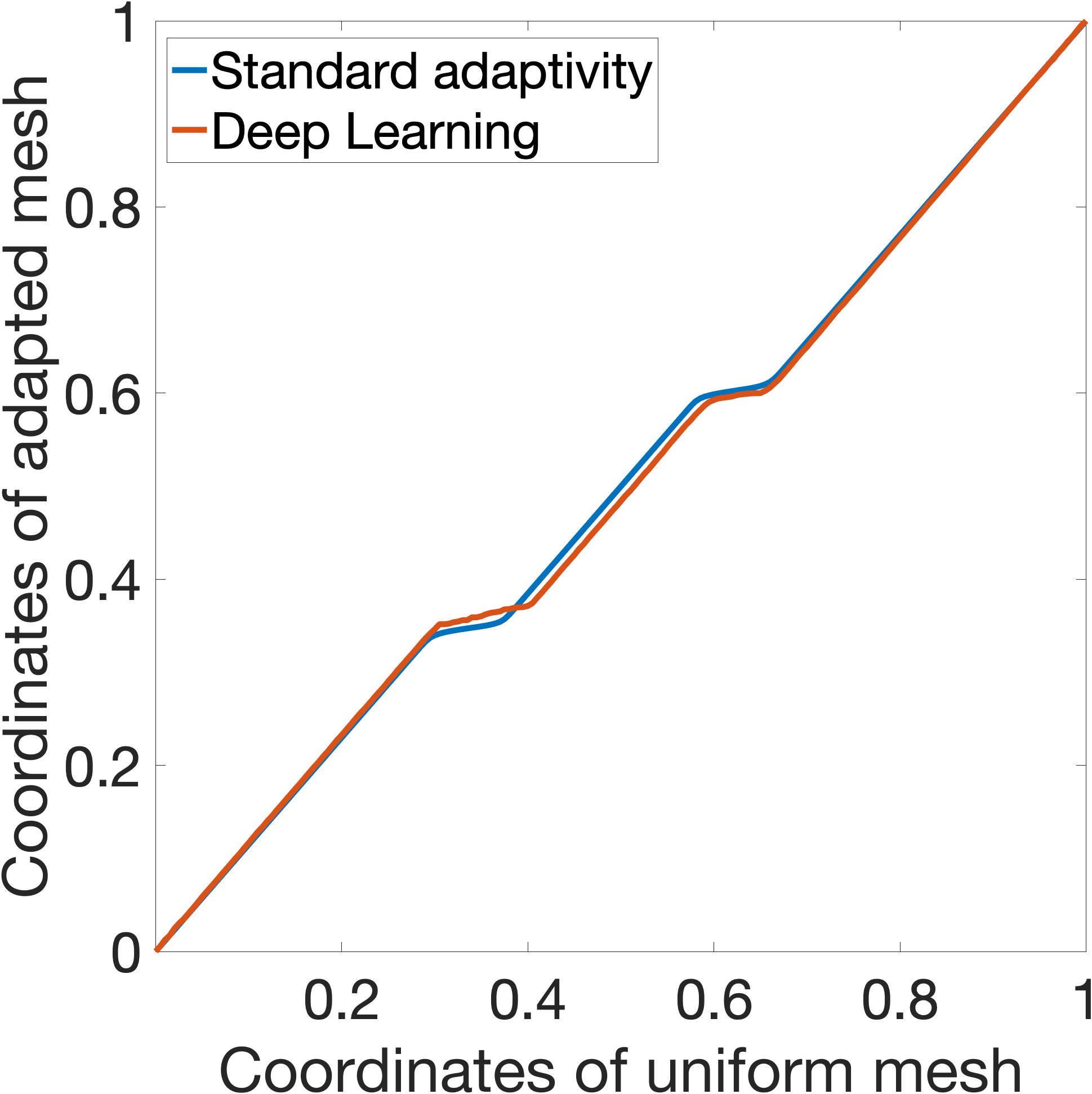}\hspace{0.17cm}
    \includegraphics[width=0.32\textwidth]{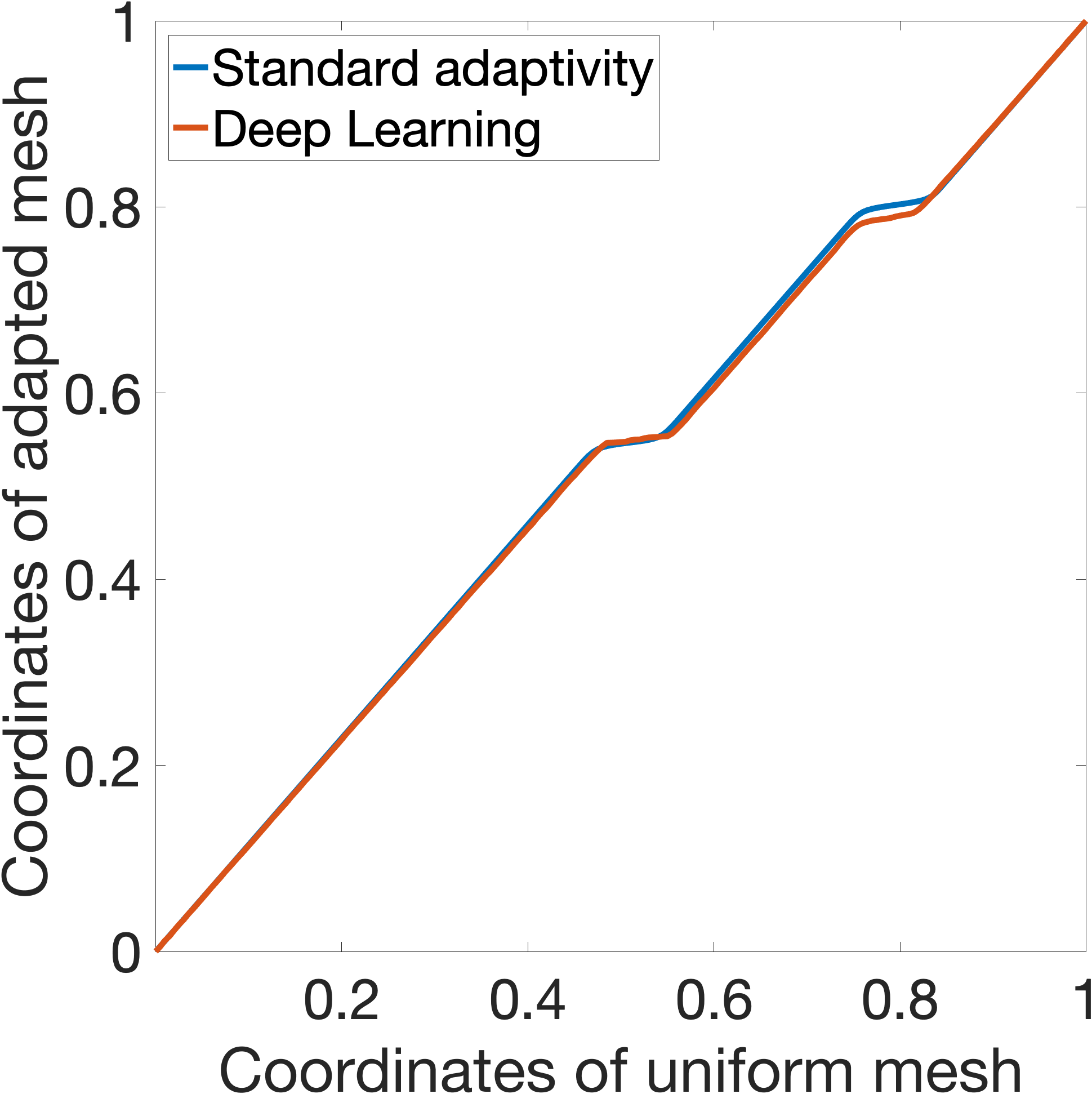}\hspace{0.2cm}
    \includegraphics[width=0.32\textwidth]{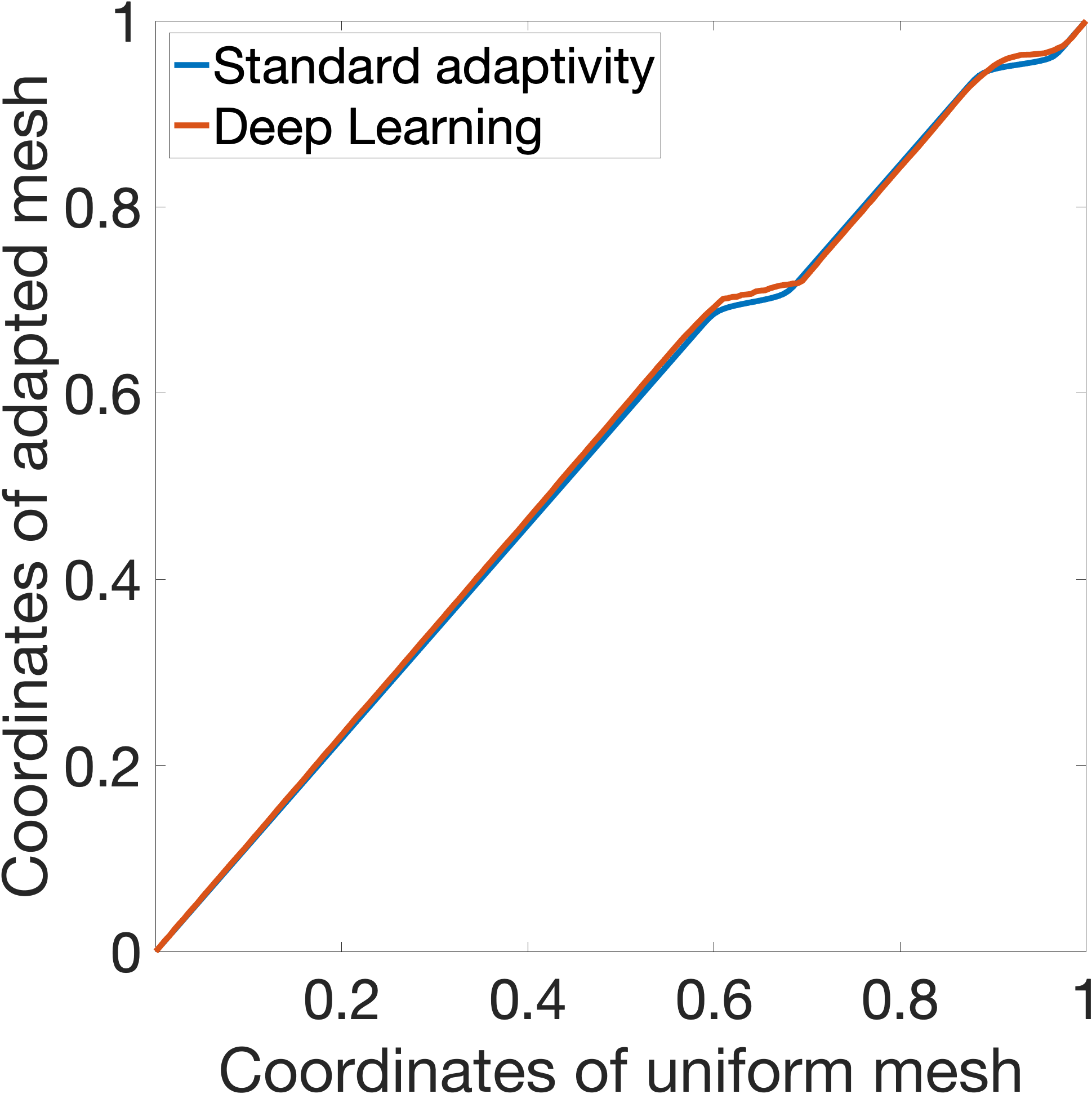}    
    \caption{Training input data: shock profile. Training output data: adapted mesh. DL driven r-adaptivity compared with standard r-adaptivity for propagating shock at time $t=0.05$ (left), $t = 0.25$ (center) and $t=0.4$ (right). The comparison is performed by training the DL model on a dataset with 500, 5,000, and 50,000 toy problems (third from top).}
    \label{mesh_changesize_appendix}
\end{figure}

\end{document}